\documentclass{imsart}[ejs]

\RequirePackage[utf8]{inputenc}
\RequirePackage[OT1]{fontenc}
\RequirePackage{amsthm,amsmath,amssymb,bm}
\RequirePackage[numbers]{natbib}
\RequirePackage[colorlinks,citecolor=blue,urlcolor=blue]{hyperref}

\startlocaldefs
\numberwithin{equation}{section}
\theoremstyle{plain}

\usepackage{cases}
\usepackage[usenames,dvipsnames]{xcolor}

\newtheorem{theorem}{Theorem}
\newtheorem{proposition}{Proposition}
\newtheorem{lemma}{Lemma}

\newtheorem{assumption}{Assumption}
\newtheorem{example}{Example}
\newtheorem{remark}{Remark}

\usepackage{nameref}
\usepackage[nameinlink]{cleveref}
\crefname{equation}{equation}{equations}
\crefname{assumption}{assumption}{assumptions}
\creflabelformat{enumi}{(#2#1#3)}

\usepackage{mathtools}
\usepackage{url}
\usepackage{autonum}

\usepackage{pgfplots}
\pgfplotsset{compat=newest}
\usepgfplotslibrary{groupplots}

\usepackage{booktabs}       % professional-quality tables

\providecommand\given{} % so it exists

\DeclarePairedDelimiterX\Set[1]{\lbrace}{\rbrace}%
{ \renewcommand\given{\SetSymbol[]} #1 }

\global\long\def\dist{\ \sim\ }
\global\long\def\given{\mid}
\global\long\def\distiid{\overset{iid}{\dist}}

\global\long\def\equaldist{\overset{d}{=}}

\newcommand{\EE}{\mathbb E}
\renewcommand{\Pr}{\mathbb{P}}

\global\long\def\bern{\mathrm{Bernoulli}}
\global\long\def\Reals{\mathbb{R}}

\global\long\def\Nats{\mathbb{N}}
\global\long\def\Ints{\mathbb{Z}}
\global\long\def\NNInts{\Ints_{+}}
\global\long\def\NNReals{\Reals_{+}}

\global\long\def\intd{\mathrm{d}} % differential d
\newcommand{\Ind}{\mathbf{1}}

\newcommand{\Ul}{t}
\newcommand{\uft}{x}
\newcommand{\Uft}{x}
\newcommand{\ufSpace}{\mathcal{X}}
\newcommand{\ufMeasure}{\rho}
\newcommand{\size}{t}              % Size of the graph
\newcommand{\Wmarg}{\mu}                % Graphex marginal
\newcommand{\Wcorr}{\nu}                % Graphex 2 points correlation function
\newcommand{\N}{N}                      % Number of vertices
\newcommand{\1}{\mathbf{1}}
\newcommand{\var}{\mathrm{var}}
\newcommand{\cov}{\mathrm{cov}}

\newcommand{\edg}{E}
\newcommand{\ver}{V}
\newcommand{\nedg}{\mathsf{e}}
\newcommand{\nver}{\mathsf{v}}
\newcommand{\degree}{\mathsf{d}}
\newcommand{\indsub}{\mathsf{s}}
\newcommand{\nsub}{\mathsf{n}}
\newcommand{\psamp}{\mathsf{ps}}
\newcommand{\cstat}{h}

\endlocaldefs

\begin{document}

\begin{frontmatter}
  \title{Bootstrap estimators for the tail-index and for the count statistics of
    graphex processes}%
  \runtitle{Bootstrap estimators for graphex processes}%
%  \thankstext{T1}{Footnote to the title with the ``thankstext'' command.}

  %\begin{aug}
\begin{aug}
\author{\fnms{Zacharie} \snm{Naulet}\corref{}\ead[label=e1]{zacharie.naulet@universite-paris-saclay.fr}}
\address{Universit\'{e} Paris-Saclay, Laboratoire de Math\'{e}matiques d'Orsay\\
91405, Orsay, France\\
\printead{e1}}
\author{\fnms{Daniel M} \snm{Roy}\ead[label=e4]{daniel.roy@utoronto.ca}}
\address{Department of Statistics, University of Toronto; Vector Institute\\
University Avenue, Toronto, Ontario, Canada\\
\printead{e4}}
\author{\fnms{Ekansh} \snm{Sharma}\ead[label=e2]{ekansh@cs.toronto.edu}}%
\address{Department of Computer Science, University of Toronto; Vector Institute\\
6 King's College Rd, Toronto M5S 3H3, Ontario, Canada\\
\printead{e2}}
\author{\fnms{Victor} \snm{Veitch}\ead[label=e3]{victorveitch@gmail.com}}%
\address{Department of Statistics, University of Chicago\\ 5747 S Ellis Ave, Chicago, IL 60637\\
\printead{e3}}

\runauthor{Z. Naulet et al.}
\end{aug}

  \begin{abstract}
    Graphex processes resolve some pathologies in traditional random graph
    models, notably, providing models that are both projective and allow
    sparsity. Most of the literature on graphex processes study them from a
    probabilistic point of view. Techniques for inferring the parameter of these
    processes -- the so-called \textit{graphon} -- are still marginal;
    exceptions are a few papers considering parametric families of
    graphons. Nonparametric estimation remains unconsidered. In this paper, we
    propose estimators for a selected choice of functionals of the graphon. Our
    estimators originate from the subsampling theory for graphex processes,
    hence can be seen as a form of bootstrap procedure.
  \end{abstract}



\end{frontmatter}

\section{Introduction}
\label{sec:introduction}

Many statistical models for relational data are based on vertex-exchangeable
random graphs (also called \textit{graphon} models), see \cite{orbanz:roy:2015}
for a review. These models are interpretable and easy to work with, but have the
considerable limitation that, almost surely, they generate dense graphs; that
is, the number of edges scales quadratically with the number of vertices as the
graph grows. Folk wisdom holds that real-world datasets are sparse and their
degree distributions are power-law\footnote{This second point is subject to long
  debates, which we discuss thoroughly in \cref{sec:contr-opin-about}.} implying
that graphon based approaches are unsatisfactory. Several model classes have
been proposed to address this issue. A non exhaustive list includes \textit{edge
  exchangeable} models \cite{crane:dempsey:2016}, \textit{preferential
  attachement} models \cite{barabasi:albert:1999}, \textit{sparse graphon
  models} \cite{bickel:chen:2009,bickel:chen:levina:2011}, and \textit{graphex processes}
\cite{caron:fox:2017,veitch:roy:2015,borgs:chayes:cohn:holden:2016}.

Graphex processes have been proposed as a generalization of the graphon models
that preserve the key properties, but allow for sparsity and degree distribution
with power-law behaviour
\cite{caron:fox:2017,veitch:roy:2015,caron:panero:rousseau:2017}. Thus, graphex
processes are good candidates to model \textit{empirical} graphs. Apart of
\cite{caron:fox:2017,todeschini:miscouridou:caron:2016,veitch:roy:2019,herlau:schmidt:moerup:2016},
inference on graphex processes is still marginal in the literature, most of
papers studying graphex processes from a \textit{probabilistic} point of view
only. Efficient estimation remains challenging for these new models. The present
paper intend to contribute to the study of graphex processes from a
\textit{statistical} viewpoint, driven by two main motivations described in
\cref{sec:motivations}.

On a high-level, a graphex process $\Set{G_t \given t \geq 0}$ is a
\textit{graph-valued} stochastic process. They are used as probabilitic models
for graph observations, which are seen as the realization of a sample $G_t$ of
the process at a finite but unknown \textit{size} $t$. All graphs will be
simple, and will always be considered to be finite, except in rare occasions, in
which case they will be clearly stated to be infinite. The law of the process
$\Set{G_t \given t \geq 0}$ is determined by a jointly measurable symmetric
function $W : \NNReals^2 \to [0,1]$, called the \textit{graphon}
\cite{caron:fox:2017,veitch:roy:2015}. Our goal is to do nonparametric inference
about some functionals of $W$. The functionals we consider are the
\textit{tail-index} of the process (see \cref{sec:estim-tail-index}) and the
\textit{rescaled injective homomorphism densities} (see
\cref{sec:bootstrap-motifs}). The motivations for these functionals are detailed
in \cref{sec:motivations}, and general definitions about
graphex processes are recalled in \cref{sec:undir-graph-proc}. Along the way, we
use our estimators to construct statistical tests for certain counts of
motifs. We propose some applications of our results in
\cref{sec:some-appl-exampl}.

All the estimators we propose are constructed using the subsampling invariance
of graphex processes under the $p$-sampling mechanism uncovered by
\cite{veitch:roy:2019}. In particular, we use the scale invariance of certain
statistics of $G_t$ to estimate the functionals of interest using $p$-samples of
$G_t$. For this reason, our estimators may be viewed as some form of
\textit{bootstrap} estimators
\cite{efron:tibshirani:1994,bhattacharyya:bickel:2015}.

\section{Motivations}
\label{sec:motivations}

\subsection{Sparsity and power-law}
\label{sec:motiv-1:-spars}

Since graphex processes intend to improve on graphon processes in term of
sparsity and power-law behaviour, these notions need to be clarified. Sparsity,
is usually understood as a property of a sequence of graphs, \textit{i.e.} as
the growth of the number of edges $\nedg(G_t)$ in term of the number of vertices
$\nver(G_t)$. In many data analysis situations we just have a single observation
at some particular size. This begs the question: does sparsity have practical
relevance in this case? Regarding the power-law behaviour of the degree
distribution, there is a long historical debate about whether or not this is
desirable to model \textit{empirical} graphs. We come back to this controversy
in \cref{sec:contr-opin-about}. Anyway, it is fair to ask if it makes sense to
account for sparsity and power-law in modeling. The main issue is that we often
think about sparsity and power-law as a property of the sample $G_t$ rather than
a property of the \textit{population}. This is misleading since without modeling
assumptions there is no reason that sparsity and degree distribution of $G_t$
inform us about the population.

Regarding graphex processes, both sparsity and power-law behaviour arise as a
structural property of $W$, \textit{i.e.} the tail-index $\sigma \in [0,1]$, as
explained in \cref{sec:graph-marg-tail}. Consequently, $\sigma$ a distinguished
feature of the population, constituting an interesting summarizing statistic, as it
explains the growth of the network and part of its structure. Whence, rather
than talking about sparsity and power-law, we shall rather talk about the
tail-index. We show in \cref{sec:estim-tail-index} that sparse graphex samples
have a readily identifiable signature which can be uncovered using the
subsampling invariance described in \cref{sec:subs-graph-proc}. This can be
leveraged to visually diagnose if a given graph can be reasonably assumed to be
sample of a sparse graphex process (see \cref{sec:assess-pert-exch}).

Also, we note that it is of practical interest to construct graphex models that
are close analogues of graphon models that have been proven to be useful in
practice; this is the
approach of \cite{todeschini:miscouridou:caron:2016,herlau:schmidt:moerup:2016}, and a number
of models of this kind are described in
\cite{borgs:chayes:cohn:holden:2016}. Then, one wants to combine the efficient
procedures already established for the graphon model with some new procedure for
estimating the additional parameters governing the sparse behaviour, such as the
tail-index.

\subsection{Nonparametric inference, bootstrap, and hypothesis testing}
\label{sec:motiv-2:-nonp}

Apart of the work of
\cite{todeschini:miscouridou:caron:2016,caron:2012,caron:fox:2017,herlau:schmidt:moerup:2016},
algorithms form statistical inference of graphex processes are still marginal in
the literature. Moreover, even the most advanced algorithm of
\cite{todeschini:miscouridou:caron:2016} assumes a parametric structure for
$W$. Full nonparametric estimation of $W$ is to the best of our knowledge still
a challenge. One possibility is to adapt the \textit{method of moments}
\cite{bickel:chen:levina:2011} to graphex processes. The main idea is that $W$
is entirely determined by the collection of all rescaled injective homomorphism
densities
\cite{bickel:chen:levina:2011,bhattacharyya:bickel:2015,diaconis:janson:2008},
which can be further estimated using the \textit{count statistics},
\textit{i.e.} the number of occurrence of certain motifs in the graph. We
elaborate on this in \cref{sec:bootstrap-motifs}.

Besides interest for nonparametric estimation of $W$, the count statistics also
give nice hints on the structure of the graph. They have been used in testing
equality of features of networks and finding confidence intervals of the count
feature \citep{middendorf:ziv:wiggins:2005,shen-orr:milo:mangan:2002}, and
similarly to determine if the dearth of certain motifs is statistically
significant \cite{bhattacharyya:bickel:2015}.

While in principle count statistics are exactly computable, in practice the
problem becomes rapidly intractable for large motifs. In
\cite{bhattacharyya:bickel:2015}, the authors propose to use two different types
of bootstrap techniques to make computations feasible. The idea is to count the
motifs in a possibly small subsampled subgraph to estimate the count in the
whole graph. Here we propose an alternative using $p$-sampling. Though for
counting motifs the performance is equivalent (and thus for inference too),
there are some advantages when it comes to test some hypotheses, or simply if
one is interested in also estimating the variance of the count statistics.

As in \cite{bhattacharyya:bickel:2015}, we apply our method to construct
statistical tests about the count statistics. However, in contrast with
\cite{bhattacharyya:bickel:2015} our null hypotheses is a graphex process,
instead of sparse graphon process. Though it might seem anecdotal\footnote{In
  particular these two models are tightly connected, see
  \cite{borgs:chayes:cohn:holden:2016}.}, this has a major advantage:
bootstraped samples of sparse graphon processes do not have distinguished
distributional properties, while $p$-sampling preserves the distribution of
graphex processes. We exploit this invariance to construct estimators for the
variances of count statistics that are based on recycling the bootstrap samples
used to estimate them, while \cite{bhattacharyya:bickel:2015} have to rely on
the construction of ad-hoc estimators. Our estimators are simply rescaled
version of the empirical variance of the bootstrap samples, and thus
straightforward to compute.

\section{Simple graphex processes}
\label{sec:undir-graph-proc}

\subsection{Notations}
\label{sec:notations}

Most of our notations are taken or inspired from
\cite{diaconis:janson:2008}. We denote the vertex and edge sets of a graph $G$
by $\ver(G)$ and $\edg(G)$, and the numbers of vertices and edges by $\nver(G)
\coloneqq |\ver(G)|$ and $\nedg(G) \coloneqq |\edg(G)|$. Also for every
$j \in \Nats \coloneqq \Set{1,2,\dots}$, we let $\degree_j(G)$ be the number of
vertices with degree $j$ in $G$. We consider both labeled and unlabeled graphs;
the labels will always be the integers $1,\dots,n$ where $n$ is the number of
vertices in the graph. A labeled graph is thus a graph with vertex set
$\Nats_n \coloneqq \Set{1,\dots,n}$ for some $n\geq 1$; we let $\mathcal{L}_n$
denote the set of the $2^{\binom{n}{2}}$ labeled graphs on $\Nats_n$ and let
$\mathcal{L} \coloneqq \bigcup_{n=1}^{\infty}\mathcal{L}_n$. An unlabeled graph
can be regarded as a labeled graph where we ignore the labels; formally, we
define $\mathcal{U}_n$, the set of unlabeled graphs or order $n$, as the
quotient $\mathcal{L}_n / \cong$ of labeled graphs modulo isomorphisms. We let
$\mathcal{U} \coloneqq \bigcup_{n=1}^{\infty}\mathcal{U}_n = \mathcal{L} /
\cong$, the set of all unlabeled graphs.

We write $X_t \sim Y_t$ for $\lim_t X_t/Y_t = 1$, $X_t = O(Y_t)$ for
$\limsup_t|X_t/Y_t| < \infty$, and $X_t = o(Y_t)$ for $\limsup_t|X_t/Y_t| =
0$. We also write $X_t \asymp Y_t$ if there is a constant $c \in \Reals$ such
that $\lim_tX_t/Y_t = c$. If $X_t$ and $Y_t$ are random variables, the limits
are understood almost-surely. In case where the limits hold in probability, we
use $X_t = O_p(Y_t)$ and $X_t = o_p(Y_t)$.  Inequalities up to generic constants
are denoted by the symbols $\lesssim$ and $\gtrsim$. The symbol $\EE$ denote
expectations, and when it is needed to emphasize that
$\Set{G_t \given t \geq 0}$ is a graphex process with parameter $W$, we use
$\EE_W$.

\subsection{Model}
\label{sec:model}

We now define the (simple) graphex process. The word simple refers to the fact
that we consider undirected graphs with no loops. The model has been introduced
in \cite{caron:fox:2017} and further studied in
\cite{veitch:roy:2015,veitch:roy:2019,todeschini:miscouridou:caron:2016,borgs:chayes:cohn:holden:2016,borgs:chayes:cohn:veitch:2017,janson:2017:2,caron:panero:rousseau:2017}. We
recall that $W : \NNReals^2 \to [0,1]$ is a jointly measurable symmetric
function, usually referred to as the \textit{graphon} \cite{veitch:roy:2015}.

A graphex process $\Set{G_t \given t \geq 0}$ is a \textit{graph-valued}
stochastic process which can be constructed as follows. Sample a unit-rate
Poisson random measure $\mathcal{V}$ on $\NNReals^2$
\citep[see][]{daley:vere-jones:2003,daley:vere-jones:2007}, and identify
$\mathcal{V}$ with the collection of points $(t,x) \in \NNReals^2$ such that
$\mathcal{V}$ has a point mass, that is
$\mathcal{V} \equiv \Set{(t_1,x_1),(t_2,x_2),\dots}$. Let
$\tilde{G} \in \mathcal{L}$ be a (eventually infinite) graph with vertex set
$\mathcal{V}$, such that for each pair of vertices $v_i = (t_i,x_i)$ and
$v_j = (t_j,x_j)$, there is an edge between $v_i$ and $v_j$ with probability
$W(x_i,x_j)$, independently for any two $i\ne j$. Let $G\in \mathcal{L}$ be the
subgraph of $\tilde{G}$ consisting only of the vertices with degree at least
one. Note that $G$ is almost-surely a graph with countably infinitely many
vertices, and that the set of edges is also countably infinite except if $W$ is
equal to $0$ almost-everywhere. For each $t \geq 0$ let $G_t\in \mathcal{U}$ be
the unlabeled induced subgraph of $G$ consisting only on the vertices $(t',x)$
such that $t' \leq t$. The process $\Set{G_t \given t \geq 0}$ is a graphex
process with parameter $W$. For any fixed $t \geq 0$, $G_t \in \mathcal{U}$ is
called a sample of size $t$. We insist that $G_t$ is unlabeled. In other words,
the labels $(t_i,x_i)$ are not observed, and the statistician may think of them
as latent variables.

If
$\int_{\NNReals^2} W(x,y)\,\intd x\intd y <
\infty$, then at each $t \geq 0$ the sample $G_t$ has finitely many vertices and
edges almost-surely \cite{caron:fox:2017,veitch:roy:2015}.

\begin{assumption}
  \label{ass:finiteness}
  The graphon satisfies
  $0 < \int_{\NNReals^2} W(x,y)\, \intd x \intd y < \infty$.
\end{assumption}

\begin{remark}
  \label{rmk:3}
  In this paper, we only consider graphex processes with no loops. We note that
  in the seminal papers \cite{caron:fox:2017,veitch:roy:2015} graphex processes
  are defined with the possibility of having loops, and hence the definition
  here is a bit less general. Indeed, the processes in this paper are a special
  case of \cite{caron:fox:2017,veitch:roy:2015} corresponding to let
  $W(x,x) = 0$ for all $x \in \NNReals$ in their model.
\end{remark}

\subsection{Subsampling of graphex processes}
\label{sec:subs-graph-proc}

Various sampling designs have been proposed in the statistical and computer
science literature to derive representative samples of a given network
\cite{kolaczyk:2009,thompson:2012,leskovec:kleinberg:faloutsos:2005,orbanz:roy:2015}. Many
of these sampling designs have been analyzed from the design-based sampling
point of view \cite{thompson:frank:2000,carrington:scott:wasserman:2005}. It is
common in the statistical network literature to oppose model-based and
design-based sampling. In the context of graphex processes,
\cite{veitch:roy:2019,borgs:chayes:cohn:veitch:2017} make the bridge between the
two points of view by identifying the sampling design naturally associated with
graphex processes.

Of particular interest, \cite{veitch:roy:2019,borgs:chayes:cohn:veitch:2017}
defines the \textit{$p$-sampling} procedure. Given a graph $G$, a $p$-sample of
$G$, written later $\psamp(G,p)$, is the random graph obtained from $G$ by
sampling the vertices of $G$ independently with probability $p$ and returning
the induced subgraph where only non-isolated vertices are conserved. The
following lemma is taken from \cite{veitch:roy:2019} and will be at the core of
our inference technique.

\begin{lemma}
  \label{lem:1}
  Let $\Set{G_t \given t \geq 0}$ be a graphex process and $p \in [0,1]$. Then
  $\psamp(G_t,p) \equaldist G_{pt}$ for any $t \geq 0$.
\end{lemma}

From \cite{borgs:chayes:cohn:veitch:2017} we also know the converse result that
if $\psamp(G_t,p) \equaldist G_{pt}$, then $\Set{G_t \given t \geq 0}$ is a
graphex process. This result shed light on the sampling design associated with
graphex processes: we may also view $G_t$ as a sample obtained from a
\textit{population} graph $G_{t'}$ with $t' > t > 0$ and
$G_t = \psamp(G_{t'},t/t')$.

\subsection{Graphex marginal and tail-index}
\label{sec:graph-marg-tail}

Recently, \cite{caron:panero:rousseau:2017} have studied asymptotics for an
important class of graphex models (encompassing most known examples). They
derive the remarkable result that, for this model class, the sparsity and the
degree distribution of the graph are governed by the tail behaviour of the
graphex marginal $\Wmarg : \NNReals \rightarrow \NNReals$ defined by
\begin{equation}
  \Wmarg(\uft) \coloneqq \int_{\NNReals}W(\uft,\uft')\,\intd\uft'.
\end{equation}
More specifically, assuming without loss of generality that $\Wmarg$ is
monotone, they show that if there is $\sigma \in [0,1]$ and a slowly varying
function\footnote{We recall that a function
  $\ell : \NNReals \rightarrow \NNReals$ is called slowly varying if
  $\lim_{x\rightarrow \infty} \ell(cx)/\ell(x) = 1$ for all $c > 0$.} $\ell$
such that as $z \rightarrow 0$,
\begin{equation}
  \label{eq:v1:14}
  \mu^{-1}(z) \coloneqq \inf\Set{y > 0 \given \mu(y) \leq z}%
  \sim z^{-\sigma}\ell(z^{-1}),
\end{equation}
then, under mild supplementary assumptions on $W$ (related to our
\cref{ass:v1:5} below), the asymptotic behaviour of $\nver(G_t)$, $\nedg(G_t)$
and $\degree_j(G_t)$ is determined by $\sigma$. They characterize the four
regimes listed below.
\begin{itemize}
  \item   \textit{Dense} : $\sigma = 0$ and
  $\lim_{t\rightarrow \infty} \ell(t) < \infty$. In that case $\mu^{-1}$ has compact
  support and $\nedg(G_t) \asymp \nver(G_t)^2$ almost-surely.

  \item   \textit{Almost-dense} : $\sigma = 0$ and
  $\lim_{t\rightarrow \infty}\ell(t) = \infty$. In that case
  $\nedg(G_t) / \nver(G_t)^2 \rightarrow 0$ almost-surely, but
  $\nedg(G_t) / \nver(G_t)^{2-\epsilon} \rightarrow \infty$ almost-surely for
  every $\epsilon > 0$.
  
  \item \textit{Sparse with power-law} : $\sigma \in (0,1)$. In that case,
  $\nedg(G_t) \asymp [\nver(G_t)/\ell(\size)]^{\frac{2}{1+\sigma}}$
  almost-surely. Moreover, for every $j\in\Nats$, $\degree_j(G_t)/\nver(G_t) \to
  \sigma\Gamma(j-\sigma)/j!$ almost-surely as $t \to \infty$.

  \item \textit{Very sparse} : $\sigma = 1$. In that case, $\ell$ has to go to
  zero sufficiently fast and $\nedg(G_t) / \nver(G_t) \rightarrow \infty$
  almost-surely while $\nedg(G_t) / \nver(G_t)^{1+\epsilon} \rightarrow 0$
  almost-surely for all $\epsilon > 0$.
\end{itemize}
We call $\sigma$ in \cref{eq:v1:14} the tail-index of $\Wmarg$. Here we are
mostly interested in the sparse with power-law regime, \textit{i.e.} $\sigma \in
(0,1)$. It is the most interesting regime, since the dense regime is equivalent
to the classical graphon model, and the very sparse is an extreme case unlikely
to be useful in statistical applications.

\section{Estimating the tail-index of graphex processes}
\label{sec:estim-tail-index}

\subsection{Heuristic and the $p$-sampling argument}
\label{sec:heur-some-diagn}

Before introducing our estimator for $\sigma$ and results on its risk, we review
here the heuristic we used to construct the estimator. Since it relies mostly on
\cref{lem:1}, we shall refer to this heuristic as the \textit{$p$-sampling
  argument} in the next. The idea is to use the scale-invariance of certain
statistics.

From \cite{caron:panero:rousseau:2017}, we have
$\nver(G_t) \sim t^{1+\sigma}\ell(t)\Gamma(1-\sigma)$ almost-surely as
$t \to \infty$, under mild conditions. Imagine we could observe
$\Set{G_t : t \geq 0}$ at different sizes $t < t'$, with $t/t'$ known. Then, we
could asymptotically recover $\sigma$ because
\begin{equation}
  \label{eq:6}
  \frac{\nver(G_{t'})}{\nver(G_t)} \sim \Big( \frac{t'}{t}\Big)^{1+\sigma}
  \frac{\ell(t')}{\ell(t)}%
  \approx \Big( \frac{t'}{t}\Big)^{1+\sigma}
\end{equation}
almost-surely as $t \to \infty$, at least if $t'/t$ remains constant. Of course,
we usually have only one sample $G_t$ with unknown $t > 0$. The key insight is
that $p$-sampling allows us to effectively simulate a sample $G_{t'}$ with
$t' < t$; intuitively, this is because $\psamp(G_t,p) \equaldist G_{pt}$ by
\cref{lem:1}. Of course, here we don't have access to the marginal law of
$\psamp(G_t,p)$, only to the conditional $\psamp(G_t,p) \mid G_t$. Hence some
cares have to be taken, but the intuition remains useful. In particular, we can
asymptotically recover $\sigma$ from $\nver(G_t)$ and $\nver(\psamp(G_t,p))$
for some $p\in [0,1]$, or even better, from $\nver(G_t)$ and%
\begin{equation}
  \label{eq:19}
  N_p(G_t) \coloneqq \EE[\nver(\psamp(G_t,p)) \mid G_t].
\end{equation}
It is easy to see that for any simple graph $G \in \mathcal{U}$ with no isolated vertex,
\begin{equation}
  \label{eq:33}
  N_p(G)%
  =%
  p \sum_{j\geq 1} \degree_j(G)\big(1 - (1-p)^j\big).
\end{equation}

\begin{figure}[h]
  \centering

  \includegraphics[width=\linewidth]{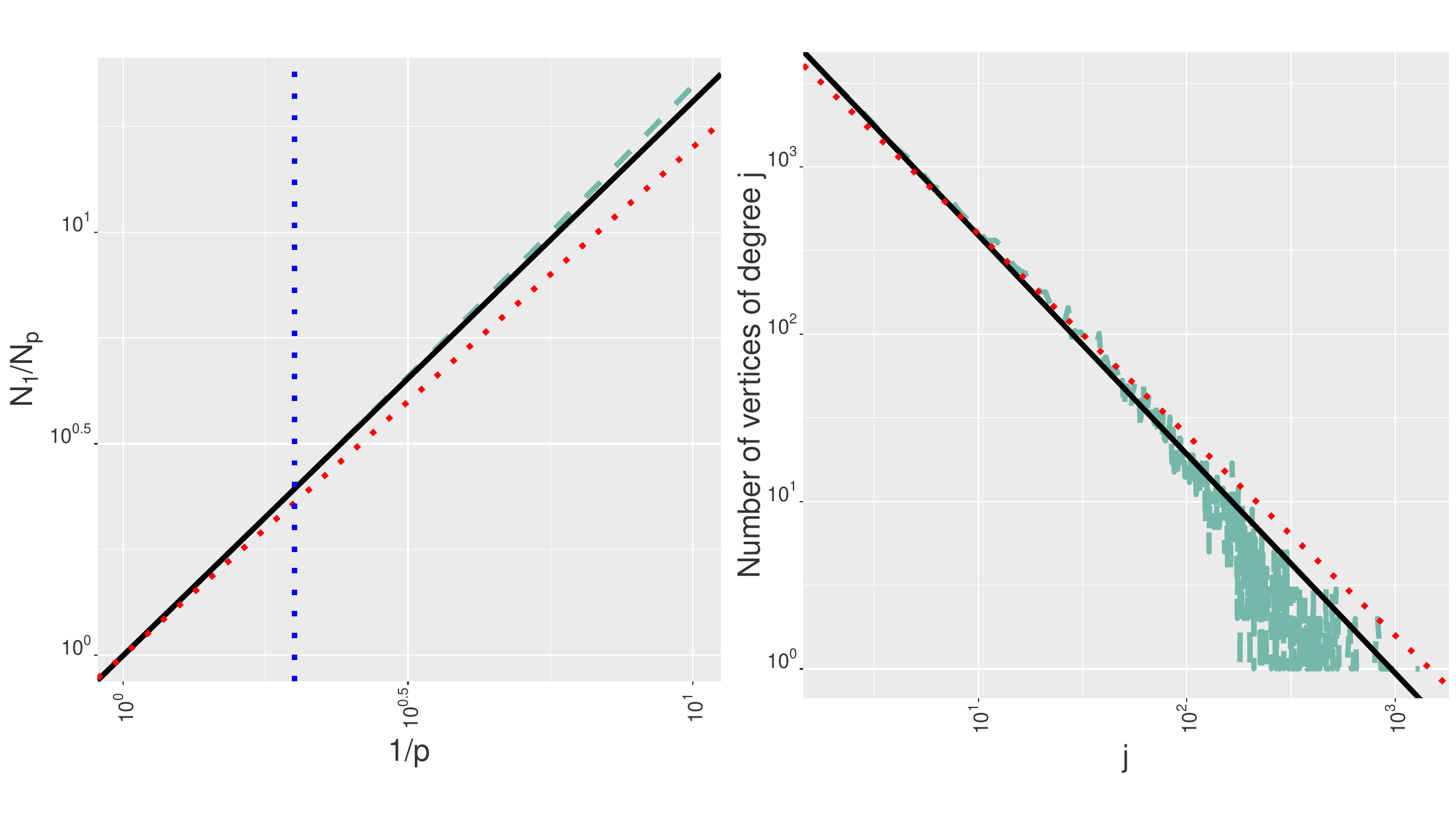}
  
  \caption{%
    On the left: the plot represent the value of $N_1(G_t)/N_p(G_t)$ as a
    function of $p$ on a log-log scale. Here $G_t$ is a sample from the GGP
    model of \cite{caron:fox:2017} with parameters $(\sigma,\tau) = (0.2,5)$ and
    $t = 2000$ (see also \cref{ex:v1:5}). The sample has 26671 vertices and
    293923 edges. With these parameters $\mu$ has a tail-index equal to 0.2. The
    red dotted line is $y = (1+\sigma) x$, which is the expected asymptotic
    behaviour of $\log \frac{N_1(G_t)}{N_p(G_t)}$ as a function of $-\log
    p$. The plain black line is $y = (1+\hat{\sigma}_p) x$, where
    $\hat{\sigma}_p$ is the estimated value of $\sigma$ using $p = 0.5$
    (represented by the blue dotted vertical line). On the right: the green
    curve represents the degree distribution of the same sample. In plain black
    (respectively dotted red) is represented
    $\degree_j \propto j^{-1-\hat{\sigma}_p}$ (respectively
    $\degree_j \propto j^{-1-\sigma}$ to illustrate that the degree distribution
    is also governed by the tail-index.}
  \label{fig:1}
\end{figure}

In \cref{fig:1}, we illustrate on a simulation example how to infer $\sigma$
from $p\mapsto N_1(G_t)/N_p(G_t)$. Indeed, $\log(N_1(G_t)) - \log(N_p(G_t)) \to
(1+\sigma)\log(p)$ almost-surely as $t \to \infty$ under mild assumptions, by
the discussion above and the results in \cite{caron:panero:rousseau:2017}. Whence,
$\sigma$ may be estimated by determining the slope of curve
$\log(N_1(G_t)/N_p(G_t))$ as a function of $p$\footnote{At least the part of the
curve corresponding to $p$ bounded away from $0$, because the asymptotic
equivalences are true as $t \to \infty$ and $p t \to \infty$.}.%

\subsection{An estimator and its risk}
\label{sec:heur-some-diagn-1}

In view of the heuristic of the previous section, we can construct a simple
estimator for $\sigma$. For a chosen value of $p \in (0,1)$, we propose to use
as an estimator for $\sigma$,
\begin{equation}
  \label{eq:v1:11}
  \hat{\sigma}_p(G_t) \coloneqq
  \begin{cases}
    \frac{\log N_1(G_t) - \log
      N_p(G_t)}{-\log p} - 1&\mathrm{if}\ N_p(G_t) \geq 1,\\
    0 & \mathrm{otherwise}.
  \end{cases}
\end{equation}
Notice that $\hat{\sigma}_p$ is well-defined since
$N_1 \geq N_p$ almost-surely for all $0 \leq p \leq 1$. Also
note that $\hat{\sigma}_p$ does not depend explicitly on the sample size
$\size$; this is important because $\size$ is generally considered to be
unknown.

% In practice, we observe better performance by computing the estimator for values
% of $p \in (0,1)$ close to $1$, see \cref{sec:v1:risk-estim-vari} for illustration.

In order to be able to compute the risk of $\hat{\sigma}_p$, we require further
assumptions on the model. Our first assumption is analogous to the assumptions
encountered in \cite{veitch:roy:2015,caron:panero:rousseau:2017}. We first
introduce the $2$-points correlation function
$\nu : \NNReals^2 \rightarrow \NNReals$,
\begin{equation}
  \Wcorr(x,x') \coloneqq \int_{\NNReals} W(x,y)W(y,x')\,\intd y.
\end{equation}

We also make the following assumption on $\Wcorr$, which is a classical
assumption in the literature \cite{veitch:roy:2015,caron:panero:rousseau:2017}
to be able to control the variance of certain (non-linear) functionals of
graphex processes, such has the number of vertices.
\begin{assumption}
  \label{ass:v1:5}
  If $0 \leq \sigma < 1$, we assume that there is a constant
  $\eta > \max(1/2,\sigma)$ and a constant $C \geq 0$ such that
  $\Wcorr(\uft,\uft') \leq C \mu(\uft)^{\eta} \mu(\uft')^{\eta}$ for all
  $\uft,\uft' \in \ufSpace$. If $\sigma = 1$, we assume
  $\nu(x,x') \leq C \mu(x)\mu(x')$, and we set $\eta = 1$.
\end{assumption}

Also, it is well-known \cite{feller:1971,caron:panero:rousseau:2017}, that under
\cref{eq:v1:14} the asymptotic equivalence
\begin{equation}
  \label{eq:v1:7}
  \int_{\NNReals}\left(1 - e^{-\size \mu(z)}\right)\,\intd z
  \sim
  \begin{cases}
    \size^{\sigma}\ell(\size) \Gamma(1 - \sigma) &\mathrm{if}\ 0 \leq \sigma <
    1,\\
    \size \int_{\size}^{\infty} z^{-1}\ell(z)\,\intd z &\mathrm{if}\ \sigma = 1
  \end{cases}
\end{equation}
holds as $\size \rightarrow \infty$.  This relation is already enough to show
consistency of $\hat{\sigma}_p$.  However, further assumptions on $\mu$ are
required to bound the bias of the estimator. In the present paper, we shall
assume the following.
\begin{assumption}
  \label{ass:v1:2}
  We assume that \cref{eq:v1:14} holds. We furthermore assume that for all
  $p \in (0,1)$ there is a sequence $(\Gamma_{p,\size})_{\size > 0}$ such that
  \begin{equation}
    \left| \frac{\int_{\NNReals}\left(1 - e^{-p\size \mu(z)}\right)\, \intd z}
      {p^{\sigma} \int_{\NNReals}\left(1 - e^{-\size \mu(z)}\right)\, \intd z} - 1 \right| \leq \Gamma_{p,\size}.
  \end{equation}
\end{assumption}

\Cref{ass:v1:2} is closely related to second-order regular variation assumptions
used in extreme values theory \citep{haan:stadtmueller:1996}. Its purpose is the
following. The integral in \cref{eq:v1:7} plays a determining role in computing the
risk of $\hat{\sigma}_{p}$. \Cref{eq:v1:7}, however, says nothing about the
rate at which the integral approaches $\size^{\sigma}\ell(\size)$. Moreover, we
know that $\ell(p \size) / \ell(\size) \rightarrow 1$ as
$\size \rightarrow \infty$ by the assumptions on $\ell$, but again, the rate of
convergence is unknown. These two rates are encapsulated in \cref{ass:v1:2},
allowing us to quantify the risk of $\hat{\sigma}_p$. Bounds on
$\Gamma_{p,\size}$ for various examples are given in \cref{sec:v1:examples}. These
examples show that in many cases $\Gamma_{p,\size}$ has polynomial decay.

We are now in position to state the main theorem of this section, whose proof is
deferred to \cref{sec:proofs-related-tail}.%
\begin{theorem}
  \label{thm:v1:main-thm}
  Under assumption \ref{ass:v1:5} and \ref{ass:v1:2}, there is a constant $C'' > 0$
  depending only on $p$ and $W$ such that for all $\size \geq 1$ it holds
  \begin{equation}
    \EE_W[(\hat{\sigma}_p(G_t) - \sigma)^2]
    \leq
    \frac{\Gamma_{p,\size}^2}{(-\log p)^2}
    + C'' \log(\size)^2 \max\Big\{\frac{1}{\size^{1+\sigma}\ell(\size)},\,
      \size^{1 - 2\eta} \Big\}.
  \end{equation}
\end{theorem}

\begin{remark}
  \label{rmk:2}
  The computation of the estimator requires to choose a value for
  $p$. Theoretical guidelines are that any value of $p$ bounded away from zero
  yields a consistent estimator. But, this leaves unclear how to choose in
  practice. Optimally, we could build a selection procedure that selects the
  ``optimal'' value of $p$, in the sense that bias and variance are suitably
  balanced. This would however require to make finer assumptions on the
  structure of the bias. We believe the gain would be marginal as in many cases
  the value of $\hat{\sigma}_p$ is relatively robust to the choice of $p$, see
  \cref{fig:1,fig:6}. In practice, we recommend to do a Hill plot to assess the
  pertinence of estimating $\sigma$ and to determine a reasonable value of $p$,
  as detailed in \cref{sec:assess-pert-exch}.
\end{remark}

\subsection{Discussion about the bias}
\label{sec:v1:examples}

In \cite{caron:panero:rousseau:2017}, the authors suggest to estimate the tail-index
using
\begin{equation}
  \hat{\sigma}_{\mathrm{CPR}} \coloneqq%
  \frac{2 \log \nver(G_t)}{\log \nedg(G_t)} - 1.
\end{equation}
It is easily seen that $\EE_W[(\hat{\sigma}_{\mathrm{CPR}} - \sigma)^2]$ is in
general a $O(1 / \log \size)$. Although the variance of $\hat{\sigma}_p$ has
polynomial decay in $\size$, the bias is proportional to $\Gamma_{p,\size}$,
whose decay can be anything. Hence it is not clear if $\hat{\sigma}_p$ dominates
$\hat{\sigma}_{\mathrm{CPR}}$ asymptotically. Here we present some examples where
the bias of $\hat{\sigma}_p$ has also polynomial decay in $\size$. We consider
the examples which were originally given in \cite{caron:panero:rousseau:2017},
because they cover a large spectrum of behaviours.
\begin{example}
  \label{ex:v1:1}
  \textit{Dense graphs}. Consider the function
  $W(x,y) = (1-x)(1-y)\Ind_{x\leq 1}\Ind_{y\leq 1}$. Then
  $\mu(x) = 0.5 (1-x)\Ind_{x\leq 1}$ and the graph is almost-surely
  dense. \Cref{ass:v1:5} is trivially satisfied with $\eta = 1$. Also, it follows
  from \cref{sec:v1:dense-graph-example} that \cref{ass:v1:2} is satisfied with
  $\Gamma_{p,\size} \propto \size^{-1}$: the risk is polynomially decreasing on
  this example.
\end{example}

% \subsubsection{Sparse, almost dense graphs without power law}
% \label{sec:v1:sparse-almost-dense}
\begin{example}
  \label{ex:v1:2}
  \textit{Sparse, almost dense graphs without power law}. Consider the function
  $W(x,y) = \exp(-x - y)$. Then $\mu(x) = \exp(-x)$ and \cref{ass:v1:5} is
  trivially satisfied with $\eta = 1$ and $\sigma = 0$. Moreover, it follows
  from \cref{sec:v1:sparse-almost-dense-1} that \cref{ass:v1:2} is satisfied with
  $\Gamma_{p,\size} \propto 1/\log(\size)$: the risk is logarithmically
  decreasing on this example.
\end{example}

% \subsubsection{Sparse graphs with power law, separable}
% \label{sec:v1:sparse-graphs-with-sep}
\begin{example}
  \label{ex:v1:3}
  \textit{Sparse graphs with power law, separable}. For $0< \sigma < 1$,
consider the function $W(x,y) = (1+x)^{-1/\sigma}(1+y)^{-1/\sigma}$. Then
obviously $\mu(x) = \sigma(1+x)^{-1/\sigma}/(1-\sigma)$ and \cref{ass:v1:5} holds
trivially with $\eta = 1$. Moreover, from \cref{sec:v1:sparse-graphs-with},
\cref{ass:v1:2} holds with $\Gamma_{p,\size} \propto \size^{-\sigma}$: the risk is
polynomially decreasing on this example.
\end{example}

% \subsubsection{Sparse graphs with power law, non separable}
% \label{sec:v1:sparse-graphs-with-nonsep}
\begin{example}
  \label{ex:v1:4}
  \textit{Sparse graphs with power law, non separable}. For $0 < \sigma < 1$,
  consider the function $W(x,y) = (1+x+y)^{-1/\sigma - 1}$. It is shown in
  \cite{caron:panero:rousseau:2017} that \cref{ass:v1:5} holds with
  $\eta = (1+\sigma)/2$. We have $\mu(x) = \sigma(1+x)^{-1/\sigma}$, which is up
  to a constant the same marginal as the previous example, so that \cref{ass:v1:2}
  is verified with $\Gamma_{p,\size} \propto \size^{-\sigma}$: the risk is
  polynomially decreasing on this example.
\end{example}

% \subsubsection{Generalized Gamma Process}
% \label{sec:v1:gener-gamma-proc}
\begin{example}
  \label{ex:v1:5}
  GGP model. For $0 \leq \sigma < 1$, let
  $\ufMeasure_{\sigma,\tau}(\intd x) = \frac{x^{-1-\sigma}\exp(-x) \intd
    x}{\Gamma(1-\sigma)}$,
  $\bar{\rho}_{\sigma,\tau}(x) \coloneqq \int_x^{\infty}\rho_{\sigma,\tau}(\intd
  x)$, and
  $W(x,y) = 1 - \exp(-2 \bar{\rho}_{\sigma,\tau}^{-1}(x)
  \bar{\rho}_{\sigma,\tau}^{-1}(y) )$. This corresponds to the model in
  \cite{caron:fox:2017}. It is shown in \cite{caron:panero:rousseau:2017} that
  \cref{ass:v1:5} holds with $\eta = 1$. We establish in
  \cref{sec:v1:gener-gamma-proc-1} that \cref{ass:v1:2} is satisfied with
  $\Gamma_{p,\size} \propto \size^{-\sigma}$: the risk is polynomially
  decreasing on this example.
\end{example}

Although the asymptotic behaviour of the graphs are relatively close for
examples \ref{ex:v1:1} and \ref{ex:v1:2}, the bias decay is radically different
and estimation is harder in the second example. We also estimated
$\EE[(\hat{\sigma}_{0.5} - \sigma)^2]^{1/2}$ and
$\EE[(\hat{\sigma}_{0.5} - \sigma)^2]^{1/2}$ on the previous examples using
$10000$ Monte Carlo samples. In the simulations, we picked $\sigma=0.3$ for the
sparse separable and sparse non-separable example, while the GGP example
corresponds to a choice of $\sigma = 0.5$. The raw numbers are given in tables
\ref{tab:v1:summary-risk} and \ref{tab:v1:summary-risk-CPR}.

\begin{table}[!tp]
  \centering
  \caption{Estimation of $\EE[(\hat{\sigma}_{0.5} - \sigma)^2]^{1/2}$
    using 10 000 MC samples for the examples described in \cref{sec:v1:examples}
    and various values of $\size$.}
  \label{tab:v1:summary-risk}
  \begin{tabular}{lcccccc}
%    \toprule
    $\size$ & 25 & 50 & 100 & 200 & 400 & 800\\% & 1600\\
%    \midrule
    Dense   & 0.1620 & 0.0698 & 0.0327 & 0.0150 & 0.0079 & 0.0039\\% & 0.0019 \\
    Almost-dense & 0.2971 & 0.2440 & 0.2081 & 0.1814 & 0.1608 & 0.1443\\
    Sparse sep. ($\sigma=0.3$) & 0.2440 & 0.1658 & 0.1186 & 0.0875 & 0.0654 & 0.0488 \\
    Sparse non-sep. ($\sigma = 0.3$) & 0.2947 & 0.1945 & 0.1356 & 0.0966 & 0.0686 & 0.0462 \\
    GGP ($\sigma=0.5$) & 0.1219 & 0.0801 & 0.0502 & 0.0283 & 0.0110 & 0.0050\\
%    \bottomrule
  \end{tabular}
\end{table}

\begin{table}[!tp]
  \centering
  \caption{Estimation of $\EE[(\hat{\sigma}_{\mathrm{CPR}} - \sigma)^2]^{1/2}$
    using 10 000 MC samples for the examples described in \cref{sec:v1:examples}
    and various values of $\size$.}
  \label{tab:v1:summary-risk-CPR}
  \begin{tabular}{lcccccc}
%    \toprule
    $\size$ & 25 & 50 & 100 & 200 & 400 & 800\\% & 1600\\
%    \midrule
    Dense   & 0.2719 & 0.2125 & 0.1757 & 0.1500 & 0.1305 & 0.1154\\% & 0.1036 \\
    Almost-dense & 0.4240 & 0.3882 & 0.3590 & 0.3349 & 0.3144 & 0.2967\\
    Sparse sep. ($\sigma=0.3$) & 0.2963 & 0.2535 & 0.2230 & 0.1994 & 0.1802 & 0.1643\\
    Sparse non-sep. ($\sigma=0.3$) & 0.3483 & 0.2947 & 0.2575 & 0.2285 & 0.2051 & 0.1851\\
    GGP ($\sigma=0.5$) & 0.1805 & 0.1577 & 0.1392 & 0.1239 & 0.1103 & 0.0982\\
%    \bottomrule
  \end{tabular}
\end{table}

We pictured in \cref{fig:v1:our-vs-cr} the risk of $\hat{\sigma}_{0.5}$
and $\hat{\sigma}_{\mathrm{CPR}}$ as functions of $\size$ for all the simulated
examples. As expected, the risk of $\hat{\sigma}_{0.5}$ has polynomial
decay every-time but the almost dense example, while the risk of
$\hat{\sigma}_{\mathrm{CPR}}$ remains logarithmically decreasing on all
examples. We note, however, that even in the situation when the risk has
logarithm decay, $\hat{\sigma}_{0.5}$ seems to perform slightly better
than $\hat{\sigma}_{\mathrm{CPR}}$.

\begin{figure}[!tp]
  \centering

  \begin{tikzpicture}
    \begin{groupplot}[%
      xlabel=$\size$,%
      ylabel=Risk,%
      height=5.5cm,%
      minor x tick num=1,%
      xtick pos=left,%
      ytick pos=left,%
      enlarge x limits=true,%
      every x tick/.style={color=black, thin},%
      every x tick label/.append style={font=\tiny, yshift=0.5ex},%
      every y tick/.style={color=black, thin},%
      every y tick label/.append style={font=\tiny, xshift=0.5ex},%
      tick align=inside,%
      xlabel near ticks,%
      ylabel near ticks,%
      axis on top,%
      legend pos=south west,%
      group style={group size=2 by 3, ylabels at=edge left, vertical sep=1.8cm},%
      ]

      % DENSE EXAMPLE
      % ------------------------------------------------------------------------
      \nextgroupplot[ymode=log,xmode=log,title=Dense]%

      % Our estimator
      \addplot+[color=blue,mark=*] coordinates {%
        (25, 0.1620)%
        (50, 0.0698)%
        (100, 0.0327)%
        (200, 0.0150)%
        (400, 0.0079)%
        (800, 0.0039)%
        (1600, 0.0019)%
      };

      % Caron and Rousseau
      \addplot+[color=red,mark=*] coordinates {%
        (25, 0.2719)%
        (50, 0.2125)%
        (100, 0.1757)%
        (200, 0.1500)%
        (400, 0.1305)%
        (800, 0.1154)%
        (1600, 0.1036)%
      };

      \legend{ $\hat{\sigma}_{0.5}$, $\hat{\sigma}_{\mathrm{CPR}}$ }

      % ALMOST DENSE EXAMPLE
      % ------------------------------------------------------------------------
      \nextgroupplot[xmode=log,ymode=log,title=Almost dense]%

      % Our estimator
      \addplot[color=blue,mark=*] coordinates {%
        (25, 0.2971)%
        (50, 0.2440)%
        (100, 0.2081)%
        (200, 0.1814)%
        (400, 0.1608)%
        (800, 0.1443)%
      };

      % Caron and Rousseau
      \addplot[color=red,mark=*] coordinates {%
        (25, 0.4240)%
        (50, 0.3832)%
        (100, 0.3590)%
        (200, 0.3349)%
        (400, 0.3144)%
        (800, 0.2967)%
      };

      \legend{ $\hat{\sigma}_{0.5}$, $\hat{\sigma}_{\mathrm{CPR}}$ }

      % SPARSE SEPARABLE EXAMPLE (sigma=0.3)
      % ------------------------------------------------------------------------
      \nextgroupplot[xmode=log,ymode=log,title={Sparse sep., $\sigma=0.3$}]%

      % Our estimator
      \addplot[color=blue,mark=*] coordinates {%
        (25, 0.2440)%
        (50, 0.1658)%
        (100, 0.1186)%
        (200, 0.0875)%
        (400, 0.0654)%
        (800, 0.0488)%
      };

      % Caron and Rousseau
      \addplot[color=red,mark=*] coordinates {%
        (25, 0.2963)%
        (50, 0.2535)%
        (100, 0.2230)%
        (200, 0.1994)%
        (400, 0.1802)%
        (800, 0.1643)%
      };

      \legend{ $\hat{\sigma}_{0.5}$, $\hat{\sigma}_{\mathrm{CPR}}$ }

      % SPARSE NON-SEPARABLE EXAMPLE (sigma=0.3)
      % ------------------------------------------------------------------------
      \nextgroupplot[xmode=log,ymode=log,title={Sparse non-sep., $\sigma=0.3$}]%

      % Our estimator
      \addplot[color=blue,mark=*] coordinates {%
        (25, 0.2947)%
        (50, 0.1945)%
        (100, 0.1356)%
        (200, 0.0966)%
        (400, 0.0686)%
        (800, 0.0462)%
      };

      % Caron and Rousseau
      \addplot[color=red,mark=*] coordinates {%
        (25, 0.3483)%
        (50, 0.2947)%
        (100, 0.2575)%
        (200, 0.2285)%
        (400, 0.2051)%
        (800, 0.1851)%
      };

      \legend{ $\hat{\sigma}_{0.5}$, $\hat{\sigma}_{\mathrm{CPR}}$ }

      % GGP EXAMPLE (sigma=0.5)
      % ------------------------------------------------------------------------
      \nextgroupplot[xmode=log,ymode=log, title={GGP, $\sigma=0.5$}]%

      % Our estimator
      \addplot[color=blue,mark=*] coordinates {%
        (25, 0.1219)%
        (50, 0.0801)%
        (100, 0.0502)%
        (200, 0.0283)%
        (400, 0.0110)%
        (800, 0.0050)%
      };

      % Caron and Rousseau
      \addplot[color=red,mark=*] coordinates {%
        (25, 0.1805)%
        (50, 0.1577)%
        (100, 0.1392)%
        (200, 0.1239)%
        (400, 0.1103)%
        (800, 0.0982)%
      };

      \legend{ $\hat{\sigma}_{0.5}$, $\hat{\sigma}_{\mathrm{CPR}}$ }
    \end{groupplot}
  \end{tikzpicture}

  \caption{Comparative plot of
    $\EE[(\hat{\sigma}_{0.5} - \sigma)^2]^{1/2}$ against
    $\EE[(\hat{\sigma}_{\mathrm{CPR}} - \sigma)^2]^{1/2}$ using 10 000 MC samples
    for the examples described in \cref{sec:v1:examples}.}
  \label{fig:v1:our-vs-cr}
\end{figure}

\subsection{On the controversy about scale-free nature of real graphs}
\label{sec:contr-opin-about}

There is no actual consensus on the scale-free nature of empirical
graphs\footnote{Here empirical refers to typical datasets encountered in the
  literature. This is of course very subjective and thus one of the source of
  the controversy.}. Here, we briefly review both sides of the debate. This is
also an opportunity for us to review what $\sigma$ is and what it isn't.

According to the general opinion, a finite graph $G$ is called scale-free if there are
$c > 0$, $\gamma > 1$ such that for $j$ large enough%
\begin{equation}
  \label{eq:34}
  \frac{\degree_j(G)}{\nver(G)} \approx cj^{-\gamma}.
\end{equation}
We left the previous definition intentionally non-rigorous to emphasize one of
the origin of the controversy: the lack of consensus on the definition of
scale-freeness
\citep{holme:2019,voitalov:hoorn:hofstad:2019}. % In fact, according to this definition,
% a sample $G_t$ from a graphex process with $\sigma \in (0,1)$ might be considered as not having a
% power-degree distribution, see \cref{fig:2} for illustration.

\citet{barabasi:albert:1999} found that for three empirical networks
\cref{eq:34} seemed true, and they attributed this property to a growth
mechanism they called \textit{preferential attachment}. Since then, many papers
claimed to have discovered new types of power-law networks and new mechanism
creating them, see for instance \cite{newman:2005,holme:2019}. But not everybody
agrees, and some authors believe that power-law degree distribution is rarely
encountered in practice
\cite{jones:handcock:2003,clauset:shalizi:newman:2009}. The heated debate has
culminated recently with \cite{broido:clauset:2019}, claiming that
\textit{scale-free networks are rare} after reviewing a large collection of
empirical networks.

The methodology in \cite{broido:clauset:2019} consists on testing as a null
hypothesis a power-law model for the degree distribution against an alternative
which is not power-law. On a corpus of 927 empirical networks, which they claim
is representative of \textit{real} networks, they found that only a very few
failed to reject the null hypothesis.

Since its first appearance as a preprint \cite{broido:clauset:2019} has
generated a lot of activity, whether it be in agreement or disagreement. We
share the same opinion as \cite{holme:2019,voitalov:hoorn:hofstad:2019} that
the source of the debate is the lack of a rigorous
definition. \cite{voitalov:hoorn:hofstad:2019} makes an effort toward a clean
definition and they arrive at the opposite conclusion of
\cite{broido:clauset:2019}. Another argument against \cite{broido:clauset:2019},
which seems to be the most widespread, is that scale-freeness is only
well-defined in the infinite-size limit \cite{holme:2019}. According to this
view, a network is scale-free if its degree distribution approaches a power-law
as the network keeps growing following the same mechanisms. In other word,
instead of \cref{eq:34}, we shall say that $\Set{G_t \given t \geq 0}$ is
scale-free if
\begin{equation}
  \label{eq:35}
  \lim_{j\to \infty}\lim_{t\to \infty}\frac{j^{\gamma}\degree_j(G_t)}{\nver(G_t)} = c.
\end{equation}
Then, $\gamma$ is not defined as property of the sample $G_t$, but rather as a
joint property of the limiting graph and the growth mechanism, equivalently of
the sampling design. Here we insist on the importance of the sampling design in
the definition: if $\Set{G_t \given t \geq 0}$ is a graphex process, then
$\gamma = 1 + \sigma$ almost-surely under mild conditions
\cite{caron:panero:rousseau:2017}, though the degree distribution of $G_t$ at
finite $t$ might have a different index, as illustrated in \cref{fig:2}. Another
famous example is the \textit{traceroute} sampling which is known to produce
samples with power-law degree distributions, even when the population is not
power-law
\cite{lakhina:byers:crovella:2003,achlioptas:clauset:kempe:2009}. Another
example of this type is provided in \cref{sec:test-train-split}. Hence, if one
is ready to accept \cref{eq:35} as a definition, talking about scale-freeness
without knowledge of the sampling design is irrelevant.

One of the criticism of the definition in \cref{eq:35} is made on the basis that
we never observe $\Set{G_t \given t \geq 0}$ but only $G_t$, and that we need
tools adapted to finite networks. We have mixed feelings about this. On one
hand, the same criticism could be made about any statistical analysis. Indeed,
one of the primary goal of statistics is to infer stuff about a population given
a sample. In the case of i.i.d observations, we have accepted for a while to
idealize the population as an abstract propability distribution, and this has
been proven to be very effective in summarizing properties of the population as
parameters of the distribution. When the observation is a graph, we have no
reason to proceed differently (see in particular in \cite{holme:2019} the
example of epidemic threshold). On the other hand \cite{holme:2019} points
out the difficulty in determining the right sampling mechanism, and it is true
that models of random networks are commonly misspecified. Then, it is less clear
that the parameters make sense, albeit this depends on the degree of
misspecification and the task to achieve.

Finally, we note that the test in \cite{broido:clauset:2019} only involves the
degree distribution, but graphs are more than just their degree distribution
(see in particular \cref{sec:test-count-stat}). For these reason, we don't
believe in the existence of a universal model that can explain all empirical
graphs, and it is the work of the statistician to determine a good model for
\textit{a given task}; whence it is useful to study as many models as
possible. This opinion seems to be shared by the authors of
\cite{broido:clauset:2019}. We add that in many cases the focus should be put on
achieving the task at hand with accuracy rather than trying to explain all the
structure of the sample. Besides, the assertion that graphex processes are not
desirable to model empirical graphs because power-law graphs are rare is too
simplistic: (i) no serious methodology can objectively prove this rarity (what
is a representative set of empirical graphs? what definition of scale-freeness
to use?), (ii) it is the degree distribution of graphex processes in the
infinite size limit which is power-law\footnote{And in fact, they can
  accommodate for other behaviours, by either satisfying $\sigma = 0$, or
  violating \cref{ass:v1:5,ass:v1:2}.}, and (iii) there are empirical graphs for
which some graphex processes have been found to explain well the degree
distribution (see \cite{caron:fox:2017,todeschini:miscouridou:caron:2016} and
also \cref{sec:some-appl-exampl}). Overall, one should not forget to look at
other aspects of the modeling that might be more important than just the degree
distribution.

\subsection{Relation to Hill estimator}
\label{sec:relat-hill-estim}

\begin{figure}[h]
  \centering

  \includegraphics[width=\linewidth]{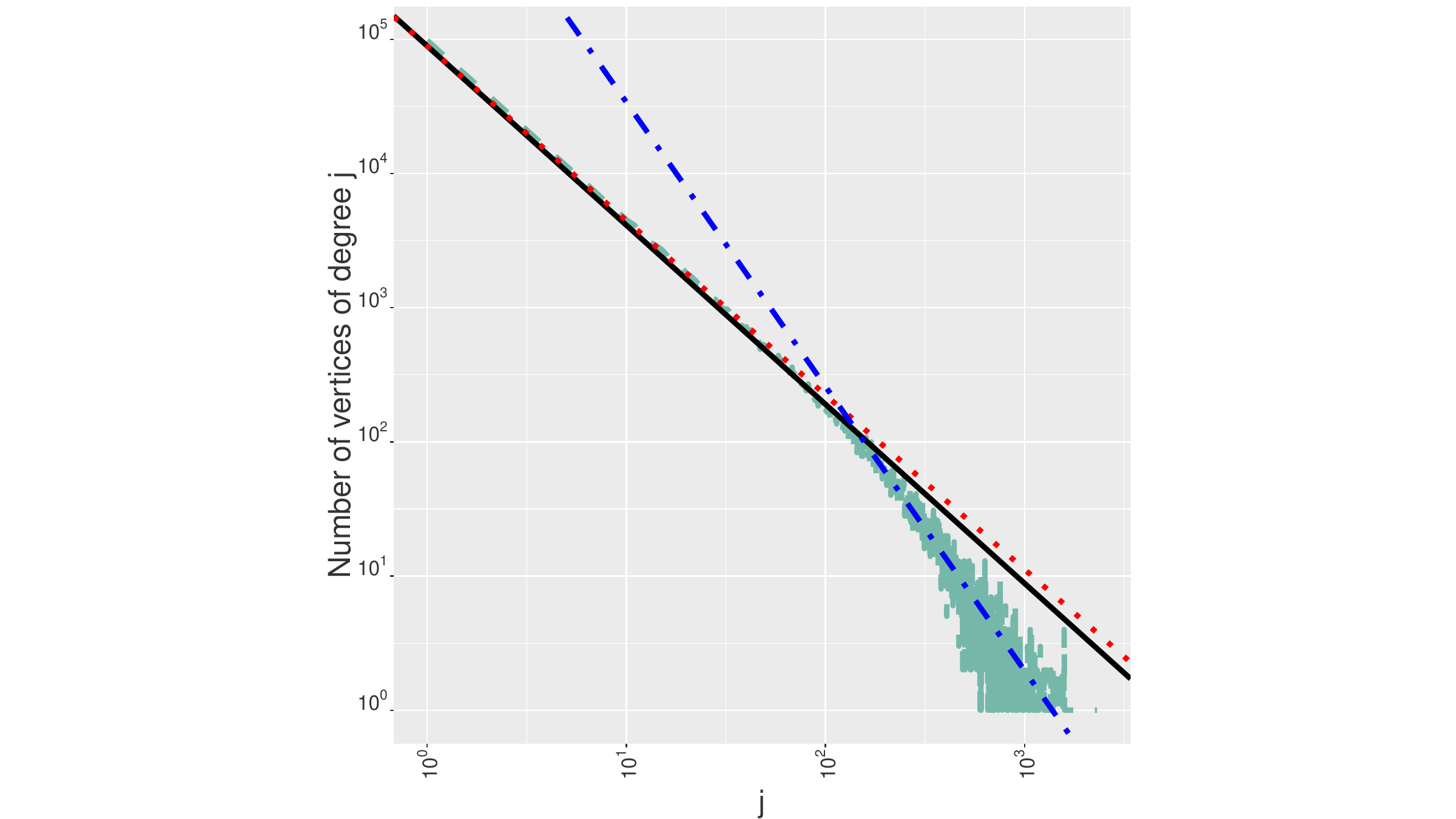}
  
  \caption{%
    In green is the degree distribution $G_t$ is a sample from the GGP model of
    \cite{caron:fox:2017} with parameters $(\sigma,\tau) = (0.3,10)$ and
    $t = 10000$ (see also \cref{ex:v1:5}). The sample has 306791 vertices and
    3947805 edges. The red line represents the
    $\degree_j \propto j^{-1-\sigma}$. The black line is
    $\degree_j \propto j^{-1-\hat{\sigma}_p}$, where $\hat{\sigma}_p = 0.333474$
    is computed using $p=0.75$. Finally the blue line is
    $\degree_j \propto j^{-\hat{\gamma}_{\mathrm{Hill}}}$, where
    $\hat{\gamma}_{\mathrm{Hill}} = 2.1181$ is computed using $d_{\min} = 300$.}
  \label{fig:2}
\end{figure}

Regardless of the controversy of \cref{sec:contr-opin-about}, tail-index of the
degree distribution of $G_t$ is often reported as an interesting summarizing
statistics \cite[Section~4]{kolaczyk:2009}. \cite{voitalov:hoorn:hofstad:2019}
plea in favor of using various estimators as the index is notoriously hard to
estimate and different estimators may give very different answers on a given
sample. In our opinion, the interest in the tail-index of $G_t$ is debatable as
it might be unrelated with the index of the population, an issue we already
discussed in \cref{sec:contr-opin-about}, and which we will emphasize here. The
most common reported value for the tail-index of the sample is based on the Hill
estimator
\begin{equation}
  \label{eq:36}
  \hat{\gamma}_{\mathrm{Hill}}(G_t) \coloneqq%
  \frac{\sum_{j=d_{\min+1}}^{\infty} \degree_j(G_t) \log \frac{j}{d_{\min}}
  }{\sum_{j=d_{\min}+1}^{\infty}\degree_j(G_t) },
\end{equation}
where $d_{\min}> 0$ is a cutoff determining where the tail of the degree
distribution starts. Most of empirical graphs are reported to have
$\hat{\gamma}_{\mathrm{Hill}} \in [2,3]$, while graphex processes have
almost-surely index $\gamma = 1 + \sigma \in [1,2]$, at least if we rely on the
definition of \cref{eq:35}. So one might think graphex processes are unrealistic
to model graphs with $\hat{\gamma}_{\mathrm{Hill}} \in [2,3]$. Once again, we
have to be careful: it is not obvious that the Hill estimator gives consistent
answers for $\gamma$. Though it is known from \cite{wang:resnick:2019} that it
is consistent in the preferential attachment model, nothing is known about
graphex processes. In particular, in \cref{fig:2} we illustrate what can go
wrong: on a simulated sample with true $\sigma = 0.3$ (hence $\gamma = 1.3$),
our estimator gives $\hat{\sigma}_{0.5} = 0.3335$, while a non carefully
calibrated Hill estimator gives
$\hat{\gamma}_{\mathrm{Hill}} = 2.1181 \in [2,3]$. Once again, this shows that
one has to be careful when relating degree distribution of the sample and the
idealized degree distribution of the population, \textit{i.e.}  the sampling
design has to be taken into account when inferring properties about the
population.

\section{Count statistics of graphex processes: definition, inference, and
  hypothesis testing}
\label{sec:bootstrap-motifs}

We briefly discuss how the $p$-sampling argument may also be used to estimate
other functionals of $W$. This essentially extends the work of
\cite{bickel:chen:levina:2011,bhattacharyya:bickel:2015} from graphon to
graphex. In particular, using $p$-sampling allows to derive a natural notion of
bootstrap for networks. Using bootstrap samples, we can then estimate the
\textit{count statistics} of the graph
(\cref{sec:statement-problem,sec:p-sampling-as-bs}), which can in turn be used
to do nonparametric inference on the (rescaled injective) homomorphism
densities (\cref{sec:p-sampling-as-bs-pop}). We also use count statistics to
construct tests for the number of occurrence of certain motifs
(\cref{sec:hypothesis-testing}). Finally, we discuss the relation of these
results with some other works (\cref{sec:relation-other-works}).
 
\subsection{Count statistics}
\label{sec:statement-problem}

Count statistics have received a lot of attention in the last few years, for
multiple reasons. A first reason is that some motifs are of interest to get
insights on structure of the network. For instance, counting the number of
triangles permit to determine the \textit{transitivity} of the graph, which is a
popular measure of the tendency of the vertices to cluster together
\citep{watts:strogatz:1998}. Another reason is that counts of motifs
characterize the distribution of the network (see \cite{diaconis:janson:2008}
and \cref{sec:p-sampling-as-bs-pop}), so that nonparametric inference on the
distribution of the network is possible using counts of certain motifs
\cite{bickel:chen:2009,bickel:chen:levina:2011,bhattacharyya:bickel:2015}. Other
related uses of motifs counts involve hypothesis testing
\cite{bhattacharyya:bickel:2015} and networks comparison
\cite{bhattacharyya:bickel:2015,ali:wegner:gaunt:2016}.

Count of motifs is notoriously hard, especially when the graph is large and the
motif has a large number of vertices. Even counting simple motifs such has
triangles might rapidly become intractable on large graphs
\cite{suri:vassilvitskii:2011}. Methods based on subsampling the graph have
emerged to bypass the complexity. A non exhaustive list of existing methods
includes methods based on vertex sampling
\cite{bhattacharyya:bickel:2015,klusowski:wu:2018,ali:wegner:gaunt:2016,lunde:sarkar:2019},
the parametric bootstrap \cite{green:shalizi:2017,levin:levina:2019}, or based
on edge sampling
\cite{assadi:kapralov:khanna:2018,eden:levi:ron:2017,goldreich:ron:2008,gonen:ron:shavitt:2011}. Here,
we take advantage of the $p$-sampling to count motifs. We discuss in
\cref{sec:relation-other-works} how that relates to already existing works, and
some of the advantages in doing $p$-sampling over other sampling design.

Before going further, we shall first clarify what me mean by count
statistics. For a finite graph $R$, we define the number of copies of $R$ in $G$
as the number of injections $\phi : \ver(R) \to \ver(G)$ such that if
$\Set{r_1,r_2} \in \edg(R)$ then $\Set{\phi(r_1),\phi(r_2)} \in \edg(G)$,
\textit{i.e.} the number of adjacency preserving map from $R$ to $G$. We call
this number $\nsub(R,G)$, the count statistic of $R$. If $\phi : A \to \ver(G)$
is injective, we let $G[\phi]$ denote the graph such that $\ver(G[\phi]) = A$
and $\Set{v,v'} \in \edg(G[\phi])$ if and only if
$\Set{\phi(v),\phi(v')} \in \edg(G)$. That is, $G[\phi]$ is isomorphic to the
subgraph of $G$ induced by $\phi$. Letting $\mathcal{A}(R,G)$ denote the set of
all injections from $\ver(R)$ to $\ver(G)$, we can rewrite
\begin{align}
  \label{eq:54}
  \nsub(R,G)%
  = \sum_{\phi \in \mathcal{A}(R,G)}\1\Set{R \subseteq G[\phi]}.
\end{align}

% --------------------------------------------------------------------------------

% For a subgraph $R$, we denote the
% number of copies of $R$ contained in $G_t$ as induced subgraph by
% $\indsub(R,G_t)$. More precisely, we say that two graphs $G,G'$ are isomorphic,
% denoted by $G \cong G'$, if there exists a bijection between
% $f : \ver(G) \to \ver(G')$ that preserves the adjacency, \textit{i.e.}
% $\Set{f(u),f(v)} \in \edg(G')$ if and only if $\Set{u,v} \in \edg(G)$. Then, for
% any $V \subseteq \ver(G)$, the restriction of $G$ to $V$ is written
% $G[V] \coloneqq (V, \edg(G) \cap V^2)$, and the number of copies of $R$
% contained in $G_t$ is
% \begin{equation}
%   \label{eq:7}
%   \indsub(R,G)%
%   \coloneqq \sum_{V\subseteq \ver(G_t)}\1\Set{G[V] \cong R}.
% \end{equation}

\subsection{Bootstrap estimators for count statistics: model-free results}
\label{sec:p-sampling-as-bs}

We first show how to derive a bootstrap estimate for $\nsub(R,G)$ using
$p$-sampling which works for any graph $G$ (so $G$ is not necessarily a graphex
process). The result is of interest by itself, since it proves a consistent
estimator for $\nsub(R,G)$ regardless of the distribution of $G$, and that
estimator is parallelizable, in contrast with the computation of $\nsub(R,G)$
which may be challenging for complex motifs $R$. Our estimator is based on the
following easy proposition.

\begin{proposition}
  \label{pro:2}
  Let $G \in \mathcal{U}$ and let $R \in \mathcal{L}_k$ be a connected graph
  such that $\nsub(R,G) < \infty$. Then,
  \begin{equation}
    \label{eq:5}
    \EE[\nsub(R, \psamp(G,p)) \mid G]%
    = p^k \cdot  \nsub(R,G).
  \end{equation}
\end{proposition}

From the previous proposition, we can obtain a natural estimator of
$\nsub(R,G)$, which we refer as the \textit{bootstrap} estimator. Indeed,
for a fixed $p \in (0,1)$ to be chosen, and for some integer $B \geq 1$, we let
$\hat{G}_1,\dots,\hat{G}_B$ be i.i.d copies of $ \psamp(G,p) \mid G$,
and we estimate $\nsub(R,G)$ by
\begin{equation}
  \label{eq:8}
  \hat{\nsub}(R,G)%
  \coloneqq \frac{1}{B p^k}\sum_{b=1}^B \nsub(R,\hat{G}_b).
\end{equation}
The main advantage is that $\nsub(R,\hat{G}_b)$ are intrinsically faster to
compute than $\nsub(R,G)$, especially when $G$ is sparse, and computations can
be run in parallel. Interestingly, the formula \eqref{eq:5} is
\textit{model-free}, in the sense that it is true for any graph $G$ that has no
isolated vertices. It follows that if $G$ has no isolated vertices, then
$\hat{\nsub}(R,G)$ is an unbiased estimator of $\nsub(R,G)$, comparable to
the Horvitz-Thompson estimator constructed in
\cite{klusowski:wu:2018}\footnote{However \cite{klusowski:wu:2018}
  consider the number of induced subgraphs of $G$ that are isomorphic to $R$,
  which is slightly different from $\nsub(R,G)$, but the arguments are comparable.}.
The next proposition quantify the risk of $\hat{\nsub}(R,G)$ as an estimator
of $\nsub(R,G)$.

\begin{proposition}
  \label{pro:3}
  Let $G \in \mathcal{U}$ and let $R \in \mathcal{L}_k$ be a connected graph
  such that $\nsub(R,G) < \infty$. Then $\hat{\nsub}(R,G)$ is an unbiased
  estimator of $\nsub(R,G)$, and%
  \begin{equation}
    \label{eq:9}
    \var(\hat{\nsub}(R,G) \mid G) =%
    \frac{\var(\nsub(R,\psamp(G,p)) \mid G)}{Bp^{2k}}%%
    \leq \frac{k }{Bp} \cdot \nsub(R,G)^2.
  \end{equation}
\end{proposition}
It is worth pointing out that by the previous, $\hat{\nsub}(R,G)$ is a
consistent estimator of $\nsub(R,G)$ as $B \to \infty$, and thus with almost
no condition on $G$. This result is of interest by itself, if one is only
interested in counting motifs in $G$ (see also \cref{sec:counting-motifs}).

\begin{remark}
  \label{rmk:4}
  In contrast with our estimator for the tail-index, we keep the extra
  randomness induced by the bootstrap samples, i.e. we have a randomized
  estimator. But, it is clear from \cref{pro:1} that taking the expectation of
  $\hat{\nsub}(R,G)$ conditional on $G$ is pointless. Moreover, the advantage of
  keeping the randomized-estimator is firstly to enable tractable computations
  as bootstrap samples can be made arbitrary small.
\end{remark}

\subsection{Nonparametric inference}
\label{sec:p-sampling-as-bs-pop}

Beyond their interest in describing the structure of the network, count
statistics permit to estimate the rescaled injective homomorphism densities
\cite{borgs:chayes:cohn:holden:2016}. Those are distinguished functionals of $W$
whose collection entirely determine $W$
\cite{diaconis:janson:2008,lovasz:2013,borgs:chayes:cohn:holden:2016}. The
rescaled injective homomorphism density of $R$ in $W$ is defined as
\begin{align}
  \label{eq:17}
  \cstat(R,W)%
  &\coloneqq%
    \frac{%
    \int_{\NNReals^{\nver(R)}} \prod_{\Set{\ell,m}\in \edg(R)
    }W(x_\ell,x_m)\,\intd x_1\dots \intd x_{\nver(R)}%
    }{%
     \{\int_{\NNReals^2}W(x_1,x_2)\,\intd x_1 \intd x_2 \}^{\nver(R)/2}%
    }.
\end{align}
It follows from
\cite{diaconis:janson:2008,lovasz:2013,borgs:chayes:cohn:holden:2016} that $W$
can be in principle recovered -- up to measure preserving transformations --
from the collection of all $\Set{\cstat(R,W) \given R \in \mathcal{L}}$. We
shall only consider motifs $R$ and processes for which the previous quantity is
finite. Indeed, we will assume that a certain number of moments of $W$ are
finite.
\begin{assumption}
  \label{ass:expecmotifs}
  A graphex process $\Set{G_t \given t \geq 0}$ with graphon function $W :
  \NNReals^2 \to [0,1]$
  satisfies this assumption for $N \geq 1$ integer if,
  \begin{equation}
    \label{eq:44}
    \max_{n=1,\dots, N}\max_{\substack{R \in \mathcal{L}_n\\R\,\mathrm{connected}}}\int_{\NNReals^n}\prod_{\Set{\ell,m}\in \edg(R)}
    W(x_{\ell},x_m)\, \intd x_1\dots \intd x_n < \infty.%
  \end{equation}
\end{assumption}
Then, we have the following result, borrowed from
\cite[Proposition~56]{borgs:chayes:cohn:holden:2016}.
\begin{proposition}
  \label{pro:4}
  Let $R \in \mathcal{L}_k$ be a connected graph. If $\Set{G_t \given t\geq 0}$
  is a graphex process satisfying \cref{ass:expecmotifs} with $N = k$, then
  almost-surely,
  \begin{equation}
    \label{eq:4}
    \lim_{t\to \infty}t^{-k}\nsub(R,G_t)%
    = \int_{\NNReals^k} \prod_{\Set{\ell,m}\in \edg(R) }W(x_\ell,x_m)\,\intd
    x_1\dots \intd x_k.
  \end{equation}
\end{proposition}

The consequence of the last proposition is that $\cstat(R,W)$ can be
consistently estimated as $t\to \infty$ by
\begin{equation}
  \label{eq:3}
  \hat{\cstat}(R,G_t)%
  \coloneqq%
 \frac{\nsub(R,G_t)}{\nedg(G_t)^{\nver(R)/2}},
\end{equation}
where in turn we can leverage the consistency of our bootstrap estimator of
$\nsub(R,G_t)$ to make computations tractable. We note that the method of
moments of \cite{bickel:chen:levina:2011} can be adapted to estimate $W$ from
$\hat{\cstat}(R_1,G_t),\hat{\cstat}(R_2,G_t),\dots$ for clever choice of motifs
$R_1,R_2,\dots$.

\subsection{Hypothesis testing}
\label{sec:hypothesis-testing}

Besides their interest for estimating $W$, the count statistics have been
used in testing equality of features of networks and finding confidence
intervals of the count feature
\citep{middendorf:ziv:wiggins:2005,shen-orr:milo:mangan:2002}. The authors of
\cite{bhattacharyya:bickel:2015} have an example where they use counts of
certain motifs to test if the dearth of short-cycles in a network is
statistically significant. As argued in \cite{bhattacharyya:bickel:2015}, the
main challenge is then the problem of finding estimates of
$s_t^2(R) \coloneqq \var_W(\nsub(R,G_t))$. Here we show how the $p$-sample
argument can be used in order to recycle the bootstrap samples to estimate
$s_t^2(R)$, which provides a somewhat simpler estimator than in
\cite{bhattacharyya:bickel:2015}. Our estimator is based on the result of the
following proposition.
\begin{proposition}
  \label{pro:1}
  Let $\Set{G_t \given t \geq 0}$ be a graphex process satisfying
  \cref{ass:expecmotifs} with $N = 2k-1$. For any connected graph
  $R \in \mathcal{L}_{k}$, as $t \to \infty$
  \begin{align}
    \label{eq:25}
    s_t^2(R)%
    &=
    \frac%
    {\EE_W[\var(\nsub(R,\psamp(G_t,p)) \mid G_t)]}
      {p^{2k-1}(1-p)}%
      + O(t^{2k-2}).
  \end{align}
\end{proposition}

This suggests that $s_t^2(R)$ may be estimated by
\begin{equation}
  \label{eq:26}
  \hat{s}^2(R,G_t)%
  \coloneqq%
  \frac%
  {\var(\nsub(R,\psamp(G_t,p)) \mid G_t)}
  {p^{2k-1}(1-p)},
\end{equation}
where in turn the variance in the numerator of the last display can be
approximated up to arbitrary precision by recycling the bootstrap samples used
for the estimation of $\nsub(R,G_t)$. The next proposition establishes that
$\hat{s}^2(R)$ is consistent.

\begin{proposition}
  \label{pro:5}
  Let $\Set{G_t \given t \geq 0}$ be a graphex process satisfying
  \cref{ass:expecmotifs} with $N = 4k-2$. For any connected graphs
  $R \in \mathcal{L}_{k}$, as $t \to \infty$
  \begin{equation}
    \label{eq:48}
    \frac{\hat{s}^2(R,G_t)}{s_t^2(R)} = 1 + O_p(t^{-1/2}).
  \end{equation}
\end{proposition}

The fact that we are able to consistently estimate $s_t^2(R)$ permit to
perform hypothesis testing for
$\bar{\nsub}(R,W,t) \coloneqq \EE_W[\nsub(R,G_t)]$. In particular, it is of
interest to construct tests for
$H_0 : \bar{\nsub}(R,W,t) = \bar{\nsub}(R,W_0,t_0)$ against
$H_a : \bar{\nsub}(R,W,t) \ne \bar{\nsub}(R,W_0,t_0)$, see for instance
\cref{sec:test-count-stat} and also \cite{bhattacharyya:bickel:2015}. We show in
the next proposition that asymptotically consistent testing for these hypotheses
can be made on the basis of the test statistic
\begin{align}
  \label{eq:58}
  T(R,G_t)%
  &\coloneqq%
    \frac{\nsub(R,G_t) - \bar{\nsub}(R,W_0,t_0)}{\sqrt{\hat{s}^2(R,G_t)}}.
\end{align}

\begin{proposition}
  \label{pro:6}
  Let $\Set{G_t \given t \geq 0}$ be a graphex process with graphon function
  $W_0$ satisfying \cref{ass:expecmotifs} with $N = 4k$. For any connected
  graphs $R \in \mathcal{L}_{k}$, $T(R,G_{t_0}) \xrightarrow{d} \mathcal{N}(0,1)$ as
  $t_0 \to \infty$.
\end{proposition}

% We illustrate in \cref{fig:4} the asymptotic normality of the test statistic on
% a simulation example.%

% \PROBLEM{zn: I think there is a mistake because the cpp code has $t=800$, and
%   furthermore the shapiro test fails, which is bad to illustrate asymptotic
%   normality...}

% \begin{figure}[h]
%   \centering

%   \includegraphics[width=\linewidth]{images/distri-teststat.pdf}
  
%   \caption{%
%     The histogram and the QQ-plot of the test statistic $T(R,G_t)$ for $R$ the
%     triangle graph and $G_t$ is a sample from the GGP model of
%     \cite{caron:fox:2017} with parameters $(\sigma,\tau) = (0.2,0.5)$ and
%     $t = 10000$ (see also \cref{ex:v1:5}). The histogram and the QQ-plot are
%     made on the basis of 1000 independent Monte Carlo realizations of $G_t$. For
%     each Monte Carlo sample, $T(R,G_t)$ has been calculated using bootstrap
%     estimates for $\nsub(R,G_t)$ and $\hat{s}^2(R,G_t)$ with $p=0.1$ and 200
%     bootstrap samples. We have superposed the density of the $\mathcal{N}(0,1)$
%     (on the histogram) distribution and the line $y=x$ (on the QQ-plot) to
%     emphasize the asymptotic normality of $T(R,G_t)$.}
%     \label{fig:4}
% \end{figure}

\begin{remark}
  \label{rmk:1}
  It would also be of interest to construct a test for the weaker hypothesis
  $H_0 : \cstat(R,W) = \cstat(R,W_0)$. It seems to us that this is a much harder
  problem. Indeed, it is true that $\hat{\cstat}(R,G_t)\sim \cstat(R,W_0)$
  almost-surely under $H_0$ by \cite{borgs:chayes:cohn:holden:2016}, and so we
  could imagine constructing a test based on $\hat{\cstat}$. However, the
  statistical fluctuations of $\nsub(R,G_t) - \EE_W[\nsub(R,G_t)]$ and of
  $\nedg(G_t) - \EE_W[\nedg(G_t)]$ are exactly of the same order, which makes the
  asymptotic normality of $\hat{\cstat}(R,G_t)$ less evident.
\end{remark}

\subsection{Relation to other works}
\label{sec:relation-other-works}

The work presented here relates to many of the paper cited previously. In our
opinion, the two most related papers are the work of
\cite{bhattacharyya:bickel:2015} and of \cite{klusowski:wu:2018}. We
now discuss the relation between these papers and our results.

In \cite{bhattacharyya:bickel:2015}, the authors consider nonparametric
estimation of certain functionals of the graphon in a \textit{sparse graphon}
model, with applications in hypotheses testing and networks comparison. The
sparse graphon model and the graphex model are different, though they are
tightly related, as explained for instance in
\cite{borgs:chayes:cohn:holden:2016}. The functionals considered in
\cite{bhattacharyya:bickel:2015,bickel:chen:levina:2011} are related to the
expectation of counts of motifs too. We note, however, that they use a slightly
different notion of counts than $\nsub(R,G_t)$, though this do not fundamentally
make difference and their procedures can be easily adapted to
$\nsub(R,G_t)$. The paper \cite{bhattacharyya:bickel:2015} consider two
subsampling procedures for counting motifs. Of particular interest for us is
their \textit{uniform subsampling bootstrap}. Given a sample $G_n$ with $n$
vertices, they consider the subgraph $\tilde{G}_{n,m}$ which consists on
uniformly sampling $m < n$ vertices without replacement of $G_n$ and returning
the induced subgraph. For moderate $m$, this is morally equivalent to
$p$-sampling. Sampling vertices without replacement is quite
natural for \textit{dense} graphon sample, but less evident in the case of a
sparse graphon sample. The issue is for instance discussed in
\cite{lunde:sarkar:2019}: in the classical dense graphon model
$\tilde{G}_{n,m} \equaldist G_m$, which is not necessarily true in the sparse
model. This is not an issue to estimate $\nsub(R,G_n)$, but it is if one wants
to estimate $\var(\nsub(R,G_n))$, to do hypothesis testing for instance. Because
the distribution of $\tilde{G}_{n,m}$ has no distinguished properties, it is
unclear that $\var(\nsub(R,G_n))$ can be estimated from
$\var(\nsub(R,\tilde{G}_{n,m})\mid G_n)$; so \cite{bhattacharyya:bickel:2015}
have to rely on rather complicated procedures to estimate
$\var(\indsub(R,G_n))$.

In \cite{klusowski:wu:2018}, the authors are also interested in estimating the
counts of certain motifs using vertex sampling. They also consider two
subsampling schemes. In particular, they consider uniform subsampling of
vertices, which is almost $p$-sampling. The only difference is that they keep
isolated vertices in the subsample while we don't. Since we only consider
connected motifs, the two methods are the same, as isolated vertices cannot
contribute to the motif. Consequently, our results in
\cref{sec:p-sampling-as-bs} are partial restatements of
\cite{klusowski:wu:2018}, though they also consider a different notion of
counts. We note that in \cite{klusowski:wu:2018} the authors are only interested
in counting motifs. They don't make probabilistic assumption on the graph and
hence do not consider inference on the distribution of the graph.

\section{Some applications}
\label{sec:some-appl-exampl}

\subsection{Assessing the pertinence of sparse graphex models on some real networks}
\label{sec:assess-pert-exch}

The tail-index of $\mu$ is a distinguished feature of sparse graphex processes
(\textit{i.e.} graphex processes satisfying
\cref{ass:finiteness,ass:v1:2,ass:v1:5} for $\sigma \ne 0$). As such, it can be
used to diagnose the pertinence of using such a model on some real datasets. In
particular, if $G_t$ is a sample of a sparse graphex process, then we expect the
degree distribution of $G_t$ to be close to a power-law with index $1+\sigma$,
and $p \mapsto \log \frac{N_1(G_t)}{N_p(G_t)}$ to be closely linear with slope
$1 + \sigma$, at least for $p$ bounded away from $0$ (the asymptotic being true
only in this region). Equivalently, we might look at the plot of
$p \mapsto \hat{\sigma}_p(G_t)$, which we call a Hill plot in analogy with the
classical Hill plots used in extreme value theory to calibrate the Hill
estimator \cite{drees:haan:resnick:2000}. If the model is correct, we expect the
Hill plot to be constant, except perhaps at the boundary close to $p = 0$. The
Hill plot is doubly interesting as it allows for both assessing the pertinence
of the model, but also to choose a reasonable value of $p$ for the estimation of
the tail-index (see also \cref{rmk:2}).

We choose to illustrate the diagnostic on two datasets from
\cite{mislove:marcon:gummadi:2007}, the Flickr and the Youtube datasets. These
datasets are available from the KONECT database \cite{kunegis:2013}.  These are
social networks datasets where vertices are users and edges represent
friendships. We note that the Flickr dataset is a directed graph, but we ignore
arrows and see it as an undirected graph. We summarize in \cref{tab:1} the main
statistics of these networks.
\begin{table}[h]
  \label{tab:1}
  \centering
  \caption{Summary of main statistics for the Flickr and Youtube networks.}
  \begin{tabular}{lcc}
%    \toprule
    Network & Number of vertices & Number of edges\\
%    \midrule
    Flickr & 2302926 & 22838277\\
    Youtube & 3223586 & 9375375\\
%    \bottomrule
  \end{tabular}
\end{table}
These two networks are interesting in the way they
have been sampled. Indeed, \cite{mislove:marcon:gummadi:2007} crawled these
networks by searching for the whole connected component reachable from a random
set of users. The methodology may be summarized as starting from a single user,
and in each step, retrieve the list of friends for a not yet visited user and
add these users to a list of users to visit. The process is continued until
exhaustion of the list. This is known as \textit{breadth-first search} (BFS)
sampling. Interestingly, \cite{mislove:marcon:gummadi:2007} crawled the Flickr
and Youtube networks using BFS sampling every day on a large period of time
(several weeks). In particular, each day they revisited all the users they had
previously discovered, and in addition all the new users that were reachable
from the previously known users. By choosing a reasonable sets of users to start
with, they have access to the whole largest connected component of these
networks. Further, provided the BFS sampling goes well, their crawling mechanism
may be viewed as observing snapshots of a graphex process at increasing times.

\begin{figure}[h]
  \centering

  \includegraphics[width=\linewidth]{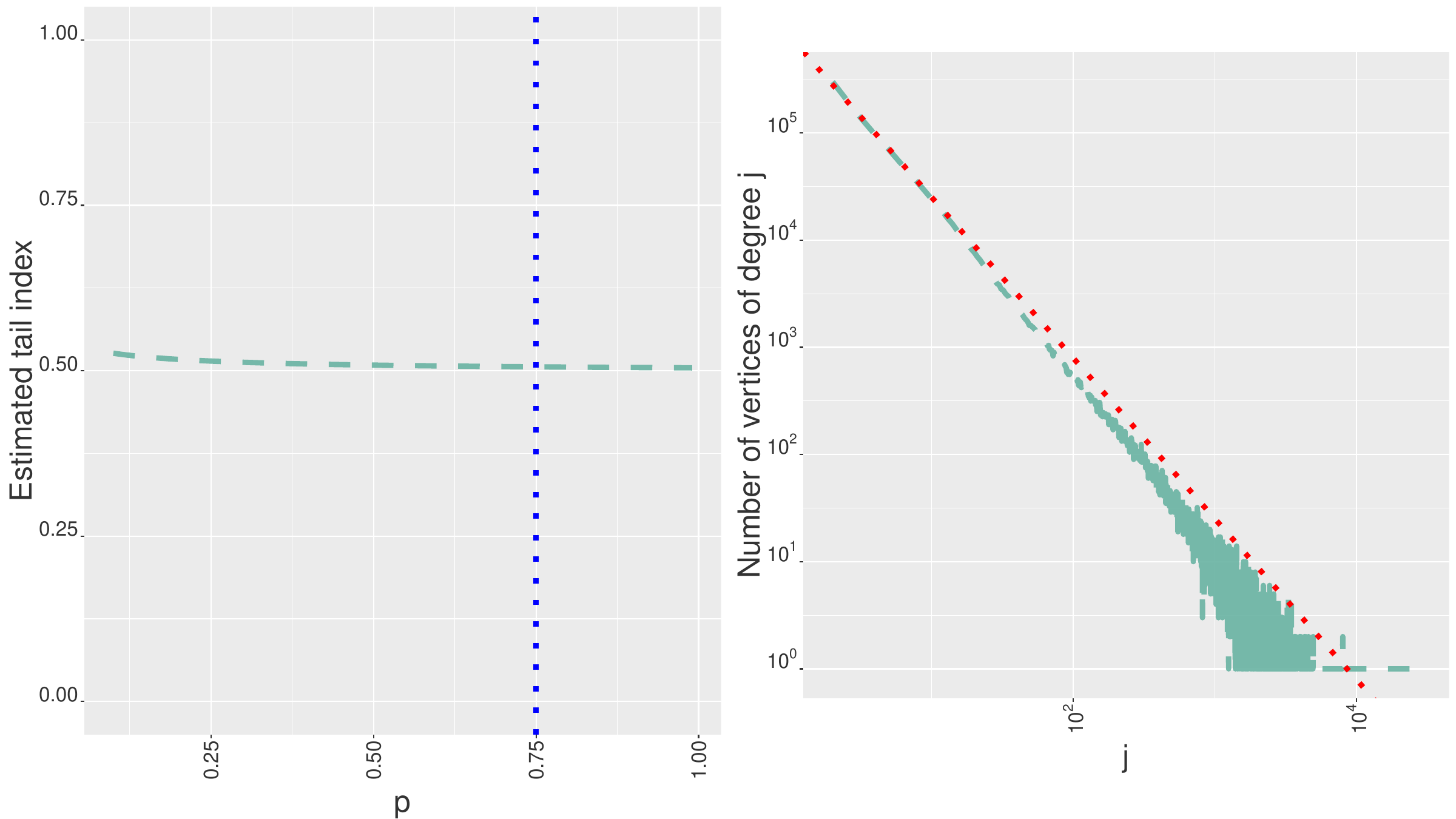}
  
  \caption{%
    On the left: The Hill plot of $\hat{\sigma}_p$ estimated on the Flickr
    dataset, that is the plot of $p\mapsto \hat{\sigma}_p$. The vertical blue
    line at $p=0.75$ corresponds to the retained value for the estimation of
    $\sigma$. We notice that the Hill plot is essentially flat in this case,
    which is the expected behaviour when the model is well-specified. On the
    right: In green is the degree distribution of the Flickr dataset. In red is
    the line corresponding to $\degree_j \propto j^{-1-\hat{\sigma}_p}$. We see
    that both curves align well.
    }
  \label{fig:6}
\end{figure}

\begin{figure}[h]
  \centering

  \includegraphics[width=\linewidth]{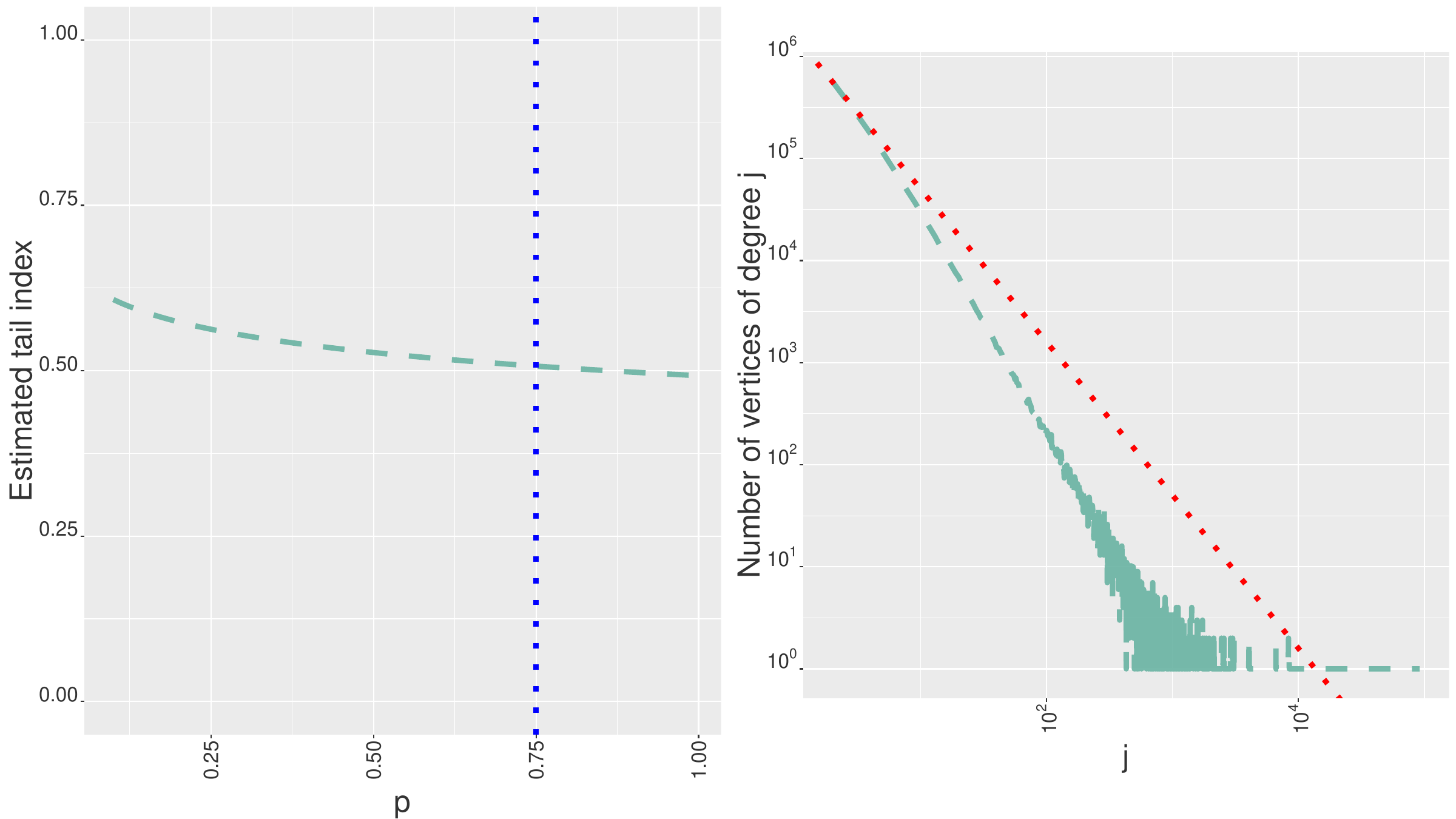}
  
  \caption{%
    On the left: The Hill plot of $\hat{\sigma}_p$ estimated on the Youtube
    dataset, that is the plot of $p\mapsto \hat{\sigma}_p$. The vertical blue
    line at $p=0.75$ corresponds to the retained value for the estimation of
    $\sigma$. We notice that the Hill plot is not flat, which should be
    considered as a warning that the model might be inadequate. On the right: In
    green is the degree distribution of the Youtube dataset. In red is the line
    corresponding to $\degree_j \propto j^{-1-\hat{\sigma}_p}$.}
  \label{fig:7}
\end{figure}

We plotted side-by-side the Hill plot and the degree distribution, respectively
for the Flickr and Youtube networks, in \cref{fig:6,fig:7}.
We plotted the degree distributions on a log-log scale, and added in red dotted
lines a linear curve with slope $-(1+\hat{\sigma}_p)$, where $\hat{\sigma}_p$
has been computed using $p=0.75$. In case the model is correct, we then expect
the degree distribution and the linear curve to align well. We can see in
\cref{fig:3} that the Flickr dataset exhibits some features of sparse graphex
processes: its Hill plot is relatively flat with value $\approx 0.5$, and the
slope of its degree distribution on a log-log scale is close to $\approx -
1.5$. Regarding the Youtube dataset, in \cref{fig:7}, this is less evident, even
the plots plea in favor of a misspecification of the sparse graphex model. This
might be because the model is not reasonable, which is likely to be caused by
the crawling procedure used in \cite{mislove:marcon:gummadi:2007}: the data they
have at their disposal make uncertain that all users in Youtube have effectively
been sampled \cite[Section 4.2.4]{mislove:marcon:gummadi:2007}\footnote{To be
  complete, the same is true for the Flickr network, but in the case of Flickr,
  the authors of \cite{mislove:marcon:gummadi:2007} were able to quantify the
  fraction of missed users \cite[Section
  4.2.1]{mislove:marcon:gummadi:2007}. They found that this fraction is likely
  to be very small. In the case of Youtube, they were unable to quantify the
  fraction of missing users.}.

Finally, we mention that these plots should only be used as visual diagnostics,
and it is hard to draw firm conclusions about the validity of the model from
these plots. One would need in particular to do a true statistical test for the
degree distribution, yet this is unclear how to perform at this time.

\subsection{Counting motifs}
\label{sec:counting-motifs}

\subsubsection{Consistent estimation of count statistics}
\label{sec:cons-estim-count}

\begin{figure}[h]
  \centering

  \includegraphics[width=\linewidth]{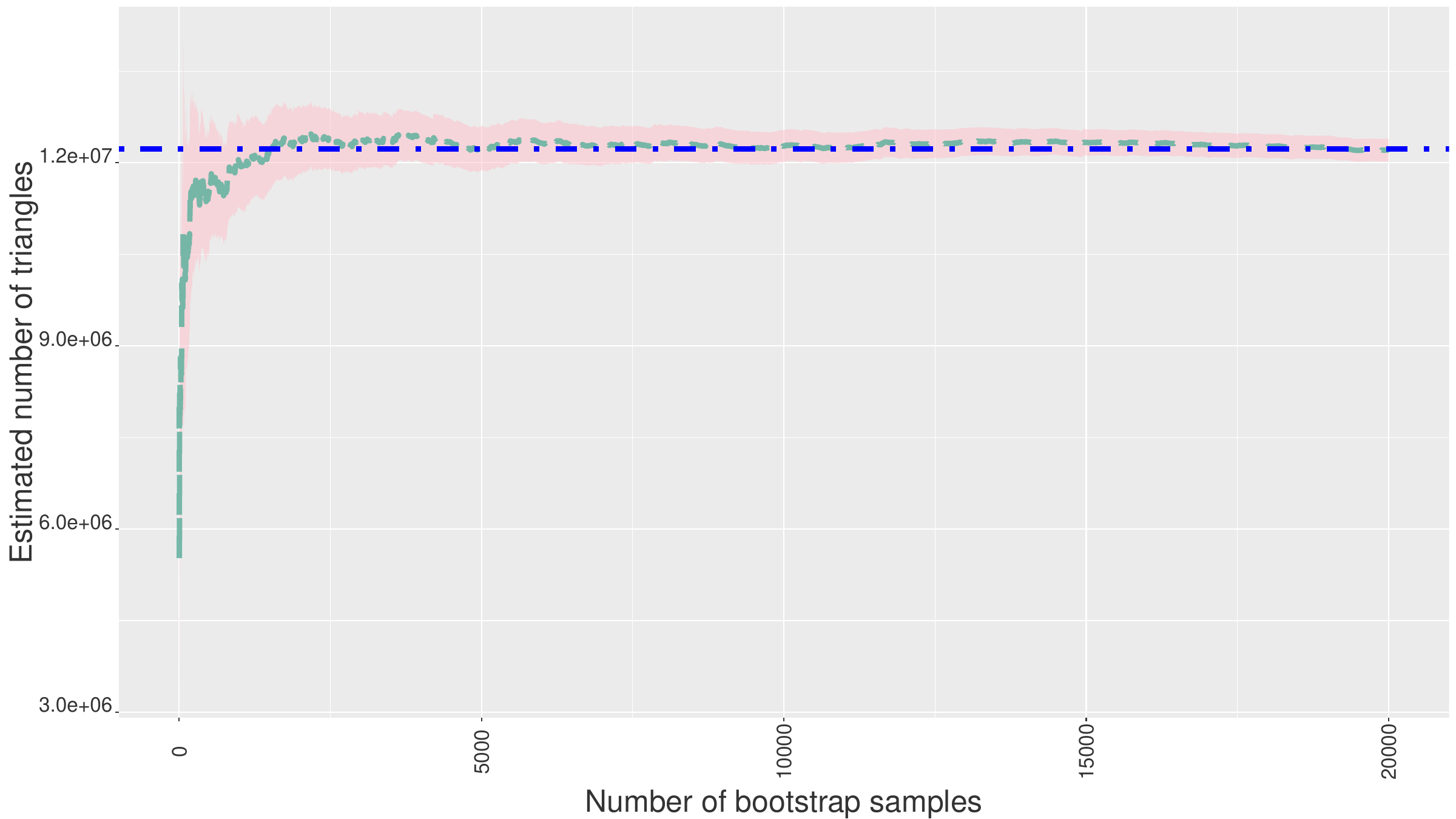}
  
  \caption{%
    In green is the plot of (one realization of) the estimated number of
    triangles using $p$-sampling of the Youtube dataset, as a function of the
    number of bootstrap sample used ($p=0.01$). The horizontal blue line
    corresponds to the true number of triangles. The pink region is
    $[\hat{\nsub}(R,G) \pm 2\delta]$, where $\delta^2$ is the empirical estimate
    of $\var(\hat{\nsub}(R,G) \mid G)$.}
  \label{fig:3}
\end{figure}

We can also use the results of \cref{sec:p-sampling-as-bs} to count the number
of motifs in a network. We illustrate this on the number of triangles in the
Youtube network of the previous section. We note here that the number of
triangles is equal to $\nsub(R,G_t)/6$ where $R$ is a triangle graph. We also
note that the method is consistent regardless of the distribution of $G_t$,
whence it does not matter that we found evidence that the Youtube dataset is
unlikely to be a sample of a sparse graphex process. We plot in \cref{fig:3} the
estimated number of triangles and the estimated standard deviation as a function
of the number of bootstrap samples used. The bootstrap samples used are small
($p=0.01$), implying that in average each sample has $\approx 12$ triangles,
which is fast to count. Yet, with only 2000 samples (which can all be computed
in parallel) the accuracy is already quite good.

\subsubsection{Testing for count statistics}
\label{sec:test-count-stat}

\begin{figure}[h]
  \centering

  \includegraphics[width=\linewidth]{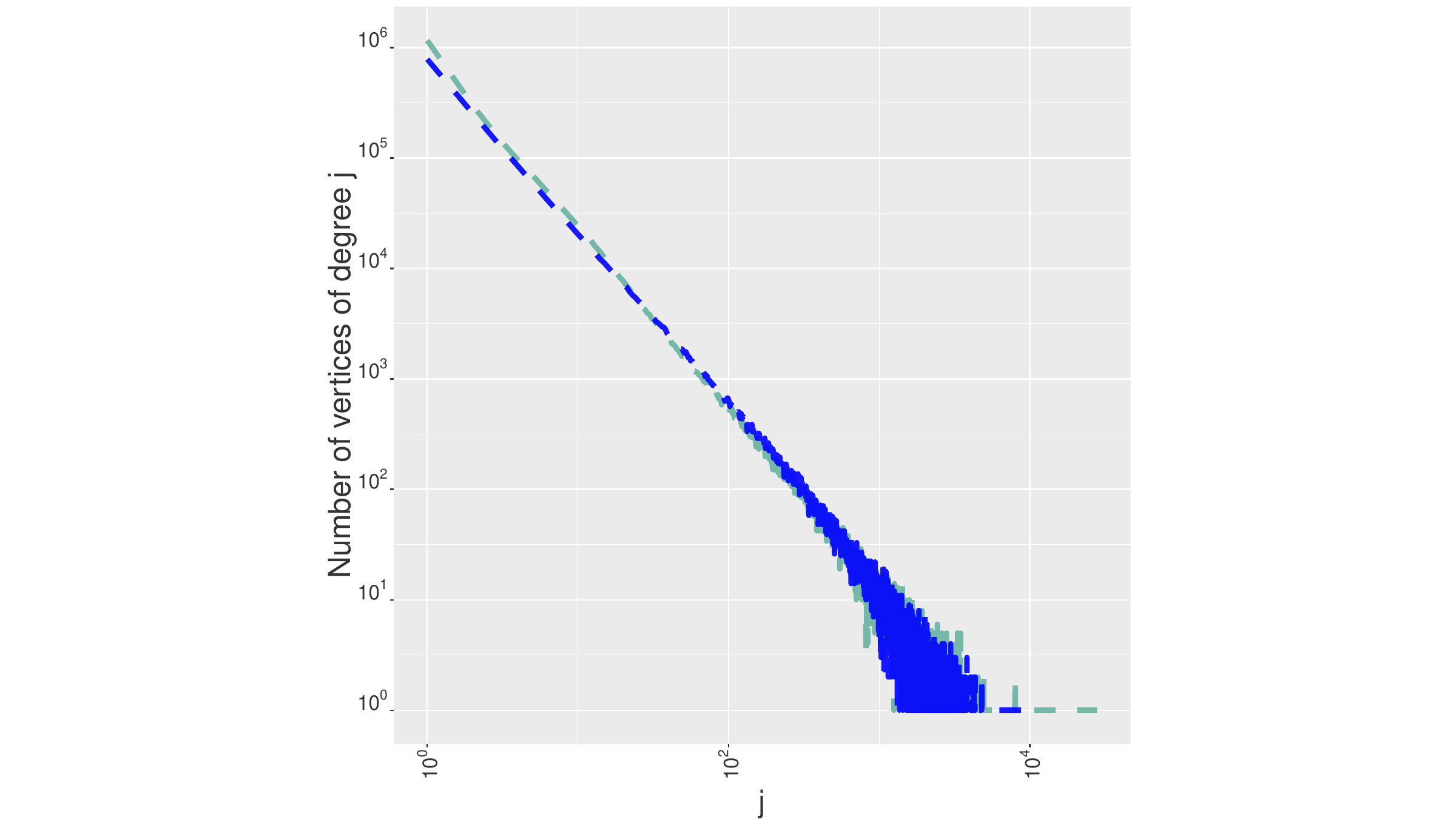}
  
  \caption{%
    In green is the degree distribution of the Flickr datasat. In blue, the
    degree distribution of one sample of a GGP model with parameters
    $(\sigma,\tau,t) = (0.5,6.3,12000)$.
  }
  \label{fig:5}
\end{figure}

We have seen in \cref{sec:assess-pert-exch} that the degree distribution of the
Flickr network was plausibly well explained as a sample from a graphex
process. In particular it is found that $\sigma = 0.5$ does explain quite well
the ``sparsity'' and the degree distribution. The current state if the art
algorithm for inference on graphex processes is the algorithm of
\cite{caron:fox:2017} and its subsequent generalization in
\cite{todeschini:miscouridou:caron:2016}. Here we focus on the simpler version
of \citet{caron:fox:2017}. Caron and Fox's algorithm is based on the GGP model
previously described in \cref{ex:v1:5}. Hence the model is a parametric family
with parameter set
$\Set{W \equiv W_{\sigma,\tau} \given \sigma < 1,\,\tau > 0}$. One of the
desirable feature of the GGP model is that the parameter $\sigma$ is exactly the
tail-index of $\mu$ described in \cref{sec:graph-marg-tail}\footnote{To be
  rigorous, this is true for $\sigma \in [0,1)$, but the GGP model also allows
  for $\sigma <0$. Then the parameter cannot be interpreted as the tail-index of
  $\mu$ anymore (which is zero).}. Yet simple, this model is able to accommodate
for $\sigma \in [0,1)$, that is from dense to sparse with power-law. The other
parameter $\tau$ acts as a cut-off for the degree distribution, controlling the
separation between the power-law regime and non power-law regime (corresponding
to larger degrees), see \cite{caron:fox:2017}. In particular, the GGP model with
$(\sigma,\tau,t) = (0.5,6.3,12000)$ satisfies $\EE_W[\nver(G_t)] \approx 2.3$
millions and $\EE_W[\nedg(G_t)] \approx 23$ millions. So, the degree
distribution and the size of Flickr are well explained by the GGP model with
$(\sigma,\tau,t) = (0.5,6.3,12000)$, as emphasized by \cref{tab:1,fig:5}. Due to
its relative simplicity, this model enables for efficient and tractable
computations \cite{caron:fox:2017,todeschini:miscouridou:caron:2016}. Then, it
makes sense to determine if other aspects of the structure of Flickr can be as
well explained by the same model. Here we test if the number of triangles
observed in Flickr is compatible with the GGP model with
$H_0 : (\sigma_0,\tau_0,t_0) = (0.5,6.3,12000)$. Under $H_0$, by numerical
integration we find that $\bar{\nsub}(R,W_0,t_0) = 3.75\cdot 10^8$ when $R$ is a
triangle. Using the bootstrap procedure of \cref{sec:bootstrap-motifs} on the
Flickr dataset (we used $p=0.1$ and $2000$ bootstrap samples), we estimated
$\nsub(R,G_t) = 5.03\cdot 10^9$ and $\sqrt{\hat{s}^2(G_t)} = 1.93\cdot
10^8$. This gives an observed value of the test statistic $T(R,G_t)$ of
approximately $24.2$, and then a $p$-value $\approx 0$. Hence $H_0$ is rejected
in this case: without surprise it does not seem reasonable that the Flickr
dataset is well-explained by a model as simple as the GGP model. This emphasize
the need of more advanced algorithm for inference, in particular for
non-parametric estimation.

\subsection{Test-train split}
\label{sec:test-train-split}

Model evaluation is a key component of data analysis.  Often, this involves
randomly splitting the available data into a test set and a training set.  There
are many seemingly natural ways to partition a graph, so some care is required
in splitting the data. The choice of partitioning scheme may induce a
significant sampling bias in the test and training sets, impeding evaluation and
complicating model comparison. Consider for instance the case of independent and
identically distributed observations $(Y_1,\dots,Y_n) \distiid P^n$: it is
natural to split the data $(Y_1,\dots,Y_n)$ by subsampling uniformly without
replacement $0 \leq m \leq n$ observations within the sample. This is because
the obtained subsample has the desirable property to be distributed according to
$P^m$. We argue that in the context of graphex processes, the natural way to
split the graph is precisely the $p$-sampling procedure described in
\cref{sec:subs-graph-proc}, because it is the only procedure leaving (up to
rescaling) the distribution of the graph invariant
\citep{veitch:roy:2019,borgs:chayes:cohn:veitch:2017}.

\begin{figure}[h]
  \centering

  \includegraphics[width=\linewidth]{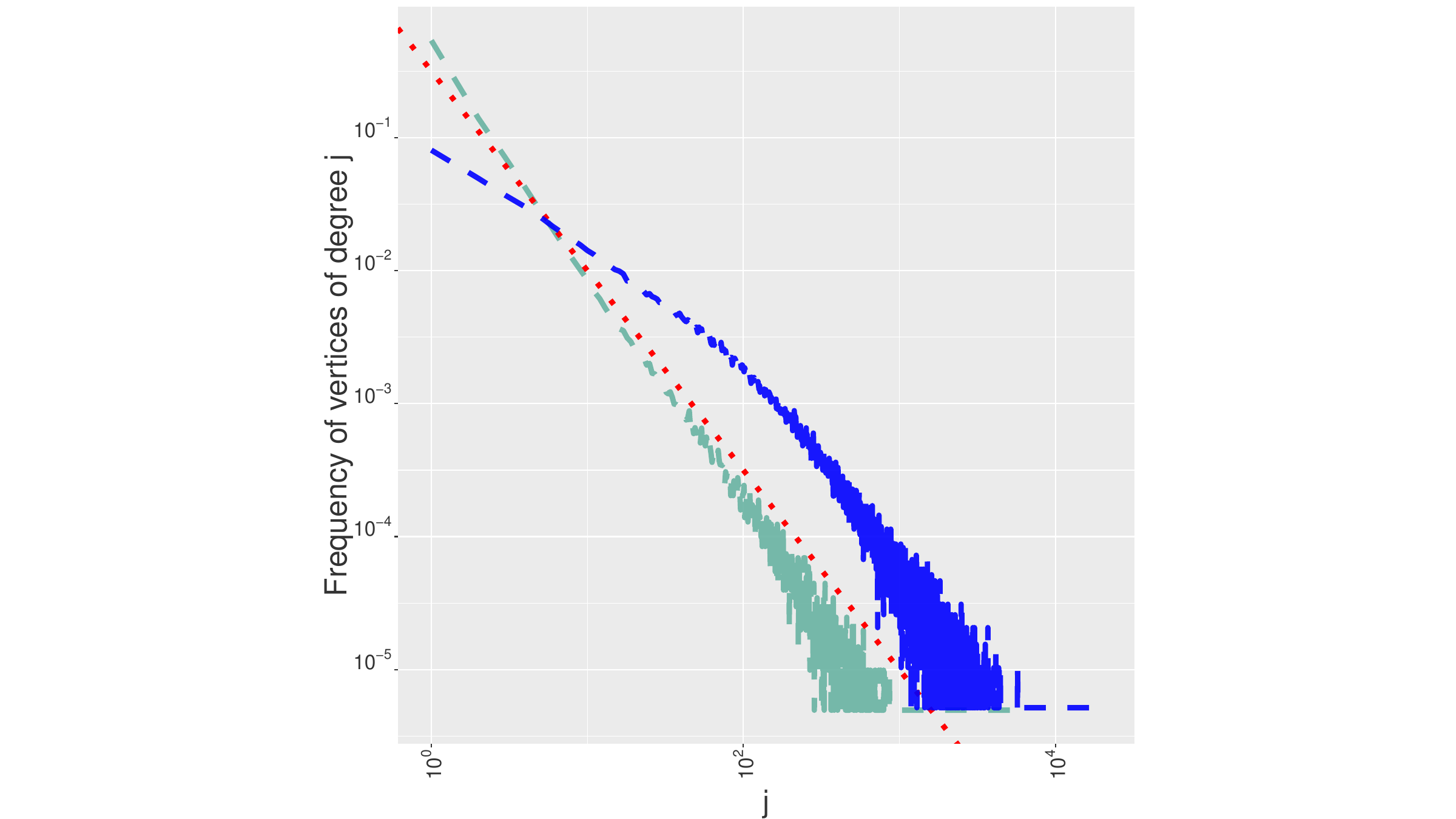}
  
  \caption{In green: the degree distribution of a $p$-sampling of the Flickr
    dataset with $p = 0.2$ (the sample has 201633 vertices and 952723 edges). In
    blue: the degree distribution of a star sampling of the Flickr dataset with
    $p' = 0.013\dots$ (the sample has 193075 vertices and 15739418 edges). In
    red is $\degree_j \propto j^{-1.5}$ which was found to be a credible
    estimate of the degree distribution of the complete Flickr dataset. }
  \label{fig:5b}
\end{figure}

To illustrate the difficulty with test--train splitting of graphs, we built two
training sets from the Flickr dataset: one using the $p$-sampling procedure, and
the other using the so-called \textit{star sampling}
\cite{capobianco:1972,kolaczyk:2009}. The star sampling we used consists on
selecting a set of vertices by uniform sampling with probability $p'$ (the stars
centers), and for each of these center select all of its neighbors (the star
points), and then returning the induced subgraph. We used $p = 0.2$ and
$p' = 0.013\dots$ to obtain training sets with a comparable number of vertices
(respectively 201633 for $p$-sampling and 193075 for star sampling). However,
the structure of the obtained training sets are very different. The $p$-sample
training set has 952723 edges, while the star sample has 15739418 edges. Whence
the star training set is much denser, as expected. The degree distributions of
the two training sets are also quite different, as illustrated in
\cref{fig:5b}. On the same figure, we drew a line with slope $-1.5$, which was
found to explain well the degree distribution of the complete dataset. We see
that it adjusts also quite well the degree distribution of the $p$-sampled
training-set. So at least in term of degree distribution, the $p$-sampling
training set captures better the structure of the complete dataset.

\section*{Acknowledgments}
\label{sec:acknowledgements}

We are thankful to François Caron for feedback on earlier versions of this
manuscript and for pointing us many interesting directions and references,
leading to a substantially improved document. We also thank the anonymous referees
for helpful comments that have also substantially improved the manuscript.  This
work was supported by U.S. Air Force Office of Scientific Research grant
\#FA9550-15-1-0074.

% Print biblio
% ------------------------------------------------------------------------------
\bibliographystyle{imsart-nameyear}
\bibliography{estimator.bib}

\begin{thebibliography}{56}
% BibTex style file: imsart-nameyear.bst, 2017-11-03
% Default style options (sort=1,type=nameyear).
% Used options (sort=1,type=nameyear).

\bibitem[\protect\citeauthoryear{Achlioptas
  et~al.}{2009}]{achlioptas:clauset:kempe:2009}
\begin{barticle}[author]
\bauthor{\bsnm{Achlioptas},~\bfnm{Dimitris}\binits{D.}},
  \bauthor{\bsnm{Clauset},~\bfnm{Aaron}\binits{A.}},
  \bauthor{\bsnm{Kempe},~\bfnm{David}\binits{D.}} \AND
  \bauthor{\bsnm{Moore},~\bfnm{Cristopher}\binits{C.}}
(\byear{2009}).
\btitle{On the Bias of Traceroute Sampling: Or, Power-Law Degree Distributions
  in Regular Graphs}.
\bjournal{Journal of the ACM (JACM)}
\bvolume{56}
\bpages{1--28}.
\end{barticle}
\endbibitem

\bibitem[\protect\citeauthoryear{Ali et~al.}{2016}]{ali:wegner:gaunt:2016}
\begin{barticle}[author]
\bauthor{\bsnm{Ali},~\bfnm{Waqar}\binits{W.}},
  \bauthor{\bsnm{Wegner},~\bfnm{Anatol~E}\binits{A.~E.}},
  \bauthor{\bsnm{Gaunt},~\bfnm{Robert~E}\binits{R.~E.}},
  \bauthor{\bsnm{Deane},~\bfnm{Charlotte~M}\binits{C.~M.}} \AND
  \bauthor{\bsnm{Reinert},~\bfnm{Gesine}\binits{G.}}
(\byear{2016}).
\btitle{Comparison of Large Networks With Sub-Sampling Strategies}.
\bjournal{Scientific reports}
\bvolume{6}
\bpages{28955}.
\end{barticle}
\endbibitem

\bibitem[\protect\citeauthoryear{Assadi, Kapralov and
  Khanna}{2018}]{assadi:kapralov:khanna:2018}
\begin{barticle}[author]
\bauthor{\bsnm{Assadi},~\bfnm{Sepehr}\binits{S.}},
  \bauthor{\bsnm{Kapralov},~\bfnm{Michael}\binits{M.}} \AND
  \bauthor{\bsnm{Khanna},~\bfnm{Sanjeev}\binits{S.}}
(\byear{2018}).
\btitle{A Simple Sublinear-Time Algorithm for Counting Arbitrary Subgraphs Via
  Edge Sampling}.
\bjournal{arXiv preprint arXiv:1811.07780}.
\end{barticle}
\endbibitem

\bibitem[\protect\citeauthoryear{Barab\'{a}si and
  Albert}{1999}]{barabasi:albert:1999}
\begin{barticle}[author]
\bauthor{\bsnm{Barab\'{a}si},~\bfnm{Albert-L\'{a}szl\'{o}}\binits{A.-L.}} \AND
  \bauthor{\bsnm{Albert},~\bfnm{R{\'e}ka}\binits{R.}}
(\byear{1999}).
\btitle{Emergence of Scaling in Random Networks}.
\bjournal{Science}
\bvolume{286}
\bpages{509-512}.
\bdoi{10.1126/science.286.5439.509}
\end{barticle}
\endbibitem

\bibitem[\protect\citeauthoryear{Bhattacharyya and
  Bickel}{2015}]{bhattacharyya:bickel:2015}
\begin{barticle}[author]
\bauthor{\bsnm{Bhattacharyya},~\bfnm{Sharmodeep}\binits{S.}} \AND
  \bauthor{\bsnm{Bickel},~\bfnm{Peter~J}\binits{P.~J.}}
(\byear{2015}).
\btitle{Subsampling Bootstrap of Count Features of Networks}.
\bjournal{Annals of Statistics}
\bvolume{43}
\bpages{2384--2411}.
\end{barticle}
\endbibitem

\bibitem[\protect\citeauthoryear{Bickel and Chen}{2009}]{bickel:chen:2009}
\begin{barticle}[author]
\bauthor{\bsnm{Bickel},~\bfnm{Peter~J.}\binits{P.~J.}} \AND
  \bauthor{\bsnm{Chen},~\bfnm{Aiyou}\binits{A.}}
(\byear{2009}).
\btitle{A Nonparametric View of Network Models and Newman-Girvan and Other
  Modularities}.
\bjournal{Proceedings of the National Academy of Sciences}
\bvolume{106}
\bpages{21068-21073}.
\bdoi{10.1073/pnas.0907096106}
\end{barticle}
\endbibitem

\bibitem[\protect\citeauthoryear{Bickel, Chen and
  Levina}{2011}]{bickel:chen:levina:2011}
\begin{barticle}[author]
\bauthor{\bsnm{Bickel},~\bfnm{Peter~J}\binits{P.~J.}},
  \bauthor{\bsnm{Chen},~\bfnm{Aiyou}\binits{A.}} \AND
  \bauthor{\bsnm{Levina},~\bfnm{Elizaveta}\binits{E.}}
(\byear{2011}).
\btitle{The Method of Moments and Degree Distributions for Network Models}.
\bjournal{Annals of Statistics}
\bvolume{39}
\bpages{2280--2301}.
\end{barticle}
\endbibitem

\bibitem[\protect\citeauthoryear{{Borgs}
  et~al.}{2016}]{borgs:chayes:cohn:holden:2016}
\begin{barticle}[author]
\bauthor{\bsnm{{Borgs}},~\bfnm{C.}\binits{C.}},
  \bauthor{\bsnm{{Chayes}},~\bfnm{J.~T.}\binits{J.~T.}},
  \bauthor{\bsnm{{Cohn}},~\bfnm{H.}\binits{H.}} \AND
  \bauthor{\bsnm{{Holden}},~\bfnm{N.}\binits{N.}}
(\byear{2016}).
\btitle{{Sparse Exchangeable Graphs and Their Limits Via Graphon processes}}.
\bjournal{ArXiv e-prints}.
\end{barticle}
\endbibitem

\bibitem[\protect\citeauthoryear{Borgs
  et~al.}{2017}]{borgs:chayes:cohn:veitch:2017}
\begin{bmisc}[author]
\bauthor{\bsnm{Borgs},~\bfnm{C.}\binits{C.}},
  \bauthor{\bsnm{Chayes},~\bfnm{J.}\binits{J.}},
  \bauthor{\bsnm{Cohn},~\bfnm{H.}\binits{H.}} \AND
  \bauthor{\bsnm{Veitch},~\bfnm{V.}\binits{V.}}
(\byear{2017}).
\btitle{Sampling perspectives on sparse exchangeable graphs}.
\bhowpublished{Preprint}.
\end{bmisc}
\endbibitem

\bibitem[\protect\citeauthoryear{Broido and
  Clauset}{2019}]{broido:clauset:2019}
\begin{barticle}[author]
\bauthor{\bsnm{Broido},~\bfnm{Anna~D}\binits{A.~D.}} \AND
  \bauthor{\bsnm{Clauset},~\bfnm{Aaron}\binits{A.}}
(\byear{2019}).
\btitle{Scale-Free Networks Are Rare}.
\bjournal{Nature communications}
\bvolume{10}
\bpages{1--10}.
\end{barticle}
\endbibitem

\bibitem[\protect\citeauthoryear{Capobianco}{1972}]{capobianco:1972}
\begin{barticle}[author]
\bauthor{\bsnm{Capobianco},~\bfnm{Michael}\binits{M.}}
(\byear{1972}).
\btitle{Estimating the Connectivity of a Graph}.
\bjournal{Lecture Notes in Mathematics}
\bpages{65-74}.
\bdoi{10.1007/bfb0067358}
\end{barticle}
\endbibitem

\bibitem[\protect\citeauthoryear{Caron}{2012}]{caron:2012}
\begin{bincollection}[author]
\bauthor{\bsnm{Caron},~\bfnm{Francois}\binits{F.}}
(\byear{2012}).
\btitle{Bayesian nonparametric models for bipartite graphs}.
In \bbooktitle{Advances in Neural Information Processing Systems 25}
(\beditor{\bfnm{F.}\binits{F.}~\bsnm{Pereira}},
  \beditor{\bfnm{C.~J.~C.}\binits{C.~J.~C.}~\bsnm{Burges}},
  \beditor{\bfnm{L.}\binits{L.}~\bsnm{Bottou}} \AND
  \beditor{\bfnm{K.~Q.}\binits{K.~Q.}~\bsnm{Weinberger}}, eds.)
\bpages{2051--2059}.
\bpublisher{Curran Associates, Inc.}
\end{bincollection}
\endbibitem

\bibitem[\protect\citeauthoryear{Caron and Fox}{2017}]{caron:fox:2017}
\begin{barticle}[author]
\bauthor{\bsnm{Caron},~\bfnm{Fran{\c{c}}ois}\binits{F.}} \AND
  \bauthor{\bsnm{Fox},~\bfnm{Emily~B.}\binits{E.~B.}}
(\byear{2017}).
\btitle{Sparse Graphs Using Exchangeable Random Measures}.
\bjournal{J. Royal Statist. Society Series B}
\bvolume{79}
\bpages{1295--1366}.
\bdoi{10.1111/rssb.12233}
\end{barticle}
\endbibitem

\bibitem[\protect\citeauthoryear{Caron, Panero and
  Rousseau}{2017}]{caron:panero:rousseau:2017}
\begin{barticle}[author]
\bauthor{\bsnm{Caron},~\bfnm{F.}\binits{F.}},
  \bauthor{\bsnm{Panero},~\bfnm{F.}\binits{F.}} \AND
  \bauthor{\bsnm{Rousseau},~\bfnm{J.}\binits{J.}}
(\byear{2017}).
\btitle{On Sparsity, Power-Law and Clustering Properties of Graphex Processes}.
\bjournal{arXiv preprint arXiv:1708.03120}.
\end{barticle}
\endbibitem

\bibitem[\protect\citeauthoryear{Carrington, Scott and
  Wasserman}{2005}]{carrington:scott:wasserman:2005}
\begin{bbook}[author]
\bauthor{\bsnm{Carrington},~\bfnm{Peter~J}\binits{P.~J.}},
  \bauthor{\bsnm{Scott},~\bfnm{John}\binits{J.}} \AND
  \bauthor{\bsnm{Wasserman},~\bfnm{Stanley}\binits{S.}}
(\byear{2005}).
\btitle{Models and methods in social network analysis}
\bvolume{28}.
\bpublisher{Cambridge university press}.
\end{bbook}
\endbibitem

\bibitem[\protect\citeauthoryear{Clauset, Shalizi and
  Newman}{2009}]{clauset:shalizi:newman:2009}
\begin{barticle}[author]
\bauthor{\bsnm{Clauset},~\bfnm{Aaron}\binits{A.}},
  \bauthor{\bsnm{Shalizi},~\bfnm{Cosma~Rohilla}\binits{C.~R.}} \AND
  \bauthor{\bsnm{Newman},~\bfnm{Mark~EJ}\binits{M.~E.}}
(\byear{2009}).
\btitle{Power-Law Distributions in Empirical Data}.
\bjournal{SIAM review}
\bvolume{51}
\bpages{661--703}.
\end{barticle}
\endbibitem

\bibitem[\protect\citeauthoryear{Crane and Dempsey}{2016}]{crane:dempsey:2016}
\begin{barticle}[author]
\bauthor{\bsnm{Crane},~\bfnm{H.}\binits{H.}} \AND
  \bauthor{\bsnm{Dempsey},~\bfnm{W.}\binits{W.}}
(\byear{2016}).
\btitle{Edge Exchangeable Models for Network Data}.
\bjournal{arXiv preprint arXiv:1603.04571}.
\end{barticle}
\endbibitem

\bibitem[\protect\citeauthoryear{Daley and
  Vere-Jones}{2003}]{daley:vere-jones:2003}
\begin{bbook}[author]
\bauthor{\bsnm{Daley},~\bfnm{D.~J.}\binits{D.~J.}} \AND
  \bauthor{\bsnm{Vere-Jones},~\bfnm{D.}\binits{D.}}
(\byear{2003}).
\btitle{An Introduction to the Theory of Point Processes. Volume I: Elementary
  Theory and Methods},
\bedition{second} ed.
\bpublisher{Springer}.
\end{bbook}
\endbibitem

\bibitem[\protect\citeauthoryear{Daley and
  Vere-Jones}{2007}]{daley:vere-jones:2007}
\begin{bbook}[author]
\bauthor{\bsnm{Daley},~\bfnm{D.~J.}\binits{D.~J.}} \AND
  \bauthor{\bsnm{Vere-Jones},~\bfnm{D.}\binits{D.}}
(\byear{2007}).
\btitle{An introduction to the theory of point processes: volume II: general
  theory and structure}
\bvolume{2}.
\bpublisher{Springer Science \& Business Media}.
\end{bbook}
\endbibitem

\bibitem[\protect\citeauthoryear{de~Haan and
  Stadtm{\"u}ller}{1996}]{haan:stadtmueller:1996}
\begin{barticle}[author]
\bauthor{\bparticle{de} \bsnm{Haan},~\bfnm{Laurens}\binits{L.}} \AND
  \bauthor{\bsnm{Stadtm{\"u}ller},~\bfnm{Ulrich}\binits{U.}}
(\byear{1996}).
\btitle{Generalized Regular Variation of Second Order}.
\bjournal{Journal of the Australian Mathematical Society}
\bvolume{61}
\bpages{381--395}.
\end{barticle}
\endbibitem

\bibitem[\protect\citeauthoryear{Diaconis and
  Janson}{2008}]{diaconis:janson:2008}
\begin{barticle}[author]
\bauthor{\bsnm{Diaconis},~\bfnm{P.}\binits{P.}} \AND
  \bauthor{\bsnm{Janson},~\bfnm{S.}\binits{S.}}
(\byear{2008}).
\btitle{Graph Limits and Exchangeable Random Graphs}.
\bjournal{Rendiconti di Matematica e delle sue Applicazioni. Serie VII}
\bpages{33--61}.
\end{barticle}
\endbibitem

\bibitem[\protect\citeauthoryear{Drees, de~Haan and
  Resnick}{2000}]{drees:haan:resnick:2000}
\begin{barticle}[author]
\bauthor{\bsnm{Drees},~\bfnm{Holger}\binits{H.}}, \bauthor{\bparticle{de}
  \bsnm{Haan},~\bfnm{Laurens}\binits{L.}} \AND
  \bauthor{\bsnm{Resnick},~\bfnm{Sidney}\binits{S.}}
(\byear{2000}).
\btitle{How to make a Hill plot}.
\bjournal{Annals of Statistics}
\bvolume{28}
\bpages{254–274}.
\bdoi{10.1214/aos/1016120372}
\end{barticle}
\endbibitem

\bibitem[\protect\citeauthoryear{Eden et~al.}{2017}]{eden:levi:ron:2017}
\begin{barticle}[author]
\bauthor{\bsnm{Eden},~\bfnm{Talya}\binits{T.}},
  \bauthor{\bsnm{Levi},~\bfnm{Amit}\binits{A.}},
  \bauthor{\bsnm{Ron},~\bfnm{Dana}\binits{D.}} \AND
  \bauthor{\bsnm{Seshadhri},~\bfnm{C}\binits{C.}}
(\byear{2017}).
\btitle{Approximately Counting Triangles in Sublinear Time}.
\bjournal{SIAM Journal on Computing}
\bvolume{46}
\bpages{1603--1646}.
\end{barticle}
\endbibitem

\bibitem[\protect\citeauthoryear{Efron and
  Tibshirani}{1994}]{efron:tibshirani:1994}
\begin{bbook}[author]
\bauthor{\bsnm{Efron},~\bfnm{Bradley}\binits{B.}} \AND
  \bauthor{\bsnm{Tibshirani},~\bfnm{Robert~J}\binits{R.~J.}}
(\byear{1994}).
\btitle{An introduction to the bootstrap}.
\bpublisher{CRC Press}.
\end{bbook}
\endbibitem

\bibitem[\protect\citeauthoryear{Feller}{1971}]{feller:1971}
\begin{bbook}[author]
\bauthor{\bsnm{Feller},~\bfnm{William}\binits{W.}}
(\byear{1971}).
\btitle{An introduction to probability theory and its applications, volume II}
\bvolume{2}.
\bpublisher{Wiley, New York}.
\end{bbook}
\endbibitem

\bibitem[\protect\citeauthoryear{Goldreich and Ron}{2008}]{goldreich:ron:2008}
\begin{barticle}[author]
\bauthor{\bsnm{Goldreich},~\bfnm{Oded}\binits{O.}} \AND
  \bauthor{\bsnm{Ron},~\bfnm{Dana}\binits{D.}}
(\byear{2008}).
\btitle{Approximating Average Parameters of Graphs}.
\bjournal{Random Structures \& Algorithms}
\bvolume{32}
\bpages{473--493}.
\end{barticle}
\endbibitem

\bibitem[\protect\citeauthoryear{Gonen, Ron and
  Shavitt}{2011}]{gonen:ron:shavitt:2011}
\begin{barticle}[author]
\bauthor{\bsnm{Gonen},~\bfnm{Mira}\binits{M.}},
  \bauthor{\bsnm{Ron},~\bfnm{Dana}\binits{D.}} \AND
  \bauthor{\bsnm{Shavitt},~\bfnm{Yuval}\binits{Y.}}
(\byear{2011}).
\btitle{Counting Stars and Other Small Subgraphs in Sublinear-Time}.
\bjournal{SIAM Journal on Discrete Mathematics}
\bvolume{25}
\bpages{1365--1411}.
\end{barticle}
\endbibitem

\bibitem[\protect\citeauthoryear{Green and Shalizi}{2017}]{green:shalizi:2017}
\begin{barticle}[author]
\bauthor{\bsnm{Green},~\bfnm{Alden}\binits{A.}} \AND
  \bauthor{\bsnm{Shalizi},~\bfnm{Cosma~Rohilla}\binits{C.~R.}}
(\byear{2017}).
\btitle{Bootstrapping Exchangeable Random Graphs}.
\bjournal{arXiv preprint arXiv:1711.00813}.
\end{barticle}
\endbibitem

\bibitem[\protect\citeauthoryear{Herlau, Schmidt and
  M{\o}rup}{2016}]{herlau:schmidt:moerup:2016}
\begin{bincollection}[author]
\bauthor{\bsnm{Herlau},~\bfnm{Tue}\binits{T.}},
  \bauthor{\bsnm{Schmidt},~\bfnm{Mikkel~N}\binits{M.~N.}} \AND
  \bauthor{\bsnm{M{\o}rup},~\bfnm{Morten}\binits{M.}}
(\byear{2016}).
\btitle{Completely random measures for modelling block-structured sparse
  networks}.
In \bbooktitle{Advances in Neural Information Processing Systems 29}
(\beditor{\bfnm{D.~D.}\binits{D.~D.}~\bsnm{Lee}},
  \beditor{\bfnm{M.}\binits{M.}~\bsnm{Sugiyama}},
  \beditor{\bfnm{U.~V.}\binits{U.~V.}~\bsnm{Luxburg}},
  \beditor{\bfnm{I.}\binits{I.}~\bsnm{Guyon}} \AND
  \beditor{\bfnm{R.}\binits{R.}~\bsnm{Garnett}}, eds.)
\bpages{4260--4268}.
\bpublisher{Curran Associates, Inc.}
\end{bincollection}
\endbibitem

\bibitem[\protect\citeauthoryear{Holme}{2019}]{holme:2019}
\begin{barticle}[author]
\bauthor{\bsnm{Holme},~\bfnm{Petter}\binits{P.}}
(\byear{2019}).
\btitle{Rare and everywhere: Perspectives on scale-free networks}.
\bjournal{Nature communications}
\bvolume{10}
\bpages{1--3}.
\end{barticle}
\endbibitem

\bibitem[\protect\citeauthoryear{Janson}{2017}]{janson:2017:2}
\begin{barticle}[author]
\bauthor{\bsnm{Janson},~\bfnm{S.}\binits{S.}}
(\byear{2017}).
\btitle{{On convergence for graphexes}}.
\bjournal{ArXiv e-prints}.
\end{barticle}
\endbibitem

\bibitem[\protect\citeauthoryear{Jones and
  Handcock}{2003}]{jones:handcock:2003}
\begin{barticle}[author]
\bauthor{\bsnm{Jones},~\bfnm{James~Holland}\binits{J.~H.}} \AND
  \bauthor{\bsnm{Handcock},~\bfnm{Mark~S}\binits{M.~S.}}
(\byear{2003}).
\btitle{Sexual contacts and epidemic thresholds}.
\bjournal{Nature}
\bvolume{423}
\bpages{605--606}.
\end{barticle}
\endbibitem

\bibitem[\protect\citeauthoryear{Kingman}{1993}]{kingman:1993}
\begin{bbook}[author]
\bauthor{\bsnm{Kingman},~\bfnm{J.~F.~C.}\binits{J.~F.~C.}}
(\byear{1993}).
\btitle{Poisson Processes}.
\bpublisher{Oxford University Press}.
\end{bbook}
\endbibitem

\bibitem[\protect\citeauthoryear{Klusowski and Wu}{2018}]{klusowski:wu:2018}
\begin{barticle}[author]
\bauthor{\bsnm{Klusowski},~\bfnm{Jason~M}\binits{J.~M.}} \AND
  \bauthor{\bsnm{Wu},~\bfnm{Yihong}\binits{Y.}}
(\byear{2018}).
\btitle{Counting motifs with graph sampling}.
\bjournal{arXiv preprint arXiv:1802.07773}.
\end{barticle}
\endbibitem

\bibitem[\protect\citeauthoryear{Kolaczyk}{2009}]{kolaczyk:2009}
\begin{bbook}[author]
\bauthor{\bsnm{Kolaczyk},~\bfnm{E.~D.}\binits{E.~D.}}
(\byear{2009}).
\btitle{Statistical Analysis of Network Data. Methods and models.}
\bpublisher{Springer}.
\end{bbook}
\endbibitem

\bibitem[\protect\citeauthoryear{Kunegis}{2013}]{kunegis:2013}
\begin{barticle}[author]
\bauthor{\bsnm{Kunegis},~\bfnm{Jérôme}\binits{J.}}
(\byear{2013}).
\btitle{KONECT}.
\bjournal{Proceedings of the 22nd International Conference on World Wide Web -
  WWW ’13 Companion}.
\bdoi{10.1145/2487788.2488173}
\end{barticle}
\endbibitem

\bibitem[\protect\citeauthoryear{Lakhina
  et~al.}{2003}]{lakhina:byers:crovella:2003}
\begin{binproceedings}[author]
\bauthor{\bsnm{Lakhina},~\bfnm{Anukool}\binits{A.}},
  \bauthor{\bsnm{Byers},~\bfnm{John~W}\binits{J.~W.}},
  \bauthor{\bsnm{Crovella},~\bfnm{Mark}\binits{M.}} \AND
  \bauthor{\bsnm{Xie},~\bfnm{Peng}\binits{P.}}
(\byear{2003}).
\btitle{Sampling biases in IP topology measurements}.
In \bbooktitle{IEEE INFOCOM 2003. Twenty-second Annual Joint Conference of the
  IEEE Computer and Communications Societies (IEEE Cat. No. 03CH37428)}
\bvolume{1}
\bpages{332--341}.
\bpublisher{IEEE}.
\end{binproceedings}
\endbibitem

\bibitem[\protect\citeauthoryear{Leskovec, Kleinberg and
  Faloutsos}{2005}]{leskovec:kleinberg:faloutsos:2005}
\begin{binproceedings}[author]
\bauthor{\bsnm{Leskovec},~\bfnm{Jure}\binits{J.}},
  \bauthor{\bsnm{Kleinberg},~\bfnm{Jon}\binits{J.}} \AND
  \bauthor{\bsnm{Faloutsos},~\bfnm{Christos}\binits{C.}}
(\byear{2005}).
\btitle{Graphs over time: densification laws, shrinking diameters and possible
  explanations}.
In \bbooktitle{Proceedings of the eleventh ACM SIGKDD international conference
  on Knowledge discovery in data mining}
\bpages{177--187}.
\end{binproceedings}
\endbibitem

\bibitem[\protect\citeauthoryear{Levin and Levina}{2019}]{levin:levina:2019}
\begin{barticle}[author]
\bauthor{\bsnm{Levin},~\bfnm{Keith}\binits{K.}} \AND
  \bauthor{\bsnm{Levina},~\bfnm{Elizaveta}\binits{E.}}
(\byear{2019}).
\btitle{Bootstrapping networks with latent space structure}.
\bjournal{arXiv preprint arXiv:1907.10821}.
\end{barticle}
\endbibitem

\bibitem[\protect\citeauthoryear{Lov{\'a}sz}{2013}]{lovasz:2013}
\begin{bbook}[author]
\bauthor{\bsnm{Lov{\'a}sz},~\bfnm{L.}\binits{L.}}
(\byear{2013}).
\btitle{Large Networks and Graph Limits}
\bvolume{60}.
\bpublisher{American Mathematical Society}.
\end{bbook}
\endbibitem

\bibitem[\protect\citeauthoryear{Lunde and Sarkar}{2019}]{lunde:sarkar:2019}
\begin{barticle}[author]
\bauthor{\bsnm{Lunde},~\bfnm{Robert}\binits{R.}} \AND
  \bauthor{\bsnm{Sarkar},~\bfnm{Purnamrita}\binits{P.}}
(\byear{2019}).
\btitle{Subsampling sparse graphons under minimal assumptions}.
\bjournal{arXiv preprint arXiv:1907.12528}.
\end{barticle}
\endbibitem

\bibitem[\protect\citeauthoryear{Middendorf, Ziv and
  Wiggins}{2005}]{middendorf:ziv:wiggins:2005}
\begin{barticle}[author]
\bauthor{\bsnm{Middendorf},~\bfnm{Manuel}\binits{M.}},
  \bauthor{\bsnm{Ziv},~\bfnm{Etay}\binits{E.}} \AND
  \bauthor{\bsnm{Wiggins},~\bfnm{Chris~H}\binits{C.~H.}}
(\byear{2005}).
\btitle{Inferring network mechanisms: the Drosophila melanogaster protein
  interaction network}.
\bjournal{Proceedings of the National Academy of Sciences}
\bvolume{102}
\bpages{3192--3197}.
\end{barticle}
\endbibitem

\bibitem[\protect\citeauthoryear{Mislove
  et~al.}{2007}]{mislove:marcon:gummadi:2007}
\begin{barticle}[author]
\bauthor{\bsnm{Mislove},~\bfnm{Alan}\binits{A.}},
  \bauthor{\bsnm{Marcon},~\bfnm{Massimiliano}\binits{M.}},
  \bauthor{\bsnm{Gummadi},~\bfnm{Krishna~P.}\binits{K.~P.}},
  \bauthor{\bsnm{Druschel},~\bfnm{Peter}\binits{P.}} \AND
  \bauthor{\bsnm{Bhattacharjee},~\bfnm{Bobby}\binits{B.}}
(\byear{2007}).
\btitle{Measurement and Analysis of Online Social Networks}.
\bjournal{Proceedings of the 7th ACM SIGCOMM conference on Internet measurement
  - IMC '07}.
\bdoi{10.1145/1298306.1298311}
\end{barticle}
\endbibitem

\bibitem[\protect\citeauthoryear{Newman}{2005}]{newman:2005}
\begin{barticle}[author]
\bauthor{\bsnm{Newman},~\bfnm{Mark~EJ}\binits{M.~E.}}
(\byear{2005}).
\btitle{Power laws, Pareto distributions and Zipf's law}.
\bjournal{Contemporary physics}
\bvolume{46}
\bpages{323--351}.
\end{barticle}
\endbibitem

\bibitem[\protect\citeauthoryear{Orbanz and Roy}{2015}]{orbanz:roy:2015}
\begin{barticle}[author]
\bauthor{\bsnm{Orbanz},~\bfnm{Peter}\binits{P.}} \AND
  \bauthor{\bsnm{Roy},~\bfnm{Daniel~M}\binits{D.~M.}}
(\byear{2015}).
\btitle{Bayesian models of graphs, arrays and other exchangeable random
  structures}.
\bjournal{IEEE transactions on pattern analysis and machine intelligence}
\bvolume{37}
\bpages{437--461}.
\end{barticle}
\endbibitem

\bibitem[\protect\citeauthoryear{Reitzner and
  Schulte}{2013}]{reitzner:schulte:2013}
\begin{barticle}[author]
\bauthor{\bsnm{Reitzner},~\bfnm{Matthias}\binits{M.}} \AND
  \bauthor{\bsnm{Schulte},~\bfnm{Matthias}\binits{M.}}
(\byear{2013}).
\btitle{Central Limit Theorems for $U$-statistics of Poisson Point Processes}.
\bjournal{Annals of Probability}
\bvolume{41}
\bpages{3879--3909}.
\end{barticle}
\endbibitem

\bibitem[\protect\citeauthoryear{Shen-Orr
  et~al.}{2002}]{shen-orr:milo:mangan:2002}
\begin{barticle}[author]
\bauthor{\bsnm{Shen-Orr},~\bfnm{Shai~S}\binits{S.~S.}},
  \bauthor{\bsnm{Milo},~\bfnm{Ron}\binits{R.}},
  \bauthor{\bsnm{Mangan},~\bfnm{Shmoolik}\binits{S.}} \AND
  \bauthor{\bsnm{Alon},~\bfnm{Uri}\binits{U.}}
(\byear{2002}).
\btitle{Network motifs in the transcriptional regulation network of Escherichia
  coli}.
\bjournal{Nature genetics}
\bvolume{31}
\bpages{64--68}.
\end{barticle}
\endbibitem

\bibitem[\protect\citeauthoryear{Suri and
  Vassilvitskii}{2011}]{suri:vassilvitskii:2011}
\begin{binproceedings}[author]
\bauthor{\bsnm{Suri},~\bfnm{Siddharth}\binits{S.}} \AND
  \bauthor{\bsnm{Vassilvitskii},~\bfnm{Sergei}\binits{S.}}
(\byear{2011}).
\btitle{Counting triangles and the curse of the last reducer}.
In \bbooktitle{Proceedings of the 20th international conference on World wide
  web}
\bpages{607--614}.
\end{binproceedings}
\endbibitem

\bibitem[\protect\citeauthoryear{Thompson}{2012}]{thompson:2012}
\begin{bbook}[author]
\bauthor{\bsnm{Thompson},~\bfnm{Steven~K.}\binits{S.~K.}}
(\byear{2012}).
\btitle{Sampling},
\bedition{3rd} ed.
\bpublisher{Wiley}.
\end{bbook}
\endbibitem

\bibitem[\protect\citeauthoryear{Thompson and
  Frank}{2000}]{thompson:frank:2000}
\begin{barticle}[author]
\bauthor{\bsnm{Thompson},~\bfnm{Steven~K}\binits{S.~K.}} \AND
  \bauthor{\bsnm{Frank},~\bfnm{Ove}\binits{O.}}
(\byear{2000}).
\btitle{Model-based estimation with link-tracing sampling designs}.
\bjournal{Survey methodology}
\bvolume{26}
\bpages{87--98}.
\end{barticle}
\endbibitem

\bibitem[\protect\citeauthoryear{Todeschini, Miscouridou and
  Caron}{2016}]{todeschini:miscouridou:caron:2016}
\begin{barticle}[author]
\bauthor{\bsnm{Todeschini},~\bfnm{Adrien}\binits{A.}},
  \bauthor{\bsnm{Miscouridou},~\bfnm{Xenia}\binits{X.}} \AND
  \bauthor{\bsnm{Caron},~\bfnm{Francois}\binits{F.}}
(\byear{2016}).
\btitle{Exchangeable Random Measures for Sparse and Modular Graphs with
  Overlapping Communities}.
\end{barticle}
\endbibitem

\bibitem[\protect\citeauthoryear{Veitch and Roy}{2015}]{veitch:roy:2015}
\begin{barticle}[author]
\bauthor{\bsnm{Veitch},~\bfnm{Victor}\binits{V.}} \AND
  \bauthor{\bsnm{Roy},~\bfnm{Daniel~M}\binits{D.~M.}}
(\byear{2015}).
\btitle{The class of random graphs arising from exchangeable random measures}.
\bjournal{arXiv preprint arXiv:1512.03099}.
\end{barticle}
\endbibitem

\bibitem[\protect\citeauthoryear{Veitch and Roy}{2019}]{veitch:roy:2019}
\begin{barticle}[author]
\bauthor{\bsnm{Veitch},~\bfnm{V.}\binits{V.}} \AND
  \bauthor{\bsnm{Roy},~\bfnm{D.~M.}\binits{D.~M.}}
(\byear{2019}).
\btitle{Sampling and Estimation for (Sparse) Exchangeable Graphs}.
\bjournal{Annals of Statistics}
\bvolume{47}
\bpages{3274-3299}.
\end{barticle}
\endbibitem

\bibitem[\protect\citeauthoryear{Voitalov
  et~al.}{2019}]{voitalov:hoorn:hofstad:2019}
\begin{barticle}[author]
\bauthor{\bsnm{Voitalov},~\bfnm{Ivan}\binits{I.}}, \bauthor{\bparticle{van~der}
  \bsnm{Hoorn},~\bfnm{Pim}\binits{P.}}, \bauthor{\bparticle{van~der}
  \bsnm{Hofstad},~\bfnm{Remco}\binits{R.}} \AND
  \bauthor{\bsnm{Krioukov},~\bfnm{Dmitri}\binits{D.}}
(\byear{2019}).
\btitle{Scale-free networks well done}.
\bjournal{Physical Review Research}
\bvolume{1}
\bpages{033034}.
\end{barticle}
\endbibitem

\bibitem[\protect\citeauthoryear{Wang and Resnick}{2019}]{wang:resnick:2019}
\begin{barticle}[author]
\bauthor{\bsnm{Wang},~\bfnm{Tiandong}\binits{T.}} \AND
  \bauthor{\bsnm{Resnick},~\bfnm{Sidney~I}\binits{S.~I.}}
(\byear{2019}).
\btitle{Consistency of Hill estimators in a linear preferential attachment
  model}.
\bjournal{Extremes}
\bvolume{22}
\bpages{1--28}.
\end{barticle}
\endbibitem

\bibitem[\protect\citeauthoryear{Watts and
  Strogatz}{1998}]{watts:strogatz:1998}
\begin{barticle}[author]
\bauthor{\bsnm{Watts},~\bfnm{Duncan~J}\binits{D.~J.}} \AND
  \bauthor{\bsnm{Strogatz},~\bfnm{Steven~H}\binits{S.~H.}}
(\byear{1998}).
\btitle{Collective dynamics of 'small-world' networks}.
\bjournal{nature}
\bvolume{393}
\bpages{440--442}.
\end{barticle}
\endbibitem

\end{thebibliography}

\appendix

\section{Proofs related to tail-index estimation}
\label{sec:proofs-related-tail}

\subsection{Preliminaries}
\label{sec:preliminaries}

We let $\bm{B} \coloneqq (B_{i,j})_{(i,j)\in \Nats^2}$ be the symmetric array of
binary variables that indicate that $v_i$ and $v_j$ are connected if
$B_{i,j} = 1$, not connected if $B_{i,j} = 0$. All other quantities are defined
in the main document. We drop out the subscript $W$ for convenience, and
$\EE$ shall be interpreted as $\EE_W$. Also, to simplify notations, we write
$N_{1,t} \equiv N_1(G_t)$, $N_{p,t} \equiv N_p(G_t)$, and $\hat{\sigma}_{p,t}
\equiv \hat{\sigma}_p(G_t)$.

\subsection{Bias-variance decomposition}
\label{sec:v1:bias-vari-decomp}

The starting point of the proof is to decompose the risk onto a deterministic
bias term and some stochastic variance terms. Then, when $N_{p,t} = 0$ we
have
\begin{equation}
  \label{eq:v1:5}
  (\hat{\sigma}_{p,t} - \sigma)^2 = \sigma^2 \leq 1.
\end{equation}
In the situation where $N_{p,t} \geq 1$, we have
\begin{multline}
  - \log (p) \hat{\sigma}_{p,\size}
  = \log\frac{\EE N_{1,\size}}{\EE N_{p,\size}} + \log p%
  + \log\left(1 + \frac{N_{1,\size} - \EE N_{1,\size}}{\EE N_{1,\size}} \right)\\%
  - \log\left(1 + \frac{N_{p,\size} - \EE N_{p,\size}}{\EE N_{p,\size}} \right).
\end{multline}
To ease notations, we define,
\begin{gather}
  b_{\sigma,\size} \coloneqq \frac{\log(\EE N_{1,\size} / \EE
    N_{p,\size})}{-\log p} - 1
  - \sigma,\quad%
  Z_p \coloneqq \frac{N_{p,\size} - \EE N_{p,\size}}{\EE N_{p,\size}}.
\end{gather}
Note that $Z_1,Z_p > -1$ because $N_{1,\size} \geq N_{p,\size} \geq 1$. Then
when $N_{p,\size}\geq 1$,
\begin{align}
  (\hat{\sigma}_{p,t} - \sigma)^2
  &= \left( b_{\sigma,\size} + \frac{\log(1 + Z_1)}{-\log p} - \frac{\log(1 +
    Z_p)}{-\log p} \right)^2\\
  &=b_{\sigma,\size}^2 + 2 b_{\sigma,\size}\frac{\log(1 + Z_1)}{-\log p} - 2
    b_{\sigma,\size} \frac{\log(1 + Z_p)}{-\log p}\\
  &\quad + \frac{\log^2(1+Z_1)}{(-\log p)^2} + \frac{\log^2(1 + Z_p)}{(-\log
    p)^2} -2 \frac{\log(1 + Z_1)\log(1 + Z_p)}{(-\log p)^2}.
\end{align}
We now introduce functions
$\varphi_1,\varphi_2 : (-1,\infty) \rightarrow \NNReals$ such that for every $z
\ne 0$,
\begin{gather}
  \varphi_1(z) \coloneqq \frac{\log(1+z)}{z},\qquad%
  \varphi_2(z) \coloneqq -\frac{\log(1 + z) - z}{z^2}.
\end{gather}
For $z= 0$ the functions $\varphi_1$ and $\varphi_2$ are extended by continuity.
The functions $\varphi_1$ and $\varphi_2$ are non-negative and monotonically
decreasing on $(-1,\infty)$. We can therefore write when $N_{p,\size} \geq 1$,
\begin{multline}
  \label{eq:v1:6}
  (\hat{\sigma}_{p,\size} - \sigma)^{2} = b_{\sigma,\size}^2
  +\frac{Z_1^2\varphi_1(Z_1)^2}{(-\log p)^2} +
  \frac{Z_p^2\varphi_1(Z_p)^2}{(-\log p)^2} - 2\frac{Z_1Z_p\varphi_1(Z_1)
    \varphi_1(Z_p)}{(-\log p)^2}\\
  +  \frac{2b_{\sigma,\size} Z_1}{-\log p} - \frac{2b_{\sigma,\size} Z_1^2
    \varphi_2(Z_1)}{-\log p} - \frac{2b_{\sigma,\size} Z_p}{-\log p} +
  \frac{2b_{\sigma,\size}
    Z_p^2\varphi_2(Z_p)}{-\log p}.
\end{multline}
But, by Young's inequality we have $2|Z_1Z_p\varphi_1(Z_1)\varphi_1(Z_p)| \leq
Z_1^2\varphi_1(Z_1)^2 + Z_p^2\varphi_2(Z_p)^2$. Combining the Young inequality
estimate with equations \eqref{eq:v1:5} and \eqref{eq:v1:6} gives
\begin{align}
  \EE[(\hat{\sigma}_{p,\size} - \sigma)^2]
  &= \EE[(\hat{\sigma}_{p,\size} - \sigma)^2 \Ind_{N_{p,\size} = 0}] +
    \EE[(\hat{\sigma}_{p,\size} - \sigma)^2 \Ind_{N_{p,\size} \geq 1}]\\
  &\leq\Pr(N_{p,\size} = 0) + b_{\sigma,\size}^2%
    + \frac{2 \EE[Z_1^2 \varphi_1(Z_1)^2\Ind_{N_{p,\size}\geq 1}]}{(-\log p)^2}%
    + \frac{2 b_{\sigma,\size}\EE[Z_1 \Ind_{N_{p,\size} \geq 1}]}{-\log p}\\
  &\quad + \frac{2\EE[Z_p^2 \varphi_1(Z_p)^2 \Ind_{N_{p,\size} \geq
    1}]}{(-\log p)^2} - \frac{2 b_{\sigma,\size} \EE[Z_p
    \Ind_{N_{p,\size}\geq 1}]}{-\log p}\\
  &\quad + \frac{2
    b_{\sigma,\size}\EE[Z_p^2\varphi_2(Z_p)\Ind_{N_{p,\size}\geq 1}]}{-\log
    p}%
    - \frac{2 b_{\sigma,\size} \EE[Z_1^2\varphi_2(Z_1) \Ind_{N_{p,\size} \geq
    1}]}{-\log p}.
\end{align}
Moreover, because $\EE[Z_1] = \EE[Z_p] = 0$ we get by Hölder's inequality the
following estimates.
\begin{gather}
  \EE[Z_1 \Ind_{N_{p,\size} \geq 1}]
  = \EE[Z_1(1 - \Ind_{N_{p,\size}= 0})] \leq \EE[|Z_1| \Ind_{N_{p,\size} = 0}]
  \leq \sqrt{\EE[Z_1^2] \Pr(N_{p,\size} = 0)}\\
  \EE[Z_p \Ind_{N_{p,\size} \geq 1}] = \EE[Z_p(1 - \Ind_{N_{p,\size}} = 0)] =
  \Pr(N_{p,\size} = 0).
\end{gather}
Also, because $\EE[N_{p,\size}] > 0$ (see \cref{pro:v1:8} below), by Chebychev
inequality, for any $t < 1$,
\begin{equation}
  \Pr(N_{p,\size} = 0)
  \leq \Pr( N_{p,\size} < \EE N_{p,\size}(1 - t))
  \leq \Pr(Z_p < - t)
  \leq \frac{\EE[Z_p^2]}{t^2}.
\end{equation}
Since this is true for all $t < 1$, we certainly have
$\Pr(N_{p,\size} = 0) \leq \EE[Z_p^2]$.  Furthermore,
$N_{1,\size} \geq N_{p,\size} \geq 1$ implies that
$Z_1 \geq -1 + (\EE N_{1,\size})^{-1}$, then
\begin{equation}
  \varphi_1(Z_1)\Ind_{N_{p,\size} \geq 1} \leq \frac{\log \EE N_{1,\size}}{1 -
    (\EE N_{1,\size})^{-1}}.
\end{equation}
Also,
\begin{equation}
  \varphi_2(Z_1)
  = \varphi_1(Z_1) - \frac{(1 + Z_1)\log(1 + Z_1) - Z_1}{Z_1^2}
  \leq \varphi_1(Z_1).
\end{equation}
Obviously, the same estimates hold for $Z_p$. Combining all the previous
estimates yield the bound (where $(x)_+ = x$ if $x\geq 0$ and $(x)_+ = 0$ if $x
< 0$),
\begin{multline}
  \EE[(\hat{\sigma}_{p,\size} - \sigma)^2] \leq b_{\sigma,\size}^2 +
  \left[1 + \frac{2(-b_{\sigma,\size})_+}{-\log p} \right] \EE[Z_p^2] + \frac{2
    (b_{\sigma,\size})_+}{-\log p} \sqrt{\EE[Z_1^2] \EE[Z_p^2] }\\
  \begin{aligned}
    &+ \left[ \frac{\log \EE N_{1,\size}}{1 - (\EE N_{1,\size})^{-1}} \right]^2
    \frac{2\EE[Z_1^2]}{(-\log p)^2}%
    + \left[ \frac{\log \EE N_{1,\size}}{1 - (\EE N_{1,\size})^{-1}} \right]
    \frac{2 (-b_{\sigma,\size})_+ \EE[Z_1^2]}{-\log p}\\
    & + \left[ \frac{\log \EE N_{p,\size}}{1 - (\EE N_{p,\size})^{-1}} \right]^2
    \frac{2\EE[Z_p^2]}{(-\log p)^2}%
    + \left[ \frac{\log \EE N_{p,\size}}{1 - (\EE N_{p,\size})^{-1}} \right]
    \frac{2 (b_{\sigma,\size})_+ \EE[Z_p^2]}{-\log p}.
  \end{aligned}
\end{multline}

Then the conclusion of \cref{thm:v1:main-thm} follows from the estimates in
sections \ref{sec:v1:study-bias-term} and \ref{sec:v1:variance--1}.

\subsection{Study of the bias term $b_{\sigma,\size}$}
\label{sec:v1:study-bias-term}

The estimate for the bias term follows from the expectation of $N_{p,\size}$ for
$p \in [0,1]$ given in the next proposition. Remark that from the
definition of $N_{p,\size}$, recalling $\mathcal{V} = \Set{(t_1,x_1),\dots}$, we have
\begin{equation}
  N_{p,\size} \equaldist
  \begin{cases}
    p \sum_i \Ind_{\Ul_i\leq \size}\left(1 - (1-p)^{\sum_{j\ne
          i}B_{i,j}\Ind_{\Ul_j \leq \size}} \right) &\mathrm{if}\ 0 \leq p <
    1,\\
    \sum_i \Ind_{\Ul_i\leq \size} \Ind\Set{\sum_{j\ne i}
      B_{i,j}\Ind_{\Ul_j \leq \size} > 0} &\mathrm{if}\ p=1.
  \end{cases}
\end{equation}

\begin{proposition}
  \label{pro:v1:8}
  For any $p \in [0,1]$ and $\size > 0$ it holds,
  \begin{equation}
    \EE[N_{p,\size}] = p \size \int_{\NNReals}\left(1 - e^{-p\size \Wmarg(\uft)}
    \right)\, \intd \uft.
  \end{equation}
\end{proposition}
\begin{proof}
  We assume here that $p \ne 1$. The case $p = 1$ follows from similar steps, or by
  taking cautiously the limit $p \rightarrow 1$ in the proof. Then we have,
  \begin{align}
    \EE[N_{p,\size}]
    &= p \EE\left[%
      \sum_i \Ind_{\Ul_i \leq \size}\left(1 - (1-p)^{\sum_{j\ne
      i}B_{i,j}\Ind_{\Ul_j \leq \size}} \right) \right]\\
    &= p\EE\left[%
      \sum_i \Ind_{\Ul_i \leq \size}\left(1 - \EE\left[(1-p)^{\sum_{j\ne
      i}B_{i,j}\Ind_{\Ul_j \leq \size}} \mid \mathcal{V} \right]  \right) \right]\\
    &= p\EE\left[%
      \sum_i \Ind_{\Ul_i \leq \size}\left(1 - \textstyle\prod_{j\ne
      i}\EE\left[(1-p)^{B_{i,j}\Ind_{\Ul_j \leq \size}} \mid \mathcal{V}
      \right] \right) \right],
  \end{align}
  where the last line follows from independence and dominated convergence
  theorem to take care of the infinite product. It is easily seen from the
  definition of $\bm{B}$ that for all $i \ne j$,
  \begin{equation}
    \EE[(1-p)^{B_{i,j}\Ind_{\Ul_j\leq \size}} \mid \mathcal{V}]
    = 1 - p W(\Uft_i,\Uft_j)\Ind_{\Ul_j \leq \size}.
  \end{equation}
  It follows, invoking the Slivnyak-Mecke formula
  \citep[Chapter~13]{daley:vere-jones:2007},
  \begin{align}
    \EE[N_{p,\size}]
    &=  p\EE\left[%
      \sum_i\Ind_{\Ul_i\leq \size}\left(1 - \textstyle\prod_{j\ne
      i}\left(1 - p W(\Uft_i,\Uft_j)\Ind_{\Ul_j \leq \size} \right) \right)
      \right]\\
    &= p\size \int_{\NNReals}\left(1 - \EE\left[\textstyle\prod_j\left(1 - p
      W(\uft,\Uft_j)\Ind_{\Ul_j \leq \size} \right) \right]
      \right)\,\intd\uft.
  \end{align}
  Then the conclusion of the proposition follows by Campbell's formula
  \citep[Section~3.2]{kingman:1993} because,
  \begin{align}
    \EE\left[\textstyle\prod_j\left(1 - p
    W(\uft,\Uft_j)\Ind_{\Ul_j \leq \size} \right) \right]
    &= \EE\left[\exp\left\{ \textstyle\sum_j\Ind_{\Ul_j \leq \size}\log(1 - p
      W(\uft,\Uft_j)) \right\}\right]\\
    &= \exp\left\{-p\size
      \textstyle\int_{\NNReals}W(\uft,\uft')\,\intd \uft'
      \right\}. &\qedhere
  \end{align}
\end{proof}

It is clear from the result of the previous proposition that we have,
\begin{align}
  \log \frac{\EE[N_{1,\size}]}{\EE[N_{p,\size}]}
  &= \log\left[\frac{\int_{\NNReals} \left(1 - e^{-\size \mu(z)} \right)\,\intd
    z}{p \int_{\NNReals}\left(1 - e^{-p\size \mu(z)}\right)\, \intd z} \right]\\
  &= - (1+\sigma) \log p - \log\left[\frac{\int_{\NNReals}\left(1 - e^{-p\size
    \mu(z)}\right)\, \intd z}{p^{\sigma}\int_{\NNReals}\left(1 - e^{-\size
    \mu(z)}\right)\, \intd z} \right].
\end{align}
It then follows from \cref{ass:v1:2} that,
\begin{equation}
  |b_{\sigma,\size}|
  \leq \frac{\Gamma_{p,\size}}{-\log p}.
\end{equation}

\subsection{Variance estimates for $N_{p,\size}$}
\label{sec:v1:variance--1}

We now compute the estimate for $\EE[Z_1^2]$ and $\EE[Z_p^2]$ in the proof of
the main theorem. We assume here that $p \ne 1$. The case $p=1$ follows by
taking cautiously the limit $p \rightarrow 1$ in the proof. To shorten coming
equations, we write $q \coloneqq 1 - p$. We then have,
\begin{multline}
  N_{p,\size}^2 = p^2 \sum_i \Ind_{\Ul_i \leq \size} \left(1 - q^{\sum_{j\ne
        i}B_{i,j}\Ind_{\Ul_j \leq \size}} \right)^2\\
  + p^2\sum_i\sum_{k\ne i} \Ind_{\Ul_i \leq \size}\Ind_{\Ul_k \leq \size}\left(1
    - q^{\sum_{j\ne i}B_{i,j}\Ind_{\Ul_j \leq \size}} \right) \left(1 -
    q^{\sum_{j\ne k}B_{k,j}\Ind_{\Ul_j \leq \size}} \right).
\end{multline}
We call $p^2 S_1$ the first term of the rhs of the previous equation, and $p^2
S_2$ the second term, so that $N_{p,\size}^2 = p^2 S_1 + p^2 S_2$. That is,
\begin{gather}
  S_1 \coloneqq \sum_i \Ind_{\Ul_i \leq \size} \left(1 - q^{\sum_{j\ne
        i}B_{i,j}\Ind_{\Ul_j \leq \size}} \right)^2\\
  S_2 \coloneqq \sum_i\sum_{k\ne i} \Ind_{\Ul_i \leq \size}\Ind_{\Ul_k \leq
    \size}\left(1 - q^{\sum_{j\ne i}B_{i,j}\Ind_{\Ul_j \leq \size}} \right)
  \left(1 - q^{\sum_{j\ne k}B_{k,j}\Ind_{\Ul_j \leq \size}} \right),
\end{gather}

% \subsubsection{Computation of $\EE[S_1]$}
% \label{sec:v1:computation-ees_1}
\par We now compute $\EE[S_1]$. Expanding the square,
\begin{equation}
  S_1= \sum_i \Ind_{\Ul_i \leq \size}\left(1 - 2q^{\sum_{j\ne i}
      B_{i,j}\Ind_{\Ul_j \leq \size}} + q^{2\sum_{j\ne i} B_{i,j}\Ind_{\Ul_j
        \leq \size}} \right)
\end{equation}
By independence, and using the dominated convergence theorem to take care of the
infinite product,
\begin{equation}
  \EE[S_1 \mid \mathcal{V}] = \sum_i \Ind_{\Ul_i \leq \size} \Big(1 - 2
  \textstyle\prod_{j\ne i}\EE\left[q^{B_{i,j}\Ind_{\Ul_j \leq \size}} \mid
    \mathcal{V}\right]
  + \textstyle\prod_{j\ne i}\EE\left[q^{2B_{i,j}\Ind_{\Ul_j \leq \size}} \mid
    \mathcal{V}\right] \Big).
\end{equation}
From the definition of $\bm{B}$ we get that,
\begin{gather}
  \label{eq:v1:3}
  \EE[q^{B_{i,j}\Ind_{\Ul_j \leq \size}} \mid \mathcal{V}]
  = 1 - pW(\Uft_i,\Uft_j) \Ind_{\Ul_j \leq \size},\\
  \label{eq:v1:4}
  \EE[q^{2B_{i,j}\Ind_{\Ul_j \leq \size}} \mid \mathcal{V}]
  = 1 - p(1 +q)W(\Uft_i,\Uft_j)\Ind_{\Ul_j \leq \size}.
\end{gather}
Therefore,
\begin{multline}
  \EE[S_1] = \EE\Big[\sum_i \Ind_{\Ul_i \leq \size} \Big\{ 1 - 2
  \textstyle\prod_{j\ne
    i}(1 - p W(\Uft_i,\Uft_j) \Ind_{\Ul_j \leq \size})\\
  + \textstyle\prod_{j\ne i} (1 - p(1+q)W(\Uft_i,\Uft_j) \Ind_{\Ul_j \leq
    \size}) \Big\} \Big].
\end{multline}
Invoking the Slivnyak-Mecke theorem
\citep[Chapter~13]{daley:vere-jones:2007},
\begin{multline}
  \EE[S_1] = \size \int_{\NNReals}\Big\{ 1 - 2\EE[\textstyle\prod_j(1 - p
  W(\uft,\Uft_j) \Ind_{\Ul_j \leq \size})]\\
  + \EE[\textstyle\prod_j(1 - p(1+q)W(\uft,\Uft_j)\Ind_{\Ul_j \leq \size})]
  \Big\}\, \intd \uft.
\end{multline}
Therefore, by Campbell's formula \citep[Section~3.2]{kingman:1993} (see also the
proof of \cref{pro:v1:8}),
\begin{equation}
  \EE[S_1]
  = \size \int_{\NNReals}\left\{ 1 - 2e^{-p\size \Wmarg(\uft)} +
    e^{-p(1+q)\size \Wmarg(\uft)} \right\}\, \intd\uft.
\end{equation}
That is,
\begin{equation}
  \EE[S_1] = \size \int_{\ufSpace} \left\{1 - e^{-p \size \Wmarg(\uft)}
  \right\}^2\,
  \ufMeasure(\intd \uft)
  + \size \int_{\ufSpace}e^{-p(1+q)\size \Wmarg(\uft)}\left\{1 - e^{-p^2\size
      \Wmarg(\uft)} \right\}\, \ufMeasure(\intd\uft).
\end{equation}

% \subsubsection{Computation of $\EE[S_2]$}
% \label{sec:v1:computation-ees_2}
We now compute $\EE[S_2]$. We start with the expansion of the product in the
expression of $S_2$,
\begin{alignat}{2}
  S_2&= \sum_i\sum_{k\ne i} \Ind_{\Ul_i \leq \size}\Ind_{\Ul_k \leq \size}\Big(1
  && - q^{\sum_{j\ne i}B_{i,j}\Ind_{\Ul_j \leq \size}} -
  q^{\sum_{j\ne k}B_{k,j}\Ind_{\Ul_j \leq \size}}\\
  & && + q^{\sum_{j\ne i}B_{i,j}\Ind_{\Ul_j \leq \size} + \sum_{j\ne k}
    B_{k,j}\Ind_{\Ul_j \leq
      \size}} \Big)\\
  &= \sum_i\sum_{k\ne i} \Ind_{\Ul_i \leq \size}\Ind_{\Ul_k \leq \size}\Big(1
  && - q^{\sum_{j\ne i}B_{i,j}\Ind_{\Ul_j \leq \size}} -
  q^{\sum_{j\ne k}B_{k,j}\Ind_{\Ul_j \leq \size}}\\
  & && + q^{2 B_{i,k} +
    \sum_{j\ne i,j\ne k}B_{i,j}\Ind_{\Ul_j \leq \size} + \sum_{j\ne k,j\ne i}
    B_{k,j}\Ind_{\Ul_j \leq \size}} \Big).
\end{alignat}
The following is justified by independence and dominated convergence theorem,
\begin{multline}
  \EE[S_2 \mid \mathcal{V}] = \sum_i \sum_{k\ne i} \Ind_{\Ul_k \leq \size}
  \Ind_{\Ul_j\leq \size} \Big( 1 - \textstyle\prod_{j\ne i}
  \EE[q^{B_{i,j}\Ind_{\Ul_j \leq \size}} \mid \mathcal{V}]\\
  - \textstyle\prod_{j\ne k} \EE[q^{B_{k,j}\Ind_{\Ul_j \leq \size}} \mid \mathcal{V}] +
  \EE[q^{2B_{i,k}} \mid \mathcal{V}] \textstyle\prod_{j\ne i,k}\EE[q^{(B_{i,j} +
    B_{k,j})\Ind_{\Ul_j\leq \size}} \mid \mathcal{V}] \Big),
\end{multline}
That is, because of equations \eqref{eq:v1:3} and \eqref{eq:v1:4},
\begin{multline}
  \EE[S_2 \mid \mathcal{V}] = \sum_i \sum_{k\ne i} \Ind_{\Ul_i \leq \size} \Ind_{\Ul_k
    \leq \size} \Big( 1 - \textstyle\prod_{j\ne i}(1 - p
  W(\Uft_i,\Uft_j)\Ind_{\Ul_j \leq \size})\\
  - \textstyle\prod_{j\ne k}(1 - p W(\Uft_k,\Uft_j)\Ind_{\Ul_j \leq
    \size}) +(1 - p(1+q)W(\Uft_i,\Uft_k)) \\
  \times\textstyle\prod_{j\ne i,k}(1 - p
  W(\Uft_i,\Uft_j)\Ind_{\Ul_j \leq \size})(1 - pW(\Uft_k,\Uft_j)\Ind_{\Ul_j \leq
    \size}) \Big).
\end{multline}
We rewrite the previous in a slightly different form in order to be able to use
the Slivnyak-Mecke theorem \citep[Chapter~13]{daley:vere-jones:2007},
\begin{multline}
  \EE[S_2 \mid \mathcal{V}] = \sum_i \sum_{k\ne i} \Ind_{\Ul_i \leq \size} \Ind_{\Ul_k
    \leq \size} \Big( 1 - \textstyle\prod_{j\ne i}(1 - p
  W(\Uft_i,\Uft_j)\Ind_{\Ul_j \leq \size})\\
  - \frac{1 - pW(\Uft_k,\Uft_i)}{1 - p W(\Uft_k,\Uft_k)} \textstyle\prod_{j\ne
    i}(1 - pW(\Uft_k,\Uft_j)\Ind_{\Ul_j\leq \size})\\
  +(1 - p(1+q)W(\Uft_i,\Uft_k))
  \frac{\textstyle\prod_{j\ne i}(1 - p W(\Uft_i,\Uft_j)\Ind_{\Ul_j \leq
    \size})(1 - pW(\Uft_k,\Uft_j)\Ind_{\Ul_j \leq \size})}{(1 -
  pW(\Uft_i,\Uft_k))(1 - p W(\Uft_k,\Uft_k))} \Big).
\end{multline}
Then we can apply the Slivnyak-Mecke theorem to find that,
\begin{multline}
  \EE[S_2]
  = \size \int_{\NNReals}\EE\Big[ \sum_k \Ind_{\Ul_k \leq \size}\Big(1 -
  \textstyle\prod_j(1 - pW(\uft,\Uft_j)\Ind_{\Ul_j \leq \size})\\
  - \frac{1 - p W(\Uft_k,\uft)}{1 - p
    W(\Uft_{k},\Uft_k)}\textstyle\prod_j(1 - p W(\Uft_k,\Uft_j)\Ind_{\Ul_j \leq
    \size})\\
  +(1 - p(1+q)W(\uft,\Uft_k))
  \frac{\textstyle\prod_j(1 - p W(\uft,\Uft_j)\Ind_{\Ul_j \leq
    \size})(1 - pW(\Uft_k,\Uft_j)\Ind_{\Ul_j \leq \size})}{(1 -
  pW(\uft,\Uft_k))(1 - p W(\Uft_k,\Uft_k))}
  \Big)\Big]\,\intd\uft.
\end{multline}
Again, we rewrite the previous in a more convenient form to apply the
Slivnyak-Mecke theorem a second time,
\begin{multline}
  \EE[S_2]
  = \size \int_{\NNReals}\EE\Big[ \sum_k \Ind_{\Ul_k \leq \size}\Big(
  1 - (1 - pW(\uft,\Uft_k) \textstyle\prod_{j\ne k}(1 -
  pW(\uft,\Uft_j)\Ind_{\Ul_j
    \leq \size})\\
  - (1 - p W(\Uft_k,\uft))\textstyle\prod_{j\ne k}(1 -
  pW(\Uft_k,\Uft_j)\Ind_{\Ul_j \leq \size})\\
  +(1 - p(1+q)W(\uft,\Uft_k)) \textstyle\prod_{j\ne k}(1 - p
  W(\uft,\Uft_j)\Ind_{\Ul_j \leq \size})(1 - pW(\Uft_k,\Uft_j)\Ind_{\Ul_j \leq
    \size}) \Big)\Big]\,\intd\uft.
\end{multline}
Applying the Slivnyak-Mecke theorem to the previous,
\begin{multline}
  \EE[S_2] = \size^2 \int_{\NNReals}\int_{\NNReals} \Big\{ 1 - (1 -
  pW(\uft,\uft')) \EE[\textstyle \prod_j(1 - p
  W(\uft,\Uft_j)\Ind_{\Ul_j \leq \size})]\\
  - (1 - p W(\uft,\uft'))\EE[\textstyle \prod_j (1 - p
  W(\uft',\Uft_k)\Ind_{\Ul_j \leq \size})]\\
  +(1 - p(1+q)W(\uft,\uft')) \EE[\textstyle\prod_j (1 - p
  W(\uft,\Uft_j)\Ind_{\Ul_j \leq \size})(1 - p W(\uft',\Uft_j)\Ind_{\Ul_j \leq
    \size})] \Big\}\\
  \times \intd\uft \intd \uft'.
\end{multline}
Using Campbell's formula \citep[Section~3.2]{kingman:1993} (see also the proof
of \cref{pro:v1:8}), and recalling
$\Wcorr(\uft,\uft') = \int_{\NNReals} W(\uft,y)W(y,\uft')\,\intd y$,
\begin{multline}
  \EE[S_2] = \size^2 \int_{\NNReals}\int_{\NNReals} \Big\{ 1 - (1 -
  pW(\uft,\uft')) e^{-p \size \Wmarg(\uft)}
  - (1 - p W(\uft,\uft'))e^{-p\size \Wmarg(\uft')}\\
  +(1 - p(1+q)W(\uft,\uft'))e^{-p \size \Wmarg(\uft) - p\size \Wmarg(\uft') +
    p^2\size \Wcorr(\uft,\uft')} \Big\}\, \intd\uft \intd \uft'.
\end{multline}

% \subsubsection{Bound on the variance}
% \label{sec:v1:bound-variance}
\par We are now in position to bound $\EE[Z_p^2]$, using the expression of
$\EE[S_1]$, $\EE[S_2]$ and \cref{pro:v1:8}. \Cref{pro:v1:8} gives
\begin{align}
  \EE[N_{p,\size}]^2
  &= (p \size)^2 \int_{\NNReals}\int_{\NNReals} \left(1 - e^{-p\size
    \Wmarg(\uft)}\right) \left(1 - e^{-p\size
    \Wmarg(\uft')}\right)\, \intd \uft \intd\uft'\\
  &= (p \size)^2\int_{\NNReals}\int_{\NNReals} \Big\{ 1 - e^{-p \size
    \Wmarg(\uft)} - e^{-p\size \Wmarg(\uft')} + e^{-p \size \Wmarg(\uft) -
    p\size \Wmarg(\uft')} \Big\} \, \intd \uft \intd\uft'.
\end{align}
Combining this with the expression for $\EE[S_1]$ and $\EE[S_2]$, we get
\begin{multline}
  \EE[N_{p,\size}^2] - \EE[N_{p,\size}]^2 = \size \int_{\NNReals}e^{-p(1+q)\size
    \Wmarg(\uft)}\left\{1 - e^{-p^2\size
      \Wmarg(\uft)} \right\}\, \intd\uft\\
  + \size \int_{\NNReals} \left\{1 - e^{-p \size \Wmarg(\uft)} \right\}^2\, \intd \uft%
  +2 p^3\size^2 \int_{\NNReals} \Wmarg(\uft) e^{-p\size
    \Wmarg(\uft)}\,\intd \uft\\
  + p^2\size^2 \int_{\NNReals}\int_{\NNReals} e^{-p \size \Wmarg(\uft) -p\size
    \Wmarg(\uft') + p^2\size \Wcorr(\uft,\uft')} \left\{1 - e^{-p^2\size
      \Wcorr(\uft,\uft')} \right\}\,\intd\uft \intd \uft'\\
  - p^3(1+q) \size^2 \int_{\NNReals}\int_{\NNReals}W(\uft,\uft') e^{-p \size
    \Wmarg(\uft) - p\size \Wmarg(\uft') + p^2\size
    \Wcorr(\uft,\uft')}\,\intd\uft\intd\uft'.
\end{multline}
By Hölder's and Young's inequality, we have
$2\Wcorr(\uft,\uft') \leq 2\sqrt{\Wmarg(\uft)\Wmarg(\uft')} \leq \Wmarg(\uft) +
\Wmarg(\uft')$. Moreover, $p^2 \leq p$, and $1- e^{-x} \leq x$ imply,
\begin{multline}
  \int_{\NNReals}\int_{\NNReals} e^{-p \size \Wmarg(\uft) -p\size \Wmarg(\uft')
    + p^2\size \Wcorr(\uft,\uft')} \left\{1 - e^{-p^2\size \Wcorr(\uft,\uft')}
  \right\}\,\intd\uft \intd \uft'\\
  \leq p^2 \size \int_{\NNReals}\int_{\NNReals} e^{-\frac{p\size}{2}
    \Wmarg(\uft)}e^{-\frac{p\size}{2}\Wmarg(\uft')}
  \Wcorr(\uft,\uft')\, \intd \uft \intd \uft'.
\end{multline}
If $\sigma = 1$ it follows from \cref{ass:v1:5} that,
\begin{multline}
  \int_{\NNReals}\int_{\NNReals} e^{-p \size \Wmarg(\uft) -p\size \Wmarg(\uft')
    + p^2\size \Wcorr(\uft,\uft')} \left\{1 - e^{-p^2\size \Wcorr(\uft,\uft')}
  \right\}\, \intd\uft \intd \uft'\\
  \leq Cp^2\size \left\{\int_{\NNReals} \Wmarg(\uft) \, \intd \uft
  \right\}^2,
\end{multline}
which is finite because of \cref{ass:finiteness}. Now if $0 \leq \sigma < 1$, it follows
from \cref{ass:v1:5} and \cite[Lemma~B.4]{caron:panero:rousseau:2017} that,
\begin{multline}
  \int_{\NNReals}\int_{\NNReals} e^{-p \size \Wmarg(\uft) -p\size \Wmarg(\uft')
    + p^2\size \Wcorr(\uft,\uft')} \left\{1 - e^{-p^2\size \Wcorr(\uft,\uft')}
  \right\}\,\intd\uft \intd \uft'\\ \leq C
  p^2\size^{1-2\eta} \left\{\int_{\NNReals} \alpha^{\eta} \Wmarg(\uft)^{\eta}
    e^{-\frac{p\size}{2}\Wmarg(\uft)}\, \intd \uft \right\}^2
  \lesssim \size^{1 - 2\eta + 2\sigma}\ell(\size)^2.
\end{multline}
Moreover, using \cite[Lemma~B.4]{caron:panero:rousseau:2017}, it is easily seen that
all other terms involved in $\EE[N_{p,\size}^2] - \EE[N_{p,\size}]^2$ are
bounded by a multiple constant of $\size^{1+\sigma}\ell(\size)$ for any
$\sigma \in [0,1]$. Since we set $\eta = 1$ when $\sigma = 1$, it follows the
estimate,
\begin{equation}
  \EE[N_{p,\size}^2] - \EE[N_{p,\size}]^2
  \lesssim \size^{1+\sigma}\ell(\size) + \size^{3 - 2\eta +
    2\sigma}\ell(\size)^2.
\end{equation}
The conclusion follows because
$\EE[N_{p,\size}]^2 \gtrsim \size^{2 + 2\sigma}\ell(\size)^2$ by \cref{pro:v1:8}.

\section{Proofs related to count statistics}
\label{sec:proofs-relat-bootstr}

\subsection{Proof of \cref{pro:2}}
\label{sec:proof-}

Let $(Z_v)_{v \in \ver(G)}$ be a collection of independent $\bern(p)$ random
variables. Let $G' \subseteq G$ be the subgraph of $G$ consisting on vertices
such $v \in \ver(G)$ such that $Z_v = 1$. Then $\psamp(G,p)$ is equal in distribution
to $G'$ without its isolated vertices. Since $R$ is connected $\nsub(R,\psamp(G,p))
\equaldist \nsub(R,G')$, and hence by \cref{eq:54}
\begin{align}
  \nsub(R,\psamp(G,p))
  &\equaldist \sum_{\phi \in \mathcal{A}(R,G')}\1\Set{R \subseteq G'[\phi]}.
\end{align}
Clearly $V(G') \subseteq V(G)$, and can assume wlog that
$\mathcal{A}(R,G') \subseteq \mathcal{A}(R,G)$ using canonical embedding. For
simplicity let also $Z_{\phi} \coloneqq \prod_{r \in \ver(R)}Z_{\phi(r)}$ for
any $\phi \in \mathcal{A}(R,G)$.  Then, we have the equivalence
$\phi \in \mathcal{A}(R,G) \wedge Z_{\phi} = 1 \Leftrightarrow \phi \in
\mathcal{A}(R,G')$. Therefore,
\begin{align}
  \label{eq:24}
  \nsub(R,\psamp(G,p))%
  &\equaldist \sum_{\phi \in \mathcal{A}(R,G)} Z_{\phi}
    \1\Set{R \subseteq G'[\phi]}\\
  \label{eq:31}
  &= \sum_{\phi \in \mathcal{A}(R,G)} Z_{\phi}%
    \1\Set{R \subseteq G[\phi]},
\end{align}
where the second line follows because we also have the equivalence $\phi \in
\mathcal{A}(R,G) \wedge Z_{\phi}=1 \Leftrightarrow G’[\phi] = G[\phi]$. Since
$\phi$ is injective $\EE[Z_{\phi}] = p^k$, whence the result.

\subsection{Proof of \cref{pro:3}}
\label{sec:proof-pro:3}

The first equality in \cref{pro:3} is obvious. Here we focus on computing and
bounding $\var(\nsub(R,\psamp(G,p)) \mid G)$. But this is almost immediate from
\cref{eq:31}. Indeed, $\var( \nsub(R,\psamp(G,p)) \mid G)$ equals
\begin{equation}
  \label{eq:45}
  \sum_{\phi_1 \in \mathcal{A}(R,G)}\sum_{\phi_2 \in \mathcal{A}(R,G)}
  \cov(Z_{\phi_1},Z_{\phi_2})\1\Set{R \subseteq G[\phi_1],\, R \subseteq
    G[\phi_2] }.
\end{equation}
Now, if $\phi_1(\ver(R)) \cap \phi_2(\ver(R)) = \varnothing$, then
$\cov(Z_{\phi_1},Z_{\phi_2}) = 0$. In general if
$|\phi_1(\ver(R)) \cap \phi_2(\ver(R))| = m$, then
$\cov(Z_{\phi_1},Z_{\phi_2}) \leq p^{2k-m}$. The result follows.

% \subsection{Proof of \cref{pro:4}}
% \label{sec:proof-pro:4}

% Since $\ver(R) = \Set{1,\dots,k}$, we can rewrite,
% \begin{equation}
%   \label{eq:1}
%   \indsub(R,G_t)%
%   =%
%   \frac{1}{k!}\sum_{\substack{u_1,\dots,u_k \in \ver(G_t)\\u_1\ne\dots\ne u_k}}%
%   \prod_{\Set{\ell,m}\in \edg(R)}\1_{\Set{u_{\ell},u_m}\in \edg(G_t)}\prod_{\Set{\ell,m}\in
%     \edg(\bar{R})}\1_{\Set{u_{\ell},u_{m}}\notin \edg(G_t)}.
% \end{equation}
% Then, introducing the variables $B_{i,j}$ such that $B_{i,j} = 1$ if
% $v_i = (t_i,x_i)$ and $v_j(t_j,x_j) = $ are connected, $B_{i,j} = 0$ otherwise,
% and
% $A_k \coloneqq \Set{(i_1,\dots,i_k) \in \Nats^k \given i_1\ne \dots\ne i_k}$,
% we have almost-surely,
% \begin{equation}
%   \label{eq:20}
%   \indsub(R,G_t)%
%   = \frac{1}{k!}%
%   \sum_{\bm{i} \in A_k}\Big(\prod_{j=1}^k\1_{t_{i_j} \leq t} \Big)%
%   \prod_{\Set{\ell,m}\in \edg(R)}B_{i_{\ell},i_m}\prod_{\Set{\ell,m}\in \edg(\bar{R})}(1
%   - B_{i_{\ell},i_m}).
% \end{equation}
% Then, the formula in \cref{pro:4} follows from Slivnyak-Mecke's formula, as it
% has been done numerous time, for instance in \citet{caron:panero:rousseau:2017}.
% % \begin{equation}
% %   \label{eq:2}
% %   \EE[\indsub(R,G_t)]%
% %   = \frac{t^k}{k!} \int_{S^k}\prod_{\Set{\ell,m}\in
% %     \edg(R)}W(x_{\ell},x_m)\prod_{\Set{\ell,m}\in \edg(\bar{R})}(1 -
% %   W(x_{\ell},x_m))\,\intd x_1,\dots\intd x_k
% % \end{equation}

% % --------------------------------------------------------------------------------

% % We now consider the variance of $\indsub(R,G_t)$, under the convenient
% % assumption that $W(x,x) = 0$.

\subsection{Proof of \cref{pro:1}}
\label{sec:proof-pro:1}

Let us assume for now that there is a universal $c > 0$ (\textit{i.e.} depending
only on $W$) such that as $t \to \infty$
\begin{equation}
  \label{eq:27}
  s_t^2(R) = \var_W(\nsub(R,G_t)) = c t^{2k-1} + O(t^{2k-2}).
\end{equation}
We delay the proof of this claim to later. Since
$G_{pt} \equaldist \psamp(G_t,p)$ by
\cite{borgs:chayes:cohn:veitch:2017,veitch:roy:2019}, we have
$\var_W(\nsub(R_1,G_{pt})) = \var_W(\nsub(R,
\psamp(G_t,p)))$, and hence as $t\to \infty$
\begin{equation}
  \label{eq:28}
  p^{2k-1}\var_W(\nsub(R,G_t))%
  = \var_W(\nsub(R,\psamp(G_t,p))) + O(t^{2k-2}).
\end{equation}
Then, by the law of total variance,
\begin{multline}
  \label{eq:10}
  p^{2k - 1}\var_W(\nsub(R,G_t))\\
  \begin{aligned}
    &=%
    \EE_W[\var(\nsub(R,\psamp(G_t,p)) \mid G_t)] + \var_W\big(
    \EE[\nsub(R,\psamp(G_t,p))\mid G_t] \big)%
    + O(t^{2k-2})\\
    &= \EE_W[\var(\nsub(R,\psamp(G_t,p))\mid G_t)]%
    + p^{2k}\var_W(\nsub(R,G_t))%
    + O(t^{2k-2}),
  \end{aligned}
\end{multline}
where the last line follows from \cref{pro:2}. Rearranging the last display
gives the proposition.

Hence, it remains to prove the \cref{eq:27}. Let $\tilde{G} \in \mathcal{L}$
be the graph constructed in \cref{sec:model}, and let $\tilde{G}_t \subseteq
\tilde{G}$ be the subgraph consisting only on vertices $v_i = (t_i,x_i)$ with
$t_i \leq t$. Likewise, $G_t$ is (in distribution) the subgraph of $\tilde{G}_t$
consisting only on non-isolated vertices. But, since $R$ is connected, we have
$\nsub(R,G_t) \equaldist \nsub(R, \tilde{G}_t)$. Write $\mathcal{V} =
\Set{(t_1,x_1),\dots}$ so that we can assume wlog that $V(\tilde{G}_t) = \Nats$
and $V(R) = \Nats_k$, and we can rewrite,
\begin{equation}
  \label{eq:46}
  \nsub(R,G_t)%
  \equaldist \sum_{\substack{i_1,\dots,i_k \in \Nats\\i_1\ne \dots\ne i_k}}\Big(
  \prod_{j=1}^k\1_{t_{i_j} \leq t}\Big) \prod_{\Set{\ell,m}\in E(R)}B_{i_{\ell},i_m}.
\end{equation}
where $B_{i,j} =1$ iff there is an edge between $(t_i,x_i)$ and $(t_j,x_j)$ in
$\tilde{G}$. For simplicity, we define
$\tilde{\nsub}_t(R,\mathcal{V}) \coloneqq \EE_W[\nsub(R,G_t) \mid
\mathcal{V}]$. Then, by the law of total variance, we can decompose
\begin{align}
  \label{eq:50}
  \var_W(\nsub(R,G_t))%
  &= \EE[\var_W(\nsub(R,G_t) \mid \mathcal{V})]%
    + \var(\tilde{\nsub}_t(R,\mathcal{V})),
\end{align}
We now estimate the terms involved in the previous rhs. We write
$A_k \coloneqq \Set{(i_1,\dots,i_k)\in \Nats^k \given i_1\ne \dots \ne i_k}$,
\textit{i.e.}  the set of $k$-arrangements of $\Nats$. We also let
$B_{\bm{i}} \coloneqq \prod_{\Set{\ell,m}\in
  E(R)}B_{i_{\ell},i_m}\1_{t_{i_{\ell}}\leq t}\1_{t_{i_m}\leq t}$, for any
$\bm{i} \in A_k$. Then, by \cref{eq:46},
\begin{align}
  \label{eq:55}
  \var_W(\nsub(R,G_t) \mid \mathcal{V})
  &\equaldist%
    \sum_{\bm{i}\in A_k}\sum_{\bm{i}' \in A_k}\cov_W(B_{\bm{i}},B_{\bm{i}'}\mid \mathcal{V}).
\end{align}
We remark that $\cov_W(B_{\bm{i}},B_{\bm{i}'} \mid \mathcal{V}) = 0$ whenever
$\Set{i_{\ell},i_m}$ and $\Set{i'_{\ell}, i'_m}$ are all distinct. In other
words, $\cov_W(B_{\bm{i}},B_{\bm{i}'} \mid \mathcal{V}) \ne 0$ implies that there
exists a pair $(i_{\ell},i_m)$ that is equal to $(i'_{\ell'},i'_{m'})$. So under
\cref{ass:expecmotifs}, we deduce from \cref{eq:55} and the Slivnyak-Mecke,
formula that
\begin{align}
  \label{eq:56}
  \EE[\var_W(\nsub(R,G_t) \mid \mathcal{V})]%
  = O(t^{2k-2}).
\end{align}
Then, it is enough to establish that $\var(\tilde{\nsub}_t(R,\mathcal{V})) \sim c
t^{2k-1}$ to finish the proof. We note that by \cref{eq:46}
\begin{align}
  \label{eq:57}
  \tilde{\nsub}_t(R,\mathcal{V})%
  \equaldist \sum_{\bm{i}\in A_k}\Big(
  \prod_{j=1}^k\1_{t_{i_j} \leq t}\Big) \prod_{\Set{\ell,m}\in E(R)}W(x_{i_{\ell}},x_{i_m}).
\end{align}
We let $N_k : A_k^2 \to \NNInts$ such that
$N_k(\bm{i},\bm{i}')$ counts the number of distinct elements in
$\bm{i}\cup \bm{i}'$. Then,
\begin{align}
  \label{eq:21}
  \tilde{\nsub}_t(R,\mathcal{V})^2%
  &\equaldist \sum_{m=k}^{2k}\sum_{\bm{i}\in A_k}\sum_{\bm{i}'\in
    A_k}\1\Set{N_k(\bm{i},\bm{i}') = m}\\%
  &\quad%
    \times \Big(\prod_{j=1}^k\1_{t_{i_j}\leq t}\1_{t_{i'_{j}}\leq t}\Big)%
    \prod_{\Set{\ell,m}\in \edg(R)} W(x_{i_{\ell}},x_{i_m})%\\
%  &\quad%
%    \times \Big(\prod_{j=1}^k\Big)%
%    \prod_{\Set{\ell,m}\in \edg(R)}
     W(x_{i'_{\ell}}x_{i'_m}).
\end{align}
All the terms in the previous summation have finite expectation under
\cref{ass:expecmotifs} if $N \geq 2k$. Also, using Slivnyak-Mecke's
formula, we see that the expectation of the term corresponding to $m=2k$ is
exactly $\EE[\tilde{\nsub}_t(R,\mathcal{V})]^2$ (so indeed it is enough to
have $N = 2k - 1)$. Similarly, the expectation
of the other terms is seen to be $O(t^m)$. It follows that,
\begin{multline}
  \label{eq:22}
  \var(\tilde{\nsub}_t(R,\mathcal{V}))%
  =%
  O(t^{2k-2})\\
    + \EE\Bigg[\sum_{\substack{\bm{i},\bm{i}' \in A_k\\N_k(\bm{i},\bm{i}') =
 2k-1} }
\Big(\prod_{j=1}^k\1_{t_{i_j}\leq t}\1_{t_{i'_{j}}\leq t}\Big)%
    \prod_{\Set{\ell,m}\in \edg(R)} W(x_{i_{\ell}},x_{i_m})W(x_{i'_{\ell}},x_{i'_m})\Bigg],
\end{multline}
which behave as $ct^{2k-1}$ when $t \to \infty$, for a constant $c> 0$
depending only on $W$ and $R$.

\subsection{Proof of \cref{pro:5}}
\label{sec:proof-pro:5}

We define for simplicity
$c(G_t) \coloneqq \var(\nsub(R,\psamp(G_t,p)) \mid G_t)$. So, in view of
\cref{pro:1} and its proof, we then have,
\begin{equation}
  \label{eq:49}
  \frac{\hat{s}^2(R,G_t)}{s_t^2(R)}%
  = \frac{c(G_t)}{\EE_W[c(G_t)]} \Big(1 + O(t^{-1})\Big).
\end{equation}
Thus it is enough to show that $\frac{\var(c(G_t))}{\EE_W[c(G_t)]^2} = O(t^{-1})$.
We have obtained in \cref{sec:proof-pro:1} that $\EE_W[c(G_t)] \asymp t^{2k-1}$,
thus the proposition is proved if we show that $\var_W(c(G_t)) = O(t^{4k-3})$. We
use the same notations and definitions as in \cref{sec:proof-pro:1}. Then, with
the same arguments as in \cref{sec:proof-pro:3,sec:proof-pro:1} we obtain that
\begin{align}
  \label{eq:29}
  c(G_t)%
  &\equaldist%
    \sum_{m=k}^{2k-1}p^m(1 - p^{2k-m})C_m,
\end{align}
where,
\begin{align}
  \label{eq:37}
  C_m%
  &\coloneqq%
    \sum_{\substack{\bm{i},\bm{i}' \in A_k\\N_k(\bm{i},\bm{i}') = m}}
    \Big(\prod_{j=1}^k \1_{t_{i_j}\leq t}\1_{t_{i'_j}\leq t}\Big)
    \prod_{\Set{\ell,m}\in \edg(R)}B_{i_{\ell},i_{m}}B_{i'_{\ell},i'_m}.
\end{align}
Therefore,
\begin{equation}
  \label{eq:30}
  c(G_t)^2%
  \equaldist%
  \sum_{m_1=k}^{2k-1}\sum_{m_2=k}^{2k-1}p^{m_1}(1
  - p^{2k-m_1})%
  p^{m_2}(1 - p^{2k-m_2})C_{m_1}C_{m_2}.%
\end{equation}
It is easily seen that under \cref{ass:expecmotifs} with
$N \geq 2k_1 + 2k_2 - 2$, all the terms involved in \cref{eq:30} have finite
expectations. Also, the terms corresponding to $m_1 + m_2 \leq 4k - 3$
are seen to have expectation $O(t^{4k-3})$. But the only way to have
$m_1 + m_2 > 4k - 3$ is to have $m_1 = 2k - 1$ and
$m_2 = 2k -1$. Hence, we deduce that
\begin{align}
  \label{eq:30b}
  \EE_W[c(G_t)^2]%
  &= O(t^{4k-3}) + %
    p^{4k-2}(1-p)^2\EE_W[C_{2k-1}^2].%  \\
\end{align}
We also note that by \cref{eq:29} and the same arguments that led to the
previous formula,
\begin{align}
  \label{eq:39}
  \EE_W[c(G_t)]^2%
  = p^{4k-2}(1-p)^2\EE_W[C_{2k-1}]^2%
  +O(t^{4k-3}).
\end{align}
Hence, by combining \cref{eq:30,eq:39},
\begin{align}
  \label{eq:40}
  \var_W(c(G_t))%
  &\leq \var_W(C_{2k-1}) + O(t^{4k - 3}).
\end{align}
So in fact it is  enough to bound $\var_W(C_{2k-1})$. We note that
$N_{k}(\bm{i},\bm{i}') = 2k-1$ means that there is exactly one of
the $i_{1},\dots,i_{k}$ which is equal to one of the
$i_{1}',\dots,i_{k}'$. Hence, we can rewrite,
\begin{align}
  \label{eq:41}
  C_{2k-1}%
  &= \sum_{a_1=1}^{k}\sum_{a_2=1}^{k}%
    \tilde{C}_{2k-1}^{a_1,a_2},
\end{align}
where,
\begin{align}
  \label{eq:43}
  \tilde{C}_{2k-1}^{a_1,a_2}%
  &\coloneqq%
    \sum_{\substack{\bm{i},\bm{i}'\in A_k,i_{a_1}=i'_{a_2}\\N_k(\bm{i},\bm{i}')=2k-1}}%
  \Big(\prod_{j=1}^k \1_{t_{i_j}\leq t}\1_{t_{i'_j}\leq t}\Big)
  \prod_{\Set{\ell,m}\in \edg(R)}B_{i_{\ell},i_{m}}B_{i'_{\ell},i'_m}.
\end{align}
And then,
\begin{align}
  \label{eq:42}
  \var_W(C_{2k-1})%
  &\lesssim%
    \sum_{a_1=1}^{k}\sum_{a_2=1}^{k}\var_W(\tilde{C}_{2k-1}^{a_1,a_2}).
\end{align}
To finish the proof, we note that
$\tilde{C}_{2k-1}^{a_1,a_2} \equaldist \nsub(Q,G_t)$, where $Q$ is the connected
graph with $2k - 1$ vertices obtained by joining two copies of $R$, $R_1$ and
$R_2$, in such a way that vertex $a_1$ of $R_1$ and $a_2$ of $R_2$ becomes the
same vertex in $Q$. Hence $\var_W(C_{2k-1}) = O(t^{4k-3})$ by \cref{eq:27}, and
$\var_W(c(G_t)) = O(t^{4k-3})$ too.

\subsection{Proof of \cref{pro:6}}
\label{sec:proof-pro:6}

We let $W \equiv W_0$ and $t \equiv t_0$ for simplicity. We define
$\tilde{\nsub}_t(R,\mathcal{V}) \coloneqq \EE_W[\nsub(R,G_t) \mid \mathcal{V}]$
as in \cref{sec:proof-pro:1}. Then, we rewrite the test statistic
$T \equiv T(R,G_t)$ as
\begin{equation}
  \label{eq:12}
  T = %
  \underbrace{\frac{\tilde{\nsub}_t(R,\mathcal{V}) 
      - \bar{\nsub}(R,W,t)}{\sqrt{\var(\tilde{\nsub}_t(R,\mathcal{V}))}}
  }_{\eqqcolon T_1}%
  \cdot
  \frac{\sqrt{\var(\tilde{\nsub}_t(R,\mathcal{V}))}}{\sqrt{\hat{s}^2(R,G_t)}}%
  + \underbrace{\frac{\nsub(R,G_t) -
      \tilde{\nsub}_t(R,\mathcal{V})}{\sqrt{\hat{s}^2(R,G_t)}}}_{\eqqcolon T_2}.
\end{equation}
We note that
$\frac{\hat{s}^2(R,G_t)}{s_t^2(R)} = 1 + O_p(t^{-1/2})$ by
\cref{pro:5} and
$\frac{\var(\tilde{\nsub}_t(R,\mathcal{V}))}{s_t^2(R)} = 1 + O_p(t^{-1})$ by
combining \cref{eq:50,eq:56}. Thus
\begin{align}
  \label{eq:16}
  \frac{\var(\tilde{\nsub}_t(R,\mathcal{V}))}{\hat{s}^2(R,G_t)} = 1 + O_p(t^{-1/2}).
\end{align}
Also, $\tilde{\nsub}_t(R,\mathcal{V}) = \EE_W[\nsub(R,G_t\mid \mathcal{V}]$, and
by the same argument as before, it is clear that $T_2 = O_p(T_2')$, where $T_2'
\coloneqq \frac{\nsub(R,G_t) - \EE_W[\nsub(R,G_t)]}{\sqrt{s_t^2(R)}}$. Clearly $\EE_W[T_2'] = 0$, and
\begin{align}
  \label{eq:14}
  \var_W(T_2')%
  &= \frac{\EE[\var_W(\nsub(R,G_t) \mid
    \mathcal{V})]}{s_t^2(R)}%
    = O_p(t^{-1}),
\end{align}
by \cref{eq:56} again, and consequently $T_2 = O_p(t^{-1/2})$. That entails
$T = T_1 + O_p(t^{-1/2})$, so the conclusion follows from Slutsky's lemma if
we show that $T_1 \to \mathcal{N}(0,1)$. But, we remark that
$\tilde{\nsub}_t(R,\mathcal{V})$ is a $U$-statistic of a Poisson process (see
the \cref{eq:57}). We note that it was not the case for $\nsub(R,G_t)$ because
of the randomness in the edge connections. Under the \cref{ass:expecmotifs} with
$N \geq 4k$, we have $T_1 \xrightarrow{d} \mathcal{N}(0,1)$ by
\cite[Theorem~5.2]{reitzner:schulte:2013}.

\section{Proof of examples rates bounds}
\label{sec:v1:proof-ggp-rates}

\subsection{Dense graphs}
\label{sec:v1:dense-graph-example}

For any $\size > 0$,
\begin{multline}
  \int_{\NNReals}\left(1 - e^{-\size \mu(x)}\right)\,\intd x
  - \int_{\NNReals}\left(1 - e^{-p\size \mu(x)}\right)\,\intd x\\
  = \int_0^1 e^{-p\size u /2}\left\{1 - e^{-(1-p)\size u /2}\right\}\, \intd u
  \lesssim \size^{-1}.
\end{multline}
On the other hand, we have for $\size$ large enough,
\begin{equation}
  \int_0^1 \left(1 - e^{-\size u/2}\right)\,\intd u \gtrsim 1.
\end{equation}
Hence, $\Gamma_{p,\size} \propto \size^{-1}$.

\subsection{Sparse, almost dense graphs without power law}
\label{sec:v1:sparse-almost-dense-1}

For any $\size > 0$,
\begin{multline}
  \int_{\NNReals}\left(1 - e^{-\size \mu(x)}\right)\,\intd x
  - \int_{\NNReals}\left(1 - e^{-p\size \mu(x)}\right)\,\intd x\\
  \begin{aligned}
    &= \int_0^1 e^{-p\size u}\left\{1 - e^{-(1-p)\size u}\right\}\,
    \frac{\intd u}{u}\\
    &\leq \size(1-p)\int_0^{1/\size} e^{-p\size u}\,\intd u +
    \int_{1/\size}^{1}e^{-p \size u}\frac{\intd
      u}{u}\\
    &\leq (1-p) + \int_1^{\size} \frac{e^{-p v}}{v}\,\intd v \lesssim 1.
  \end{aligned}
\end{multline}
On the other hand, when $\size > 0$ is large enough,
\begin{align}
  \int_0^1 \left(1 - e^{-\size u}\right)\,\frac{\intd u}{u}
  &= (1 - e^{-\size})\log(\size) - \int_0^{\size}\log(v)e^{-v}\,\intd v\\
  &\geq (1 - e^{-\size})\log(p\size) - \int_1^{\size} \log(v)e^{-v}\,\intd v\\
  &\geq (1 - e^{-\size})\log(\size) - \log(\size)
    \int_1^{\size}e^{-v}\,\intd v \gtrsim \log(\size).
\end{align}
Hence $\Gamma_{p,\size} \propto 1/\log\size$.

\subsection{Sparse graphs with power law}
\label{sec:v1:sparse-graphs-with}

Let $b \coloneqq 1/(1-\sigma)$, then for any $\size > 1$,

\begin{multline}
  p^{\sigma}\int_{\NNReals}\left(1 - e^{-\size \mu(x)}\right)\,\intd x -
  \int_{\NNReals}\left(1 - e^{-p\size \mu(x)}\right)\,\intd x\\
  \begin{aligned}
    &= \sigma\left(\frac{\sigma}{1-\sigma}\right)^{\sigma} \size^{\sigma}
    \left\{ p^{\sigma}\int_0^{\frac{\sigma\size}{1-\sigma}}(1 -
      e^{-u})\,u^{-1-\sigma}\intd u - \int_0^{\frac{\sigma\size}{1-\sigma}}(1 -
      e^{-p u})\, u^{-1-\sigma}\,\intd
      u \right\}\\
    &= \sigma \left(\frac{\sigma}{1-\sigma}\right)^{\sigma} p^{\sigma}
    \size^{\sigma}
    \int_{\frac{p\sigma\size}{1-\sigma}}^{\frac{\sigma\size}{1-\sigma}}(1 -
    e^{-u})\, u^{-1-\sigma}\intd u \leq \frac{\sigma(1-p)}{p}.
  \end{aligned}
\end{multline}
On the other hand, when $\size$ is large enough,
\begin{align}
  \int_{\NNReals}\left(1 - e^{-\size \mu(x)}\right)\,\intd x
  &= \size^{\sigma} \sigma\left(\frac{\sigma}{1-\sigma}\right)^{\sigma}
    \int_0^{\frac{\sigma\size}{1-\sigma}}(1 - e^{-v}) v^{-1-\sigma}\,\intd v\\
  &\geq \size^{\sigma} \sigma\left(\frac{\sigma}{1-\sigma}\right)^{\sigma}
    \int_0^{\frac{\sigma\size}{1-\sigma}} e^{-v} v^{-\sigma}\,\intd v \gtrsim
    \size^{\sigma}.
\end{align}
Hence $\Gamma_{p,\size} \propto \size^{-\sigma}$.

\subsection{Generalized Gamma Process}
\label{sec:v1:gener-gamma-proc-1}

In the proof we assume without loss of generality that $\tau = 1$. By
definition, we can see that the graphex marginal is given by
\begin{align}
  \Wmarg(x)
  &= \int_{\NNReals}\Big(1 - e^{-2
    \bar{\rho}_{\sigma,\tau}^{-1}(x')\bar{\rho}_{\sigma,\tau}^{-1}(x)}\Big) \intd
x'\\
  &=\int_{\NNReals}\Big(1 - e^{-2 y\bar{\rho}_{\sigma,\tau}^{-1}(x)}\Big)
    \rho_{\sigma,\tau}(\intd y)\\
  &= \frac{(1 + 2 \bar{\rho}_{\sigma,\tau}^{-1}(x) )^{\sigma} - 1}{\sigma}.
\end{align}
Then, for any $p \in [0,1]$,
\begin{align}
  \label{eq:13}
  \int_{\NNReals}\Big(1 - e^{-pt\mu(z)}\Big)\intd z%
  &= \int_{\NNReals}\Big(1 - \exp\Big(-pt \frac{(1 + 2
    \bar{\rho}_{\sigma,\tau}^{-1}(x) )^{\sigma} - 1}{\sigma}  \Big) \Big)\intd
 x\\
  &= \int_{\NNReals}\Big(1 - \exp\Big(-pt \frac{(1 + 2y )^{\sigma} - 1}{\sigma}
    \Big) \Big)\rho_{\sigma,\tau}(\intd y).
\end{align}
Since $\tau = 1$ by assumption, the previous can be rewritten as,
\begin{multline}
  \label{eq:v1:2}
  \int_{\NNReals}\left(1 - e^{-p\size \mu(z)}\right)\,\intd z =
  \frac{2^{\sigma}}{pt }\int_{\NNReals}(1 - e^{-u})\left(1 + \frac{\sigma u}{p
      \size}\right)^{-1+1/\sigma}\\
  \times \left(\left(1 + \frac{\sigma u}{p \size}\right)^{1/\sigma} -
    1\right)^{-1-\sigma} \exp\left\{- \frac{1}{2}\left(\left(1 + \frac{\sigma
          u}{p\size}\right)^{1/\sigma} -1\right) \right\}\, \intd u
\end{multline}

We will find a bound on $\Gamma_{p,\size}$ by first lower bounding,
\begin{equation}
  A \coloneqq \int_{\NNReals}\left(1 - e^{-\size
      \mu(z)}\right)\,\intd z,
\end{equation}
and then upper bounding,
\begin{equation}
  B = \left|p^{\sigma}\int_{\NNReals}\left(1 - e^{-\size
        \mu(z)}\right)\,\intd z - \int_{\NNReals}\left(1 -
      e^{-p\size \mu(z)}\right)\, \intd z \right|.
\end{equation}
From the two bounds computed below, we will conclude that
$\Gamma_{p,\size} \propto \size^{-\sigma}$.

% \subsubsection{Lower bound on $A$}
% \label{sec:v1:lower-bound-when}
\par \textit{Lower bound on $A$}. Since the integrand is a positive function, we
lower bound $A$ by integrating on a smaller set. it is clear that when
$u \in [0,\size]$,
\begin{equation}
  \left(1 + \frac{\sigma u}{\size}\right)^{-1+1/\sigma} \geq 1,\
  \mathrm{and}\quad%
  \left(1 + \frac{\sigma u}{\size}\right)^{1/\sigma} - 1 \lesssim
  \frac{u}{\size}.
\end{equation}
Therefore, from \cref{eq:v1:2} we deduce that,
\begin{equation}
  A \gtrsim \size^{\sigma}\int_{[0,\size]} (1 - e^{-u}) u^{-1-\sigma}\,\intd u
  \gtrsim \size^{\sigma}.
\end{equation}

% \subsubsection{Upper bound on $B$}
% \label{sec:v1:upper-bound-b}
\par \textit{Upper bound on $B$}. We first write,
\begin{multline}
  \int_{\NNReals}\left(1 - e^{-p\size \mu(z)}\right)\,\intd z\\
  = \int_0^{\size} \left(1 - e^{-p\size \mu(z)}\right)\,\intd z +
  \int_{\size}^{\infty}\left(1 - e^{-p\size \mu(z)}\right)\,\intd z%
  \eqqcolon C_{p,\size} + D_{p,\size}.
\end{multline}
It is easily seen from \cref{eq:v1:2} that for some constant $c > 0$ (eventually
depending on $\sigma$ and $p$, but not $\size$),
\begin{align}
  D_{p,\size}
  &\lesssim \int_{\size}^{\infty}(1 -
    e^{-u})\left(\frac{u}{\size}\right)^{-1+1/\sigma}
    \left(\frac{u}{\size}\right)^{(-1 - \sigma)/\sigma} e^{-c
    (u/\size)^{1/\sigma}}\, \frac{\intd u}{\size}\\
  &= \size \int_{\size}^{\infty}(1 - e^{-u})
    u^{-2}e^{-c(u/\size)^{1/\sigma}}\,\intd u\\
  &\leq \int_1^{\infty} u^{-2} e^{-c u^{1/\sigma}}\,\intd u \lesssim 1.
\end{align}
It turns out that $B \lesssim |p^{\sigma}C_{1,\size} - C_{p,\size}| + 1$. We now
consider the function,
\begin{multline}
  F(p,\sigma,\size,u)
  \coloneqq \frac{1}{p\size}\left(1 + \frac{\sigma u}{p
      \size}\right)^{-1+1/\sigma}\\
  \times\left(\left(1 + \frac{\sigma u}{p
        \size}\right)^{1/\sigma} - 1\right)^{-1-\sigma} \exp\left\{- \frac{1}{2}
    \left(\left(1 + \frac{\sigma u}{p\size}\right)^{1/\sigma} -1\right)
  \right\}.
\end{multline}
With a little bit of effort, it is seen that for any $u \in [0,\size]$,
\begin{equation}
  F(p,\sigma,\size,u) = u^{-1-\sigma}(\size p)^{\sigma} +
  u^{-\sigma}O(\size^{-1+\sigma}).
\end{equation}
Hence,
$p^{\sigma}F(1,\sigma,\size,u) - F(p,\sigma,\size u) =
u^{-\sigma}O(\size^{-1+\sigma})$ whenever $u \in [0,\size]$. It follows that,
\begin{align}
  |p^{\sigma}C_{1,\size} - C_{p,\size}|
  &\lesssim \size^{-1+\sigma}\int_0^{\size}(1 - e^{-u})u^{-\sigma}\,\intd u\\
  &= \int_0^1(1 - e^{-\size u}) u^{-\sigma}\,\intd u \leq \frac{1}{1 - \sigma}.
\end{align}
Henceforth, $B \lesssim 1$.

\end{document}